\documentclass[11pt,a4paper]{article}
\usepackage[utf8]{inputenc}

\usepackage{comment}
\usepackage{graphicx}

\DeclareUnicodeCharacter{2212}{\textendash}
\usepackage{geometry}                		
                   		
\usepackage{graphicx}				
							\usepackage{enumerate}
\usepackage{amssymb}
\usepackage{tcolorbox}
\usepackage{amsmath}	
\usepackage{changepage}
\usepackage{makecell}
\usepackage{enumitem}
\usepackage{mathtools}
\usepackage{mathtools}
\usepackage{MnSymbol}
\usepackage[numbers]{natbib} 
\usepackage[colorlinks]{hyperref}
\usepackage{setspace}
\usepackage{physics}
\usepackage{amsthm}
\usepackage{slashed}
\usepackage{dsfont}
\usepackage{mathrsfs}
\usepackage{listings}
\usepackage{tikz-cd}

\usepackage{cleveref} 
\usepackage{yfonts}
\usepackage{changepage}
\usepackage{float}
\usepackage{bbm}

\usepackage{subcaption}
\usepackage{shuffle}

\newtheorem{theorem}{Theorem}[section]
\theoremstyle{definition}
\newtheorem{definition}[theorem]{Definition}
\theoremstyle{plain}
\newtheorem{proposition}[theorem]{Proposition}
\newtheorem{corollary}[theorem]{Corollary}
\newtheorem{lemma}[theorem]{Lemma}
\theoremstyle{definition}
 \newtheorem{hyp}[theorem]{Hypothesis}
 \newtheorem{rem}[theorem]{Remark}
\newtheorem{notation}[theorem]{Notation}
\numberwithin{equation}{section}

\usepackage{enumitem}   
\usepackage{ stmaryrd }

\hypersetup{citecolor={red}} 
 \newgeometry{bottom=3cm}

 \newcommand{\one}{\mathbbm{1}}

\renewcommand{\norm}[1]{\left\lVert#1\right\rVert}
\newcommand{\Id}{\operatorname{d}}

\counterwithin{figure}{section}

\usepackage{tikz}
\usepackage{eso-pic}
\usetikzlibrary{calc}
\usetikzlibrary{decorations.pathmorphing}

\title{Trace and Regularized Index of Callias Operators on Hyperbolic Space}
\author{Jilly Kevo \footnote{\textsc{J. Kevo, Institute for applied Mathematics, Rheinische Friedrich-Wilhelms-Universit\"at Bonn, Endenicher Allee 60, 53115 Bonn}, \url{jkevo@uni-bonn.de}}}
\date{\today}

\begin{document}

\maketitle
\begin{abstract}
	Let $d \geq 3$ be odd. 
On hyperbolic space $\mathbb{H}^d$, we
consider Callias-type operators, i.e. 
Dirac-Schrödinger operators of the form $D = ic \nabla + A$, where 
$c \nabla$ is the Dirac operator and $A$  a function with values in the self-adjoint operators, and 
provide an explicit formula for the iterated trace 
$\text{tr}_{L^2(\mathbb{H}^d,H)} \text{tr}_{S}( e^{-tD^*D} - e^{-tDD^*})$, $t>0$, that appears in the formulation 
of the Witten index.
Thereby, the operator $D$ is not assumed to be Fredholm.
	\end{abstract}

\tableofcontents

\section{Introduction}
\subsection{Problem statement and outline}
To study the Witten index for 
Callias-type operators $D = ic \nabla + A$
on hyperbolic space
$\mathbb{H}^d$,
$$  \text{ind}_W(D) := \lim_{t \rightarrow \infty} 
    \text{tr}_{L^2(\mathbb{H}^d,H)} \bigg( \text{tr}_S 
    \left(
    e^{-tD^*D} - e^{-tDD^*}
    \right) \bigg),$$
we compute the 
 double trace  
$\text{tr}_{L^2(\mathbb{H}^d,H)} \text{tr}_{S}( e^{-tD^*D} - e^{-tDD^*})$, $t>0$, under assumptions 
on $D$ that do not force it to be Fredholm.\\ 
Let us start by making the above notions precise.
In what follows, $H$ will always be a complex, separable
Hilbert space. 
The self-adjoint model operator will be denoted by $A_0$,
with domain $\text{Dom}(A_0)$. For $f$ a Borel 
function of $A_0$, its domain $\text{Dom}f(A_0)$
will be considered together 
with the graph norm, as to embed it continuously into $H$. 
\begin{definition}
   Let 
    $A$ be a function on $\mathbb{H}^d$ with 
    values in the self-adjoint operators 
    in $H$ with common domains
    $\text{Dom}(A(x)) =\text{Dom}(A_0)$, 
    $x \in \mathbb{H}^d$.
    Furthermore,  denote by $M_0$ the 
    operator valued function $x \mapsto A_0$.
    By abuse of notation, we write 
\begin{equation}
    \begin{split}
        (Af)(x) &:= A(x)f(x), \quad x \in \mathbb{H}^d, \: f \in \text{Dom}(A), \\
        (M_0f)(x) &:= A_0 f(x), \quad 
        x \in \mathbb{H}^d, \: f \in \text{Dom}(M_0),
    \end{split}
\end{equation}
for the multiplication operators defined
by $M_0$ and $A$, where 
\begin{equation}
    \text{Dom}(A) := \text{Dom}(M_0) := 
    L^2(S, \text{Dom}(A_0)).
\end{equation}
\end{definition}
We postpone the discussion on the spin structure
of $\mathbb{H}^d$ as a homogeneous, symmetric space to section \ref{structureH} and 
also refer to it for a definition of the Dirac 
operator in that context. 
In all what follows,
\begin{equation}
    \text{Dom}(c\nabla) := H^1(S, H),
\end{equation}
will be taken as the domain for the Dirac operator,
 where $c$ denotes Clifford 
 multiplication.

\begin{definition}
    Define 
    \begin{equation}
    D := i c\nabla + A, \quad \text{and} \quad 
    D_0 := ic \nabla + M_0,
    \end{equation}
    where 
    \begin{equation}
        \text{Dom}(D) = \text{Dom}(D_0) = 
        \text{Dom}(c \nabla) \cap \text{Dom}(M_0).
    \end{equation}
Moreover, set 
\begin{equation}
    H_0 := D_0D_0^* = D_0^*D_0 = \slashed{\Delta} + M_0^2,
\end{equation}
with 
\begin{equation}
    \text{Dom}(H_0) = H^2(S,H) \cap \text{Dom}(M_0^2),
\end{equation}
and where $\slashed{\Delta}$ denotes the square of the Dirac operator.
\end{definition}
\begin{notation}
    Let $S^p(H)$ denote the $p$-th 
    Schatten-von Neumann class of $H$ and 
    $B(H)$ its linear, bounded operators.
    Set $\langle T \rangle_z := (T^*T +z)^{1/2}$
    for the Japanese bracket, 
    $\langle T \rangle := \langle T \rangle_1$,
    where $T$ is an operator on $H$ with resolvent
    set $\rho(T)$. For $z \in \rho(T)$, 
    let $R_z(T) := (T+z)^{-1}$ be the resolvent.
\end{notation}

\begin{definition}
    Let $B$ be a function of $\mathbb{H}^d$
    with values in the symmetric operators in
    $H$, with domains $\text{Dom} \: B(x) = \text{Dom} \: A_0$. 
    Let 
    \begin{equation}
        D := D_B = ic \nabla + B + M_0, \quad 
        \text{Dom}(D) := \text{Dom}(D_0),
    \end{equation}
    and 
    \begin{equation}
        H_B := \slashed{\Delta} + (B+M_0)^2, 
        \quad \text{Dom}(H_B)
:= \text{Dom}(H_0).
\end{equation}
\end{definition}
\begin{definition}
Let $S \rightarrow \mathbb{H}^d$ denote the spinor bundle, $x_0 \in \mathbb{H}^d$ the base point and $\tau(x_0,x) \: : \: S_x \rightarrow S_{x_0}$ spinor parallel transport.
To define the  trace $\textnormal{tr}_S$, we 
lift the parallel transport map to 
$U \: : \: L^2(M,S) \rightarrow L^2(M,S)$, 
acting as 
$(U \psi)(x) := \tau(x_0,x) \psi(x)$, 
where $\psi(x) \in S_x$. 
If $T$ is an operator on $L^2(M,S)$, we set $\widetilde T := U T U^{-1} \: : \: L^2(M,S_{x_0}) \rightarrow L^2(M,S_{x_0})$ so that
\begin{equation*}
    \text{tr}_S(T) :=\text{tr}_{S_{x_0}}(\widetilde T).
\end{equation*}
\end{definition}
\begin{rem}
For ease of notation, we suppress the identity 
operators appearing in 
$D = i \sum_{j=1}^d c_j \otimes \partial_{x^j} \otimes \one_H + \one_{\mathbb{C}^r}
\otimes M_A$ 
and $H_B := \one_{\mathbb{C}^r} 
\otimes (\slashed{\Delta} 
\otimes \one_H + M_{A^2})$ whenever there is no 
risk for confusion.
\end{rem}
Let us now mention the hypotheses we need to impose
on the operator valued function $A$.
Here, we will focus on a quantitative description of the
assumptions and delay the precise formulation 
to Section \ref{section2.2}, Hypothesis \ref{hypothesis}. 
\begin{hyp}[Hypothesis \ref{hypothesis}, heuristically]
The assumptions we impose play two conceptual roles:
\begin{enumerate}
    \item to ensure that the perturbed Dirac 
    operator $D = ic \nabla + A$ is well-behaved as 
    an unbounded operator (closedness, 
    common domain, self-adjoint squares), and, inter alia,
    $$\text{tr}_{S} \left( 
    e^{-tD^*D} - e^{-tDD^*}\right) \in S^1(L^2(\mathbb{H}^d,H)). $$
    \item to encode a controlled notion of asymptotic translation invariance of $A$  at infinity on $\mathbb{H}^d$, replacing the usual asymptotic homogeneity condition in the original Euclidean Callias theory.
\end{enumerate}
Let us describe these points in more detail:
\begin{enumerate}[label=1.\arabic*]
\item The first group of conditions requires that $B = A - A_0$ belongs to a sufficiently regular operator class $C_b^1(\mathbb{H}^d, A_0,2N+1)$ and that certain mixed resolvent-weighted semi-norms $\rho_z^{\pm 2N}(B)$ remain small when $z$ is taken far from $(-\infty,0]$.
Heuristically, this says that $A(x)$ varies smoothly in $x$ and all $A(x)$ have constant domain $\text{Dom}(A_0)$.
As a consequence (see Proposition \ref{domain}), 
$$\text{Dom}(DD^*) 
    = \text{Dom}(D^*D) = \text{Dom}(H_0),$$
    and
    $$\text{Dom}(D) = \text{Dom}(D^*) = \text{Dom}(D_0).$$
Moreover,  $D$ is closed with $D^* = -ic\nabla + A$, and  $DD^*$, $D^*D$, and $H_B$ are all non-negative and self-adjoint.
As already remarked, those assumptions are
also crucial in the proof of the 
Schatten-von-Neumann membership of the partial
trace (Proposition \ref{decomposition2}), as well as the limiting procedure
involved in the approximation scheme of the
perturbation (Proposition \ref{approximationscheme}).

\end{enumerate}
\begin{enumerate}[label=2.\arabic*]
\item The second group of assumptions regulates how $A(x)$ behaves when $x$ escapes to infinity in $\mathbb{H}^d$.
In the Euclidean setting, Callias assumed that $A(x)$ becomes asymptotically homogeneous or approaches a fixed limit along rays. Here, we do not assume pointwise convergence of $A(x)$ at infinity, but  require instead that along certain geodesics,
$A(w)$ and $A(w')$, for $w = \exp_{w'}(y)$ and $
w' =  \exp_o(Ry)$ being exponential polar coordinates, become
exponentially close, uniformly in the angular direction $y\in S_1(o)$ and in the point $x \in \mathbb{H}^d$.
The exponential factor $e^{R + \sqrt{R}}$ reflects the geometry,
see Remark \ref{rem2.17}. This  replaces the Euclidean notion of “$A$ stabilizes at infinity”.
\item The bound 
$$\norm{\Id B(x)}_{S^1(H)} \sim \: 
\text{Vol}_{\mathbb{H}^d}(S_R(o))^{-\frac{1}{d-1}},$$
for $d(x,o) = R \rightarrow \infty$,
compensates, for instance, for the 
factor $\sinh^{d-1}(R)$ that is produced when rescaling
certain integrals from $S_R(o)$ to the
geodesic unit sphere $S_1(o)$, see e.g. 
Lemma \ref{lemma6.2} or Lemma \ref{lemma3}.

\item Finally, the conditions in \eqref{hyp4.1}–\eqref{hyp4.3} strengthen the translation invariance properties for $B$ and $\Id B$, requiring that
their difference along certain geodesics escaping to infinity
decays at least exponentially in $R$.
These conditions 
 ensure the vanishing of the
terms in Lemma \ref{lemma6.2} through the symmetrization 
procedure described in Definition \ref{deltadef} 
and Lemma \ref{opproddec}.

\end{enumerate}

\end{hyp}

We recall the notion of partial Witten index as 
stated in \cite{oliver23}, which, in this form, was first
introduced in \cite{gesztesy88}.
\begin{definition}
Consider a linear, densely defined operator $T$ on 
$\mathbb{C}^r \otimes X$, where $X$ is a separable, complex
Hilbert space.
The partial, semi-group regularized Witten index of $T$
is said to exist whenever, for some $t_0 >0$,
\begin{equation*}
    \text{tr}_{\mathbb{C}^r} \left(
    e^{-tT^*T} - e^{-tTT^*}
    \right) \in S^1(X), \quad \text{for all} \: t \geq t_0,
\end{equation*}
and  the limit
\begin{equation*}
    \text{ind}_W(T) := \lim_{t \rightarrow \infty} 
    \text{tr}_X \bigg( \text{tr}_{\mathbb{C}^r} 
    \left(
    e^{-tT^*T} - e^{-tTT^*}
    \right) \bigg),
\end{equation*}
exists. Then $\text{ind}_W(T)$ is called the 
partial Witten index of $T$.

\end{definition}
\begin{rem}
Analogously to the proof of Lemma 1 in
\cite{callias}, one can 
show, that if 
$T$ is Fredholm and there is some $t_0 >0$
such that
\begin{equation*}
    \text{tr}_{\mathbb{C}^r} \left(
    e^{-tT^*T} - e^{-tTT^*}
    \right) \in S^1(X), \quad \text{for all} \: t \geq t_0,
\end{equation*}
then the partial Witten index of $T$ exists
and coincides with its Fredholm index, 
\begin{equation*}
    \text{ind}(T) = \text{dim  ker}(T) - \text{dim  coker}(T).
\end{equation*}

\end{rem}

\begin{rem}
Contrary to the Fredholm index which takes values
in $\mathbb{Z}$, the Witten index may assume any
real number. 
An explicit example illustrating this has been derived in the Euclidean case, see Proposition 1.2. in \cite{oliver25}.
There, the author considers  $(d+1)$-Dirac-Schrödinger  operators of the form
\begin{equation*}
    (D_V F) := i \sum_{j=1}^d c^j \partial_{x^j}
    f(x,y) + i \partial_yf(x,y) + i \sum_{j=1}^d
    c^j F_j(x) \phi(y) f(x,y),
\end{equation*}
where $f \in W^{1,2}(\mathbb{R}^{d+1}, \mathbb{C}^r \otimes H)$, 
$\phi \in C_c^{\infty}(\mathbb{R})$ with $\int_{\mathbb{R}} \phi(y) \Id y = 1$ and
$F \: : \: \mathbb{R}^d
\rightarrow \mathbb{R}^d$ is a bounded smooth function with bounded derivatives, 
satisfying 
\begin{equation*}
    \begin{split}
        &\partial_{x^{i}} F(x) = O( \vert x \vert^{-1}), \quad \vert x \vert \rightarrow \infty, \\
        \exists  \epsilon >0 \: : \: & \partial_R F(x) = O( \vert x \vert^{-1 - \epsilon}),
        \quad
        \vert x \vert \rightarrow \infty \\
        & \lim_{R \rightarrow \infty} 
        R^{\vert \gamma \vert} 
        \vert (\partial^{\gamma} F)(x_0 + Rz) -
        (\partial^{\gamma} F)( Rz) \vert = 0,
    \end{split}
\end{equation*}
for $z \in S_1(0)$, $x_0 \in \mathbb{R}^d$
and $\gamma \in \mathbb{N}_0^d$, with 
$\vert \gamma \vert \leq 1$. 
In this case, the Witten index assumes a particularly simple form, allowing surjectivity to be demonstrated through direct computation.

\end{rem}

\begin{rem}
Let us briefly discuss why the Callias problem is only considered
 in odd-dimensional open space. 
 There are several ways to show this, one of them
 being through the theory of multiple operator integrals (MOI)
 \footnote{I learned this from Oliver 
 Fürst's lecture V5B10 Selected  Topics in Functional Analysis and Operator Theory - \textit{Multiple Operator Integrals and Spectral Shift}, that was held in the winter term 2025/2026 at the University of Bonn.}. For a different approach,
 see for instance the discussion  after Theorem $1.2.$ in \cite{gesztesy16}. \\
Let $\mathbb{E}_X = (E_0, \dots, E_j)$ be a 
tuple of spectral measures over $(X_0, \dots X_j)$,
$H_B := \one_{\mathbb{C}^r} 
\otimes (\slashed{\Delta} 
\otimes \one_H + M_{A^2})$, 
$DD^* = H_B + ic(M_{\partial_{x^j}A}) := H_B + G$
and $D^*D = H_B -G$.
Imposing that $G$ is bounded, for instance, by requiring that $x \mapsto \partial_{x^j}A(x) 
\in L^{\infty}(\mathbb{H}^d,B(H))$, we set
$$T_j(f,H_B,G) := \mathcal{J}(\text{Dif}_j(f), \mathbb{E}_{H_B}, \mathbb{G}),$$
and 
$$R_j(f,H_B,G) := j \int_0^1 (1-t)^{j-1}
\mathcal{J}(\text{Dif}_j (f), \mathbb{E}_{H_B + tG},
G),$$
where $\mathcal{J}(\text{Dif}_j(f), \mathbb{E}_{H_B}, \mathbb{G})$ is a multiple 
operator integral 
and $\text{Dif}_j(f)$
denotes the higher divided differences of $f$, see e.g. \cite{peller}
for definitions and
an introduction to the topic of multiple 
operator integrals. Now, for 
$f \in W^d(\mathbb{R})$, where $W^d$ 
denotes the $d$-th Krein-Wiener class, we have
\begin{equation*}
    \begin{split}
        &f(D^*D) - f(DD^*) \\
        &=
      f(D^*D) - f(H_B) + f(H_B) - f(DD^*)   \\
      &= \sum_{j=0}^{d-1} \left[
      T_j(f,H_B,-G ) - T_j(f,H_B,G) 
      \right] + R_d(f,H_B,-G) 
      - R_d(f,H_B,G),
    \end{split}
\end{equation*}
where the second inequality follows from the non-commutative Taylor Theorem.
By multilinearity of the MOI,
\begin{equation*}
    T_j(f,H_B, \pm G ) = (\pm1)^j 
    \mathcal{J}(\text{Dif}_j(f), \mathbb{E}_{H_B},
    \mathbb{G}),
\end{equation*}
so
\begin{equation*}
\begin{split}
    &f(D^*D) - f(DD^*)\\
    &= 
   \sum_{j=0}^{d-1}  ((-1)^j-1)
   \mathcal{J}(\text{Dif}_j(f), \mathbb{E}_{H_B},
    \mathbb{G}))
   + R_d(f,H_B,-G) 
      - R_d(f,H_B,G) \\
      &= - 2 \sum_{\substack{j \: \text{odd}\\ 1 \leq j < d}} 
      \mathcal{J}(\text{Dif}_j(f), \mathbb{E}_{H_B},
    \mathbb{G}))
      + R_d(f,H_B,-G) 
      - R_d(f,H_B,G).
      \end{split}
\end{equation*}
Note that since $H_B = \one_{\mathbb{C}^r}
\otimes (\slashed{\Delta} \otimes \one_H + M_{A^2})$, it holds that
$\mathbb{E}_{H_B} = \one_{\mathbb{C}^r} \otimes 
\mathbb{E}_{\slashed{\Delta}  \otimes \one_H + M_{A^2}}$
and so $\mathbb{E}_{H_B}$ commutes with 
$c_j \otimes \one_{L^2(\mathbb{H}^d,H)}$. 
Given that $G = i \sum_{j=1}^{d} c_j \otimes 
M_{\partial_{x^j}A}$, 
multilinearity and commutativity of the 
MOI imply
\begin{equation*}
\begin{split}
    &\mathcal{J}(\text{Dif}_j(f), \mathbb{E}_{H_B},
    \mathbb{G}) \\
    &= \mathcal{J}(\text{Dif}_j(f), \one_{\mathbb{C}^r} \otimes 
\mathbb{E}_{\slashed{\Delta}   \otimes \one_H + M_{A^2}} ,
    \mathbb{G})        \\
    &= i^j \sum_{\lambda \in \{1, \dots,d \}^j} 
    c_{\lambda_1}  \dots 
    c_{\lambda_j} \otimes 
    \mathcal{J}(\text{Dif}_j(f), 
    \mathbb{E}_{\slashed{\Delta}  \otimes \one_H + M_{A^2}},
    (M_{\partial_{x^{\lambda_1}}A}, \dots, 
    M_{\partial_{x^{\lambda_j}}A})).
    \end{split}
\end{equation*}
Applying the trace $\text{tr}_{\mathbb{C}^r}$
yields then 
\begin{equation*}
\begin{split}
    &\text{tr}_{\mathbb{C}^r} \mathcal{J}(\text{Dif}_j(f), \mathbb{E}_{H_B},
    \mathbb{G}) \\
    &= i^j 
    \sum_{\lambda \in \{1, \dots,d \}^j} 
    \text{tr}_{\mathbb{C}^r}(c_{\lambda_1}  \dots 
    c_{\lambda_j})
    \cdot  \mathcal{J}(\text{Dif}_j(f), 
    \mathbb{E}_{\slashed{\Delta}  \otimes \one_H + M_{A^2}},
    (M_{\partial_{x^{\lambda_1}}A}, \dots, 
    M_{\partial_{x^{\lambda_j}}A})), 
    \end{split}
\end{equation*}
using the fact that 
$\text{tr}_{\mathbb{C}^r}(S_1 \otimes S_2) =
\text{tr}_{\mathbb{C}^r}(S_1) \cdot S_2$,
for $S_1 \in B(\mathbb{C}^r)$ 
and $S_2 \in B(L^2(\mathbb{H}^d,H))$.
Now, provided that $d$ is even and $f \in W^{\infty}(\mathbb{R})$, such that
$$\sum_{n=0}^{\infty} \frac{\norm{f}_{W^n}}{n!} < \infty,$$
as is the case for
$f(z) = e^{-tz}$ and any $t>0$, 
\begin{equation*}
     f(D^*D) - f(DD^*) = 
     -2 \sum_{j \: \text{odd}}^{\infty} 
     \mathcal{J}(\text{Dif}_j(f), \mathbb{E}_{H_B},
    \mathbb{G}),
\end{equation*}
converges
and given that in even dimensions, 
\begin{equation*}
    \text{tr}_{\mathbb{C}^r}(c_{\lambda_1}  \dots 
    c_{\lambda_j}) =0,
\end{equation*}
for all odd $j$ by Lemma \ref{cliffordid}, 
we thus have
\begin{equation*}
     \text{tr}_{\mathbb{C}^r}( f(D^*D) - f(DD^*))=0.
\end{equation*}

\end{rem}

\begin{rem}
    On the other hand, note that every 
    elliptic operator on a compact manifold is 
    Fredholm and that the index of any such  
    operator on an odd-dimensional, compact manifold
    is automatically zero.
\end{rem}

Let us now state the main theorem of this paper.

\begin{theorem}[Theorem \ref{maintheorem}]
 Let $t >0$ and assume Hypothesis \ref{hypothesis}.
For real-valued $\phi \in C^{\infty}_c(\mathbb{H}^d)$, 
write $A_{\phi} := A_0 + (1- \phi)(A-A_0)$
and set $ \kappa_c := (2i)^{\frac{d-1}{2}} (-i)^d$.
Then 
\begin{equation}
    \begin{split}
&\textnormal{tr}_{L^2(\mathbb{H}^d,H)}
    \textnormal{tr}_{S} \left(
    e^{-tD^*D} - e^{-tDD^*}\right) \\
    &=  \frac{2}{d} i^d \kappa_c 
     \left[ t^d p_t(0) \mathcal{X}_t -
    \int_0^t s^{d/2} \partial_s \left(s^{d/2} p_s(0)\right) \mathcal{X}_s \Id s  
     \right],
    \end{split}
\end{equation}
 where 
 \begin{equation*}
 \mathcal{X}_s = \int_{\mathbb{H}^d} X_s, \qquad      X_s := \int_{u \in \Delta_{d-1}} 
    \textnormal{tr}_H \bigwedge_{j=1}^d 
  \bigg( (\Id A_{\phi})e^{-su_j A_{\phi}^2}\bigg) \Id u,
 \end{equation*}
and 
\begin{equation*}
        \begin{split}
            p_t(0)
=
&\frac{1}{2^{d} \pi^{d/2}\Gamma(d/2)} \\
& \cdot\sum_{k=0}^{m}
\Gamma\left(k+\frac12\right)
e_{m-k}\left(
\frac14,\frac94,\ldots,\frac{(2m-1)^2}{4}
\right)
t^{-k-\frac12},
        \end{split}
    \end{equation*}
    with $e_j$ denoting the $j$-th elementary symmetric polynomial,
with $e_0=1$.

\end{theorem}

In Section \ref{section2}, we establish the main
 analytical framework, that is, fix notations concerning
 the base manifold $\mathbb{H}^d$, the spinor bundle,
 and the Dirac-Schrödinger operator $D = ic \nabla + A$. We also provide
 a brief introduction to the 
 Lie-theoretic analysis of $\mathbb{H}^d$ as 
 a symmetric space, which allows us to introduce 
 and work with the Fourier-Helgason transform
 in the next section.
 A central component is Hypothesis \ref{hypothesis},
 which imposes analytic and asymptotic constraints on 
$A$
and its derivatives, ensuring that $A$ has, 
among others, the right domain, as well as  trace-class properties needed for later sections.
 \\
Section \ref{section3} demonstrates trace-class 
membership of the partial trace, that is
$$\textnormal{tr}_{S}\left(e^{-tD^*D} -
        e^{-tDD^*} \right) \in S^1(L^2(\mathbb{H}^d,H)).$$
        This is achieved by showing the 
         corresponding statement for the associated resolvent
    identity, which is in turn  accomplished by decomposing the latter into a sum of bounded and Schatten-von-Neumann operators. Afterwards, 
    holomorphic functional calculus is used to 
deduce the exponentiated version.\\
Section \ref{section4} develops an approximation scheme for the operator-valued function $B$, 
through spectral projections, finite-rank operators and a mollification procedure. 
This will allow us to prove the final trace formula
in subsequent sections.\\
Building on the approximation scheme, Section 
\ref{section5}  constructs integral kernels for the 
main operators in order to apply a
Mercer-type theorem for the calculation of the iterated 
trace and also establishes associated decay and regularity 
estimates.\\
In Section \ref{section6}, we combine all 
previously established results to
prove 
Theorem \ref{maintheorem}.

\subsection{Brief historical overview of the problem}
Besides our work, which extends the Euclidean case treated by Fürst in \cite{oliver23} 
to $d$-dimensional hyperbolic space, the study of Callias operators has a much longer history. For an account of the mathematical development of the theory, beginning with the celebrated Atiyah–Patodi–Singer “index = spectral flow” result \cite{atiyah1, atiyah2, atiyah3}, we refer the reader to the introductory section of \cite{oliver23}. There, the evolution of the assumptions imposed on the perturbation is also carefully traced through the literature.
To avoid redundancy, we summarize only the key mathematical milestones in Tables $1$ to $4$
below and instead focus our discussion on a complementary aspect, namely the physical lineage of the theory.

The study of the Callias problem predates Callias’s classical paper \cite{callias}.  
Its origins can be traced to the 
Georgi-Glashow model, an $SU(2)$ 
Yang-Mills-Higgs theory coupled to a scalar field
in the adjoint (triplet) representation.
In this model, 't Hooft \cite{hooft} 
and Polyakov \cite{polyakov} independently
discovered smooth, finite-energy magnetic monopole
solutions, that is, 
static, spherically symmetric classical solutions of the field equations carrying non-vanishing magnetic charge.
The Lagrangian density of the Georgi-Glashow model
is
 \begin{equation*}
     \mathcal{L}
  = \frac{1}{2} (D^\mu \Phi)^2
  -\frac{1}{4g^2}(F^{\mu \nu})^2 - \frac{1}{4}
  \lambda (\Phi^2 - v^2)^2,
     \end{equation*}
 where the covariant derivative and field strength are defined by
 \begin{equation*}
     \begin{split}
         D^\mu \Phi &= \partial^\mu \Phi + A^\mu \times \Phi, \\
         F^{\mu \nu} &= \partial^\mu A^\nu
         - \partial^\nu A^\mu + A^\mu \times A^\nu.
     \end{split}
 \end{equation*}
 In the vacuum, that is $\Phi^2 = v^2$, the 
 $SU(2)$ gauge symmetry is broken down to 
 electromagnetic $U(1)$.
Accordingly, the vacuum manifold is the quotient
$$SU(2)/U(1) \simeq S^2,$$
and since finite-energy considerations require 
that the Higgs-field approaches a vacuum configuration at spatial infinity, 
the normalized Higgs field at infinity defines 
a continuous map $\Phi_{\infty} \: : \: S^2_{\infty} \rightarrow S^2$, for 
a monopole solution 
$\Phi \: : \: \mathbb{R}^3 \to \mathfrak{su}(2)$.
Such maps, as well as 
 the superselection sectors,
 are classified 
 by 
\begin{equation*}
    \pi_2(SU(2)/U(1)) = \mathbb{Z} \neq \{1 \},
\end{equation*}
guaranteeing that monopole solutions exist (which, in practice, can be found using the 
so-called "Hedgehog" Ansatz) and
where the corresponding integer is identified
with the magnetic charge of the monopole.
  \\
Subsequent work by Jackiw and Rebbi \cite{jackiw} showed in $1976$ that a Dirac fermion coupled to the  
't Hooft-Polyakov monopole
possesses a normalizable zero-energy solution. 
Each  monopole of unit charge supports exactly one such
zero mode, which is localized around the monopole's 
core,  i.e. the region where the Higgs field is small
or vanishes, and the fermion becomes massless. 
Conceptually speaking, the core acts like a potential
well that binds the fermionic zero mode, 
preventing it from spreading to infinity, which, in
 the Callias framework, 
is essential to guarantee Fredholmness of the 
involved operators.
This was  a profound discovery, since the existence of
 zero modes implies fractional fermion numbers. \\
Finally, a mathematical explanation of the
Jackiw–Rebbi zero mode 
was provided by Callias \cite{callias} in 1978, who considered Dirac-type operators of the form
\[
D_\Phi = \gamma^i(\partial_i + A_i) + \Phi(x),
\] 
where $\Phi(x)$ is a Higgs-like potential 
and is given, like $A_i$, by the monopole 
configuration. 
Callias assumed that \(\Phi(x)\) satisfies what is today often called
"Callias condition", namely, that \(\Phi(x)^2\) is in some sense uniformly invertible outside a compact region, so that under this hypothesis, \(D_\Phi\) is a Fredholm operator, with Fredholm index depending only on the asymptotic behavior of \(\Phi\). This corresponds to localization 
near the monopole core.
In three dimensions, Callias proved the topological formula
\[
\operatorname{ind}(D_\Phi) = \deg(\Phi_\infty),
\]
where \(\Phi_\infty : S^2_\infty \to S^2\) is again the normalized Higgs field at infinity. 
He concluded that, for a monopole of magnetic charge \(n\), the Dirac operator possesses precisely \(n\) chiral zero modes, with the explicit zero mode found by Jackiw and Rebbi corresponding to the case \(n=1\), hence
interpreting the existence of the zero modes
as a topological property.
Mathematically speaking, the novelty was that Fredholmness came from the Higgs field at infinity, not compactness of the manifold and it was the very first index theorem of that kind.\\
Weinberg \cite{weinberg} extended the analysis to multi-monopole 
configurations, showing that a magnetic monopole
of magnetic charge $n$ supports exactly $n$ chiral
zero modes, thereby confirming Callias's abstract theorem. \\
As for the Witten index itself, its origins can be traced back to the following formula for the Fredholm index of a Fredholm operator $T$ on a Hilbert space $H$,
\begin{equation}
\begin{split}
\operatorname{ind}(T)
&= \dim\ker T^*T-\dim\ker TT^* \\
&= \operatorname{Tr}_H\bigl(f(T^*T)-f(TT^*)\bigr),
\end{split}
\end{equation}
where $f\in C_c^\infty(-\varepsilon,\varepsilon)$
 satisfies $f(0)=1$, for sufficiently small 
 $\varepsilon>0$.
 The particular choice of the function $f$
  is irrelevant, since by functional calculus $f(T^*T)$ and $f(TT^*)$ are precisely the orthogonal projections onto $\text{ker}T$ and $\text{ker}T^*$, respectively. 
  A particularly important choice is given by the heat kernel. In the setting of a $\mathbb{Z}_2$-graded Dirac-type operator $D$, the McKean--Singer formula \cite{mckean} asserts that
  $$\operatorname{ind}(D)
=
\operatorname{Tr}_H
\left(
e^{-tD^*D}-e^{-tDD^*}
\right),
\qquad t>0$$ 
Since, in the strong operator topology,
$$e^{-tD^*D}\longrightarrow P_{\ker D},
\qquad
e^{-tDD^*}\longrightarrow P_{\ker D^*}
\qquad\text{as }t\to\infty,$$
this heat-kernel expression recovers the usual definition of the Fredholm index in the large-time limit. 
The Witten index arises by retaining this heat-kernel expression in situations where the operator is not necessarily Fredholm. One defines, whenever the limit exists,
$$\operatorname{Wind}(T)
:=
\lim_{t\to\infty}
\operatorname{Tr}_H
\left(
e^{-tT^*T}-e^{-tTT^*}
\right)$$
In contrast to the Fredholm case, the heat trace difference need not be independent of $t$, and the existence of the large-time limit is therefore a nontrivial question. 
In this sense, the Witten index may be viewed as an extension of the McKean--Singer heat-kernel formula to operators for which the ordinary Fredholm index is no longer defined. \\
However, the necessity  of introducing
a Witten index-like object is much older, dating back to Laplace and his work on celestial mechanics, 
namely his  study of planetary motions.
More precisely, given a nowhere-vanishing trigonometric polynomial $f(t)=\sum_{k=1}^n a_k e^{i\lambda_k t},
\quad a_k\in\mathbb C,\quad \lambda_k\in\mathbb R$, 
one may choose a continuous argument function $\operatorname{Arg} f(t)$ and ask whether the limit, the so-called mean motion of $f$,
$\lim_{t\to\infty}\frac{\operatorname{Arg}f(t)}{t}$ exists. 
Laplace established the result in the special case of $n=3$, while Weyl subsequently proved the conjecture using his then-recently proved ergodic theorem. 
The resulting formula, 
given by
$$\mu(f) := \frac{1}{2\pi i T}\int_{0}^T\frac{f'(s)}{f(s)}\,\Id s,$$ bears a close resemblance to the classical formula for the winding number about $0$.
This connection is further reflected in the Toeplitz index theorem, which states that any invertible, continuous function on the circle $\varphi$ is Fredholm with index given by 
$$-\operatorname{wind}(\varphi)
=-\frac{1}{2\pi i}\int_{\mathbb T}\frac{\varphi'(z)}{\varphi(z)} \Id z.$$
Now, upon closer inspection, one notices 
 that the mean motion, just like the winding number, is not only stable under small perturbations, but that the condition 
 $f(t) \neq 0$ also resembles
ellipticity.
 This analogy suggests, that one may 
 try looking for an interpretation of the mean motion as an index of some kind. \\
For instance, such an interpretation arises naturally in the theory of Wiener--Hopf operators, where, to an almost-periodic function $f$ one can associate a Wiener--Hopf operator $W_f$ acting on $L^2(\mathbb R_+)$. These operators belong to a type $\mathrm{II}_\infty$ von Neumann algebra equipped with a semifinite trace $\tau$. 
Theorem 2.2 of Coburn, Douglas, Schaeffer and Singer \cite{coburn} shows that, for a boundedly invertible differentiable almost-periodic function $f$, the operator $W_f$ is $\tau$-Fredholm and satisfies
$\operatorname{Ind}\tau(W_f)=-\mu(f)$,
where $\mu(f)$ denotes the mean motion. Thus mean motion admits an interpretation as a Breuer--Fredholm index in a semifinite von Neumann algebra.
This perspective was further developed by Carey and Pincus \cite{carey77} in their study of almost-normal operators. For an almost-normal operator $T$ in a semifinite von Neumann algebra, Carey and Pincus introduced a principal function $g_T\in L^1(\mathbb C)$ characterized by trace identities of Helton--Howe type. They conjectured that, for a Wiener--Hopf operator $W_f$ associated with a finite trigonometric sum,
$g_{W_f}(z)=-\mu(f-z)$
for almost every $z\in\mathbb C$, and proved 
it in \cite{carey79}.
This should be compared with the formulation 
of the Witten index in terms of the
(Krein-Lifshitz or higher order) spectral shift function, as is sketched in Table 
\ref{table3}. \\
Note that a notion of the Witten index  also appeared in physics,
where it was  defined in 
$1981$ by Witten \cite{witten81}, who introduced supersymmetric quantum mechanics (SUSY QM) as a toy model for studying supersymmetry breaking. 
In its original form, it was first stated as 
\begin{equation*}
    W = \text{Tr} (-1)^F e^{-\beta H},
\end{equation*}
where $H = Q^2$ is the SUSY Hamiltonian,
$Q$  a supercharge and $F$ the fermion parity operator.
Conceptually, the Witten index counts 
the difference between bosonic and fermionic zero-energy states
and, being invariant under smooth deformations of the Hamiltonian, also constitutes a topological invariant.
In the SUSY QM realization of a Dirac operator, 
$H = D^2$, $Q = D$, the Witten index 
equals the analytic index of $D$. When $D$ is a 
Callias operator, the Witten index computed from 
SUSY QM reproduces the Callias index, extending this 
correspondence to a non-compact setting. \\
Later, in $1987$, Bollé et al \cite{bolle}
related the Witten index and axial anomaly to Krein’s spectral shift function, further extending the index to continuous spectra and non-Fredholm operators. \\ \\
Below, we summarize the index formulas for Euclidean Callias operators in the Fredholm case.
We emphasize once more that, in order to ensure Fredholmness, the operator-valued family 
$A$ must satisfy a uniform invertibility condition outside a compact subset.
\begin{table}[H]\label{table1}
\centering
\begin{tabular}{|l|l|l|}
\hline
      & \makecell{$\mathbb{R}$}     & \makecell{$\mathbb{R}^d$, $d \geq 3$ odd}  \\ \hline
\textbf{Index} & \makecell{ $\text{ind}(D) =$ \\
$\text{SpecFlow}(A,0)$} & \makecell{$\text{ind}(D) =$ \\ $c_d \lim_{R \rightarrow \infty} \int_{S_R} \text{tr}_{\mathbb{C}^m}\left( U (\text{d}U)^{\wedge(d-1)}\right)$, \\
$U := \vert A \vert^{-1} A$} \\ \hline
\end{tabular}
\caption{Fredholm case.}
\end{table}
\begin{table}[H]\label{table2}
\centering
\begin{tabular}{|l|l|l|}
\hline
      & \makecell{$\mathbb{R}$}     & \makecell{$\mathbb{R}^d$, $d \geq 3$ odd}  \\ \hline
\textbf{Index} & \cite{atiyah1, atiyah2, atiyah3} & \makecell{ \cite{callias}, \cite{gesztesy16} } \\ \hline
\end{tabular}
\caption{Literature for the Fredholm case.}
\end{table}
Likewise, known formulas which are also applicable to non-Fredholm Euclidean Callias operators are collected below. \\ \\
\begin{table}[H]
\centering
\resizebox{\textwidth}{!}{%
\begin{tabular}{|l|l|l|}
\hline
 & \makecell{$\mathbb{R}$} & \makecell{$\mathbb{R}^d$, $d \ge 3$ odd} \\ \hline

\makecell{\textbf{Index}} &
\makecell{
$\text{ind}_W(D) = \tfrac{1}{2}\xi(0+,A_+,A_-) + \tfrac{1}{2}\xi(0-,A_+,A_-)$,\\
$A_{\pm} = \lim_{x \to \pm \infty} A(x)$
} &
\makecell{
$\text{ind}_W(D) = - \xi^{(d-1)}_{d,A_0,B}(0+)$
} \\ \hline

\makecell{\textbf{Principal}\\\textbf{Trace}\\\textbf{Formula}} &
\makecell{
$\text{tr}_{L^2(\mathbb{R},H)}(e^{-tD^*D}-e^{-tDD^*}) =$\\
$\sqrt{\tfrac{t}{\pi}}\int_0^1 \text{tr}_H \big( e^{-tA^2(s)}(A_+ - A_-) \big)\, \Id s$
} &
\makecell{
$\text{tr}_{L^2(\mathbb{R}^d,H)}\text{tr}_{\mathbb{C}^r}(e^{-tD^*D}-e^{-tDD^*}) =$\\
$C_d t^{d/2}\!\int_{s \in \Delta_{d-1}}\!\!\int_{\mathbb{R}^d}\!\!\text{tr}_H \bigwedge_{j=0}^{d-1}$\\
$\big(\text{d}A_{\phi} e^{-ts_jA_{\phi}^2}\big)\, \Id s$,\\
$A_{\phi} := A_0 + (1-\phi)(A-A_0), \; \phi \in C_c^{\infty}(\mathbb{R}^d)$
} \\ \hline

\makecell{\textbf{Example}} &
\makecell{
$D = \partial_x + i \partial_y + V(x,y)$,\\
$\text{ind}_W(D) = \frac{1}{2\pi}\!\int_{\mathbb{R}}\!(V(+\infty,y)-V(-\infty,y))\, \Id y$
} &
\makecell{
$D_V := i \sum_{j=1}^d c_j \partial_{x^j} + i \partial_y + V(x,y)$,\\
$\text{ind}_W(D_V) = c_d \int_{\mathbb{R}^d}\text{tr}_H (U^{-1}\text{d}U)^{\wedge d}$,\\
$U(x) := U^{V(x,.)}(\infty,-\infty)$ is the evolution system of\\
$u'(y)=iV(x,y)u(y), \; y\in\mathbb{R}$
}
\\ \hline
\end{tabular}
}
\caption{Beyond the Fredholm case.}
\label{table3}
\end{table} 
\begin{table}[H]\label{table4}
\centering
\begin{tabular}{|l|l|l|}
\hline
 & \makecell{$\mathbb{R}$} & \makecell{$\mathbb{R}^d$, $d \ge 3$ odd} \\ \hline

\makecell{\textbf{Index}} &
\makecell{
\cite{pushni}, \cite{gesztesy11}, \cite{carey14}, ...
} &
\makecell{
\cite{olivershift}
} \\ \hline

\makecell{\textbf{Principal}\\\textbf{Trace}\\\textbf{Formula}} &
\makecell{
\cite{callias}, \cite{pushni}, \cite{gesztesy11},  \cite{gesztesy16}, ...
} &
\makecell{\cite{oliver23}} \\ \hline

\makecell{\textbf{Example}} &
\makecell{
\cite{bolle}, \cite{carey16}, \cite{carey162}
} &
\makecell{
\cite{oliver25}
}
\\ \hline
\end{tabular}
\caption{Literature for the generalized  case.}
\end{table} 
Let us briefly mention that index theory for Callias-type operators has also been investigated on more general non-compact Riemannian manifolds. 
In chronological order, see, for instance, 
\cite{anghel1}, \cite{bruning}, \cite{anghel2}, 
\cite{rade}, \cite{bunke95}, \cite{kucero}, \cite{carvalho}, \cite{kottke}, \cite{braverman},
\cite{shi}, \cite{dungen}, \cite{cecchini}, 
\cite{schulz} and \cite{dungen2}.
However, all of these works assume that the operators
in question are Fredholm. In contrast, Theorem \ref{maintheorem}
provides the first trace formula applicable to the non-Fredholm case in a non-compact, curved space. 

\section{Notations and setup}\label{section2}
\subsection{Some representation theory}\label{structureH}
Hyperbolic space has been well studied and there are 
many sources treating what is summarized below.
Here, we mostly follow the exposition provided
in \cite{camporesi}, but also refer to \cite{helgason2}, 
\cite{helgason}. \\
Let $d \geq 3$ be an odd integer. 
By $d$-dimensional real hyperbolic space, 
we mean the manifold
$$\mathbb{H}^d := \{x \in \mathbb{R}^{d+1} 
\: \vert \: (x,x)_L = -1, \: x_{d+1} > 0 \},$$
with base point $o := (0, \dots, 0,1) \in \mathbb{R}^{d+1}$,
where $(.,.)_L$ denotes the 
Lorentz form 
$$(x,y)_L = x_1y_1 + \dots + x_ny_n - x_{n+1}y_{n+1}.$$
The group of isometries is 
the identity connected component 
 $G' := \text{SO}_e(d,1)$ of the Lorentz group,
 with maximal compact subgroup
 $K' :=SO(d)$. The latter fixes
 the base point $o$ so that the 
 orbit map gives rise to a 
 diffeomorphism $G'/K' \rightarrow \mathbb{H}^d, g'K' \mapsto g \cdot o$.
Using the classical two-fold coverings of $G'$ and $K'$, 
we have
$\mathbb{H}^d := G/K$, where $G := \text{Spin}_e(d,1)$ is the identity
connected component 
of $\text{Spin}(d,1)$, and $K := \text{Spin}(d)$.
Their respective 
Lie algebras will be denoted by 
$\mathfrak{g} = \mathfrak{spin}(d,1)
\simeq \mathfrak{so}(d,1)$ and 
$\mathfrak{k} = \mathfrak{spin}(d) \simeq 
 \mathfrak{so}(d)$. 
 Recall that $\mathfrak{g}$ admits the Cartan 
 decomposition 
 \begin{equation}
     \mathfrak{g} = \mathfrak{k} \oplus \mathfrak{p},
 \end{equation}
where $\mathfrak{p} \simeq \mathfrak{g}/\mathfrak{k} \simeq T_o(G/K)$ 
can be identified with $\mathbb{R}^d$ via the isomorphism 
$$\mathfrak{p} \rightarrow \mathbb{R}^d, 
\quad X := \begin{pmatrix}
0_d & x \\
x^t & 0_d 
\end{pmatrix} \mapsto x.$$
In that case, 
$$g'(X,Y) := y^tx$$
induces by translations a $G$-invariant 
Riemannian metric $g$ on $\mathbb{H}^d$.

\begin{definition}
    Let $\tau$ be a unitary, finite-dimensional
    representation of $K$ on a Hilbert space
    $V_{\tau}$. A function of 
    (right) type $\tau$ on $G$ is a function 
    $f \: : \: G \rightarrow V_{\tau}$ 
    which is $K$-equivariant in the 
    sense of 
    \begin{equation}
    f(x) = \tau(k)^{-1} f(x), \quad 
    x \in G, k \in K.
\end{equation}
\end{definition}

\begin{notation}
     The space of all functions 
    of type $\tau$ will be denoted by 
    $\Gamma(G, \tau)$, with $\Gamma$ replaced 
    by $C^{\infty}$, $L^2$, $C_c$ and so forth,  depending 
    on the function space in need.
\end{notation}
By the 
general theory of homogeneous vector 
bundles on symmetric spaces (see e.g. \cite{wallach} or \cite{camporesi}), the 
sections 
of the homogeneous vector bundle 
$G \times_K V_{\tau}$ correspond to 
functions of type $\tau$ on $G$, i.e.
\begin{equation}\label{sections} 
    \Gamma(G  \times_K V_{\tau}) = \Gamma(G, \tau).
\end{equation}
Denoting by $e_1$ the
element 
\begin{equation}
    e_1 := \begin{pmatrix}
        0 & 0 & 1 \\
        0 & 0_{d-1} & 0 \\
        1 & 0 & 0
    \end{pmatrix} \in \mathfrak{p},
\end{equation}
the algebra $\mathfrak{a} := \mathbb{R}e_1 
\subset \mathfrak{p}$ is a
Cartan subspace, whose corresponding Lie group 
$A$ is a Lie subgroup of $G$, parametrized  
by $a_t := \exp(te_1)$. If 
$\alpha \in \mathfrak{a}^*$ is such that
$\alpha(te_1) = t$, 
then $(\mathfrak{g}, \mathfrak{a})$ admits 
$R(\mathfrak{g}, \mathfrak{a}) = \{\pm \alpha \}$
as restricted root system. 
The positive subsystem is then 
$R^+(\mathfrak{g}, \mathfrak{a}) = \{ \alpha \}$, 
with Weyl group $W(\mathfrak{g}, \mathfrak{a}) =
\{ \pm \one \}$. Note that via 
\begin{equation}\label{system}
\mathfrak{a}^*_{\mathbb{C}} \ni \lambda \alpha \mapsto \lambda \in \mathbb{C},
\end{equation}
one has $\mathfrak{a}^*_{\mathbb{C}} \simeq \mathbb{C}$.
Denote by $\mathfrak{n} = \mathfrak{g}_{\alpha}$
the positive root subspace with analytic subgroup 
$N$.
The half-sum of roots in 
$R^+(\mathfrak{g}, \mathfrak{a})$ is then 
\begin{equation}
    \mathcal{\rho} = \frac{d-1}{2}\alpha,
\end{equation}
which will be 
interepreted as an element in $\mathbb{C}$ via 
\eqref{system}. \\
Recall moreover the 
Iwasawa and Cartan decompositions
of $G$, given, respectively, by 
\begin{equation}
    G = KAN
\end{equation}
and \begin{equation}
    G = K \{a_t \: \vert \: t \geq 0\}K.
\end{equation}
The minimal parabolic subgroup of $G$ is 
\begin{equation}
    P := MAN,
\end{equation}
where $M \simeq \text{Spin}(d-1)$ is the centralizer
of $A$ in $K$.

\begin{definition}
    Let $\tau$ be a finite dimensional 
    representation of $K$. A function $F \: : \: 
    G \longrightarrow \text{End} (V_{\tau})$ 
    is a $\tau$-radial function on $G$ if 
    \begin{equation}
        F(k_1 xk_2) = \tau(k_2)^{-1} F(x) \tau(k_1),
    \end{equation}
    where $x \in G$ and $k_1, k_2 \in K$. 
    The set of all $\tau$-radial functions on $G$
    is denoted by $\Gamma(G, \tau, \tau)$. 
\end{definition}
An example of a  $\tau$-radial function is the 
hyperbolic spinor heat kernel.
\begin{rem}
    As stated in \cite{camporesi}, the $\tau$-radial
    functions are related to the functions of 
    type $\tau$ as follows:
    \begin{itemize}
        \item if $F$ is a $\tau$-radial function, 
        then, for all $v \in V_{\tau}$, the 
        function $x \mapsto F(x)v$ is of 
        type $\tau$;
        \item  any $\tau$-radial function
        is determined by a function of type $\tau$ of the
        form $x \mapsto F(x) v$, where $v \neq 0$.
    \end{itemize}
    Moreover, any $\tau$-radial function is completely
    determined by its restriction to 
    $\{a_t \: \vert \: t \geq 0\}$. 
\end{rem}
With respect to the 
Cartan decomposition, the 
Haar measure on $G$ reads
\begin{equation}\label{eq:haarmeasure}
    \int_G f(x) dx = \int_K \int_0^{\infty} 
    (2 \sinh{t})^{d-1} \int_K f(k_1 a_t k_2) \: \Id k_2 \: \Id t \: \Id k_1,
\end{equation}
where $f$ is any function on $G$, $\Id t$
the standard Lebesgue measure on $\mathbb{R}$
and $\Id k$ is the Haar measure on 
$K$ normalized via $\int_K \Id k =1$.
The inner product on $L^2(G, \tau)$ is then
given by 
\begin{equation}
    (f,g) = \int_G (f_1(x), f_2(x))_{V_{\tau}} \Id x 
\end{equation}
and on $C_c^{\infty}(G, \tau, \tau)$
by
\begin{equation}
    (F,H) = \int_G (F(x), G(x))_{HS}  \Id x,
\end{equation}
where $(.,.)_{HS} := \text{tr}(F(x) H(x)^*)$ denotes the Hilbert-Schmidt 
inner product.
\begin{rem}\label{rem:haarnormalization}
    One has to be careful with the normalization employed by Camporesi-Pedon 
   \cite{camporesi} upon which the above discussion of hyperbolic space as a homogeneous symmetric space $G/K$ is based.
   Namely, the quotient measure 
   induced by the Haar measure normalization
   in Equation \eqref{eq:haarmeasure}
   is 
   \begin{equation}
       \Id \mu_{CP} = \frac{2^{d-1}}{\text{vol}(S^{d-1})} \Id \text{vol}_{\mathbb{H}^d}.
   \end{equation}
   Consequently, the hyperbolic spinor heat kernel, as described in Proposition 
   \ref{prop:spinorheatkernel}, must also be
   rescaled by  $\frac{2^{d-1}}{\text{vol}(S^{d-1})}$ in order to obtain the usual small-time asymptotics. 
   See Remark \ref{rem:smalltimespinorkernel}.
\end{rem}

\subsection{The 
perturbed Dirac operator}\label{section2.2}

Since $\mathbb{H}^d$ is a noncompact 
Riemannian symmetric space, it is
contractible and hence admits a unique spin
structure. 
The spinor bundle $S$ is then 
given by the homogeneous vector 
bundle 
\begin{equation}
S := G \times_K V_{\tau_d},
\end{equation}
where $\tau_d$ is the complex spin representation
of $K$ on $V_{\tau_d} = \mathbb{C}^{2^{[ d/2]} }$, $[d/2]$ standing for the integer part
of $d/2$. Since $d$ is odd, $\tau_d$ is 
irreducible. 
From \eqref{sections},  the spinors can be 
described 
explicitly as 
\begin{equation}
    \Gamma(G \times_K V_{\tau_d}) = \Gamma(G, \tau_d).
\end{equation}
Let now $x \in \mathbb{H}^d = G/K$. Denote by 
\begin{equation}
    c_x \: : \: T^*_x(G/K) \longrightarrow 
    \text{End}(S_x),
\end{equation}
the extension of the Clifford multiplication $c$
to fibers of the spinor bundle and by 
\begin{equation}
    \nabla \: : \: C^{\infty}(S) \longrightarrow
    C^{\infty}(T^*(G/K) \otimes S)
\end{equation}
the covariant 
derivative in the spinor bundle $S$.

\begin{definition}
    The Dirac operator on $\mathbb{H}^d$ is
    defined as the composition 
    \[\begin{tikzcd}[column sep=scriptsize]
	{ c \nabla \: : \:
    \Gamma(S)} && {\Gamma(T(G/K) \otimes S)} && {\Gamma(S),}
	\arrow["\nabla", from=1-1, to=1-3]
	\arrow["c", from=1-3, to=1-5]
\end{tikzcd}\]
where the cotangent bundle has been identified 
with the tangent bundle by means of the metric $g$.
\end{definition}

The commutation rules for the spinorial 
covariant derivative imply the 
Lichnerowicz
formula (see \cite{lichno}) for the non-compact space $\mathbb{H}^d$, relating $(c\nabla)^2$ 
to the spinor Laplacian and the scalar curvature
of $\mathbb{H}^d$.
\begin{proposition}[Lichnerowicz
formula]\label{casimir}
On hyperbolic space, the square of the Dirac operator 
acts on spinors via 
    \begin{equation}
    \begin{split}
    \slashed{\Delta} :=  (c\nabla)^2
        &= \nabla^*\nabla + \frac{\textnormal{scal}(\mathbb{H}^d)}{4},\\
        \end{split}
    \end{equation}
 where  $\textnormal{scal}(\mathbb{H}^d)= -d(d-1)$ is the scalar curvature of
    $\mathbb{H}^d$ and $\nabla^*\nabla$ denotes
    the spinor Laplacian.
\end{proposition}
We now recall some well-known properties 
of Clifford multipliers.
\begin{proposition}[\cite{oliver23}, Lemma 2.1.]\label{cliffordid}
Let $n \in \mathbb{N}$.
\begin{itemize}
    \item For $\alpha \in \{1,\dots,d\}^n$,
    \begin{equation}
        \textnormal{tr}_{S_x} \left(\prod_{j=1}^n c^{\alpha_j} \right) = 0,
    \end{equation}
    provided that ($n,d$ are odd and $n<d$)
    or ($n$ is odd and $d$ is even).
    \item For $\alpha \in \{1,\dots,d\}^n$
    and $d$ odd, 
    \begin{equation}
        \textnormal{tr}_{S_x} \left(\prod_{j=1}^d c^{\alpha^j} \right) = \kappa_c \epsilon_{\alpha},
    \end{equation}
    where $\epsilon_{\alpha} = \epsilon_{\alpha_1 \dots \alpha_d}$ is the
    Levi-Civita symbol and 
    \begin{equation}\label{kappac}
        \kappa_c := \textnormal{tr}_{S_x}
        \left( \prod_{j=1}^d c^j \right).
    \end{equation}
\end{itemize}
\end{proposition}
\begin{rem}
   The minimal choice 
   for the rank of Clifford matrices $(c^i)^d_{i=1} \subset \mathbb{C}^{r \times r}$ that allows them to satisfy 
   \begin{equation*}
       c^ic^j + c^jc^i = -2 \delta_{ij}\one_{\mathbb{C}^r},
   \end{equation*}
  is $r = 2^{\frac{d-1}{2}}$. In that case, 
  Equation \eqref{kappac} reads
  \begin{equation}
      \kappa_c = (2i)^{\frac{d-1}{2}} (-i)^d.
  \end{equation}
\end{rem}

\begin{definition}\label{familyb}
    Let $B$ be a function on $\mathbb{H}^d$ with values in 
    symmetric operators
    on $H$ and 
    common domains $\text{Dom} \: B(x) = \text{Dom} \: A_0$. 
    We say that $B \in C^k_b(\mathbb{H}^d, A_0, m)$, for $k \in 
    \mathbb{N}_0$ and $m \geq 0$, if the following two conditions 
    are satisfied:
    \begin{enumerate}
        \item For $x \in \mathbb{H}^d$, $\langle A_0 \rangle^{\beta} B(x) \langle A_0 \rangle^{-1- \beta}$ 
        is a family of closable operators with closures in 
        $B(H)$ for any $\beta \in [-m,m]$;
        \item $x \mapsto 
        \overline{\langle A_0 \rangle^{\beta} B(x) \langle A_0 \rangle^{-1- \beta}} \: \in \: C^k_b(\mathbb{H}^d, B(H))$ in the strong operator topology.
    \end{enumerate}
\end{definition}

\begin{definition}
    For operators $B \in C^1_b(\mathbb{H}^d, A_0, 1)$ and $z \in \mathbb{C} \backslash (- \infty, 0]$, denote by $\rho_z(B)$ the 
    semi-norms 
    \begin{equation}
    \begin{split}
\rho_z(B) & := \norm{B \langle M_0 \rangle_z^{-1}}_{L^{\infty}(
       \mathbb{H}^d, B(H))} \: + \:
       \norm{\langle M_0 \rangle B \langle M_0 \rangle^{-1}
       \langle M_0 \rangle_z^{-1}}_{L^{\infty}(
       \mathbb{H}^d, B(H))}  \\
 & + \: \norm{ (\Id B)\langle M_0 \rangle^{-1}
       \langle M_0 \rangle_z^{-1}}_{L^{\infty}(
       \mathbb{H}^d, B(H))}.
\end{split}  
    \end{equation} 
In particular, for $B \in C^1_b(\mathbb{H}^d, A_0, \vert \beta 
\vert + 1)$ and $\beta \in \mathbb{R}$, let 
\begin{equation}
    \rho_z^{\beta}(B) := \rho_z( 
    \langle M_0 \rangle^{\beta} B \langle M_0 \rangle^{-\beta}).
\end{equation}
\end{definition}

\begin{rem}\label{trianglerho}
For $\beta \in [-2N,2N]$, we have 
\begin{equation}
    \rho^{\beta}_z(B) \leq \rho^{-2N}_z(B) 
    + \rho_z^{2N}(B).
\end{equation}
See \cite{oliver23}, Remark 2.5. for more details.
\end{rem}

\begin{proposition}[\cite{oliver23}, Proposition 
2.6]\label{domain}
    If a family of functions $B \in C^1_b( \mathbb{H}^d, A_0, 1)$ with values in 
    symmetric operators is such that 
\begin{equation}
    \lim_{\text{dist}(z, (- \infty, 0]) \: \rightarrow \: \infty}  \rho_z(B) = 0,
\end{equation}
then $D_B$ is closed in $L^2(S, H)$ and admits 
$D^*_B = - i c\nabla + B + M_0$ as its 
adjoint. 
Additionally, the operators 
$D_B D_B^*$, $D_B^*D_B$, and $H_B$ are non-negative,
self-adjoint with common domain $\textnormal{Dom} \: H_0$.
\end{proposition}

\begin{definition}
    Consider operators $B \in C_b^1(\mathbb{H}^d, A_0,m)$,
    $m \geq 0$, and $\alpha \geq 0$, 
    $\beta \in \mathbb{R}$, with $\vert \beta -1 \vert \leq m$, $\vert \beta - \alpha \vert \leq m$. Moreover, let $p \geq 2$, $s \geq 0$ 
    and $\nu \in \Gamma(\text{T} \mathbb{H}^d)$.
    Then define a set of semi-norms through
    \begin{equation}
        \tau^{\alpha, \beta, \nu, p, s}(B) :=
        \norm{x \mapsto \langle x \rangle^s
        \langle A_0 \rangle^{-\alpha + \beta}
        (\nu B)(x) \langle A_0 \rangle^{-\beta}}_{L^{\infty}(\mathbb{H}^d, S^p(H))},
    \end{equation}
    and set $\tau^{\alpha, \beta, \nabla, p, s}
:= \sum_{i =1}^d \tau^{\alpha, \beta, \partial_i, p, s}$.
\end{definition}

\begin{definition}[Definition 2.11, \cite{oliver23}]\label{deltadef}
    Let $x,y \in \mathbb{H}^d$ and $T$ be
    an operator valued function over $\mathbb{H}^d$. Set 
    \begin{equation}
        \begin{split}
            (\tau_y T)(x) &:= T(y^{-1}x), \\
            \delta_y &:= \tau_{y^{-1}} - \tau_y \\
            \mathcal{W}^2_y(T) &:= \delta_{\frac{y}{2}}^2(T).
        \end{split}
    \end{equation}
    Moreover, for an elementary tensor
    of operator functions $\textbf{T} := T_1 
    \otimes \dots \otimes T_n$ and a vector
    $\textbf{y} = (y_1, \dots, y_n)$, let
    \begin{equation}
       \begin{split}
           \text{OpProd} ( \textbf{T}) 
           &:= T_1 \cdot ... \cdot T_n, \\
      \tau_{\textbf{y}} \textbf{T}  &:=
      \tau_{y_1} T_1 \otimes \dots \otimes 
      \tau_{y_n} T_n.
       \end{split}
    \end{equation}
\end{definition}

What follows are assumptions on $A$ which 
guarantee that the domains 
are well-behaved, that is $\text{Dom}(DD^*) = 
\text{Dom}(D^*D) = \text{Dom}(H_0)$ and 
$\text{Dom}(D) = \text{Dom}(D^*) = \text{Dom}(D_0)$, see Proposition \ref{domain}. Moreover, the conditions also
ensure that, for
$t>0$, the partial trace $\text{tr}_{S_x}(e^{-tD^*D} - e^{-tDD^*})$ is trace class 
in $L^2(\mathbb{H}^d, H)$.
\begin{hyp}\label{hypothesis}
    Let $1 \leq \alpha < 1 + \frac{2}{d}$ and $N \geq \lfloor \frac{\alpha - 1}{2}(d + 1) \rfloor$.
    Moreover, let $B$ be a family of symmetric 
    operators in $H$ with common domains 
    $\textnormal{Dom} \;  A_0$ such that 
    $B \in C^1_b(\mathbb{H}^d, A_0, 2N +1)$. 
    Assume: 
    \begin{enumerate}
        \item \begin{equation}
            \lim_{\text{dist}(z, (-\infty, 0]) \rightarrow \infty} \rho_z^{-2N}(B) + \rho_z^{2N}(B) = 0,
        \end{equation}
 and that for all $\beta \in [-2N + \alpha, 2N + 1]$,
\begin{equation}\label{eq2.39}
    \begin{split}
        \tau^{\alpha, \beta, \nabla, d, 1}(B) &< \infty. \\
    \end{split}
\end{equation}
 If, on the other hand, $N=0$ and $\alpha >1$, then assume 
 \begin{equation}
      \tau^{\alpha, \alpha, \nabla, d, 1}(B) < \infty.
 \end{equation}

\item Denote by $S_1(o)$ the geodesic unit sphere on 
$\mathbb{H}^d$. Then, 
for $\phi \in \text{Dom} \: A_0$, a.e. $y \in S_1(o)$,
a.e. $x \in \mathbb{H}^d$, and $\gamma \in 
\mathbb{N}_0^d$ with 
$\vert \gamma \vert \leq 1$, 
\begin{equation}\label{hyp2}
    \lim_{R \rightarrow \infty} e^{R + \sqrt{R}}\norm{(
    \partial^{\gamma} B(w) -
    \tau^{*,LC}(w,w')
    \partial^{\gamma}B(w')) \phi}_H = 0,
\end{equation}    
for $w = \exp_{w'}(y)$ and $
w' =  \exp_o(Ry)$,  where 
$\tau^{*,LC}(x,y) \: : \: T_y^* \mathbb{H}^d 
\rightarrow T_x^* \mathbb{H}^d$ is the induced 
Levi-Civita parallel transport on the cotangent
bundle.
\item
\begin{equation}\label{hyp3}
   \norm{\Id B(x)}_{S^1(H)} \sim \: 
\text{Vol}_{\mathbb{H}^d}(S_R(o))^{-\frac{1}{d-1}},
\end{equation}
as $R=d(x,o) \rightarrow \infty$;
\item 
Assume that, for $d=3$, there exists $\delta>0$ such that, for a.e.
$y\in S_1(o)$, $z'\in\mathbb{H}^d$, $c>0$, and
$\phi\in\operatorname{Dom}A_0^2$,
\begin{equation}\label{hyp4.1}
    \sup_{\substack{R\geq 1\\ d(z,o)\leq c\log R}}
\left[
e^{R/2}\bigl\|(B(w)-B(w'))\phi\bigr\|_H
+
e^{R}\bigl\|\mathcal{W}_z^2(B)(\exp_o(Ry))\phi\bigr\|_H
\right]<\infty,
\end{equation}
\begin{equation}\label{hyp4.2}
    \begin{split}
      &\sup_{\substack{R\geq 1\\ d(z,o)\leq c\log R}}
m_d(z)\bigg[
e^{3R/2}\bigl\|(\Id B(w) 
    - \tau^{*,LC}(w,w')
    \Id B( w'))\phi\bigr\|_H  \\
&\qquad \qquad \qquad +
e^{2R}\bigl\|\mathcal{W}_z^2((c(\Id B))\tau)(\exp_o(Ry))\phi\bigr\|_H
\bigg]<\infty,  
    \end{split}
\end{equation}
where $w = \exp_{w'}(y)$, $
w' =  \exp_o(Ry)$, and
\begin{equation}\label{hyp4.3}
   e^R \bigl\|\mathcal{W}_{z'}^2(B)(\exp_o(Ry))\phi\bigr\|_H
+
e^{2R} \bigl\|\mathcal{W}_{z'}^2((c(\Id B))\tau)(\exp_o(Ry))\phi\bigr\|_H
\longrightarrow 0,
\end{equation}
as $R\to\infty$, where
\[
m_d(z):=
\begin{cases}
\min\{1,e^{-d(z,o) \delta} \}, & d=3,\\
1, & d\geq 5.
\end{cases}
\]

    \end{enumerate}
\end{hyp}

\begin{rem}
    Contrary to the Euclidean setting treated in 
    \cite{oliver23}, it is not necessary to single
    out the radial direction in Hypothesis 
    \ref{hypothesis}, as well as in all the
    calculations that follow. This is due to the
    fact that $e^{-r} \in L^d(\mathbb{H}^d)$,
    whereby $1/r \notin L^d(\mathbb{R}^d)$.
\end{rem}

\begin{rem}\label{rem2.17}
    In the flat setting in \cite{oliver23}, where Equation \eqref{hyp2}
  corresponds to $$\lim_{R \rightarrow \infty} R^{\vert \gamma \vert}\norm{(\partial^{\gamma} B(Ry + x) - \partial^{\gamma}(Ry)) \phi}_H = 0,$$
 the additional 
factor of $R$ results from the growth 
difference between  Euclidean surface area and 
the volume of balls of
radius $R$.
However, in $\mathbb{H}^d$,
the
volume and surface area of a ball of radius $R$ grow asymptotically at the
same exponential rate, namely $e^{(d-1)r}$.
To account for this, a factor of  $e^R$ 
is needed in Equation \eqref{hyp2}.
The requirement that the difference decays slightly faster than $e^{-r}$ 
is there for technical reasons, allowing us,
for instance, to use a dyadic decomposition
of $\mathbb{H}^d$ in the proof of 
Lemma \ref{shuffleheatkernellemma}.
\end{rem}

\begin{rem}
    As in \cite{oliver23}, Hypothesis \ref{hypothesis} implies, together with 
    a Kato-Rellich argument, that for a.e. 
    $x_0 \in \mathbb{H}^d$, 
    $A(x_0) = A_0 + B(x_0)$ can 
    be taken as the model operator $A_0$.
\end{rem}

\begin{rem}
    Let us discuss the importance of the
    $\alpha$ parameter that appears in Hypothesis 
    \ref{hypothesis}.
    First of all, specifying $N$ prescribes how 
    regular the function $B$ needs to be on the 
    Hilbert scale generated by $A_0$. 
   Choosing $\alpha$ large, weakens the condition in Equation \ref{eq2.39} on $B$, although, it also
    simultaneously forces one to increase 
    $N$. 
   The possibility of choosing $\alpha \geq 1$ is particularly significant in view of applications. For instance, one may wish to treat first-order differential operators on $\mathbb{R}^d$ as concrete realizations of $A_0$, in which case the underlying Hilbert space is $H = L^2(\mathbb{R}^d)$. 
   Such examples cannot, in general, be accommodated when $\alpha = 1$. 
   Indeed, in this setting, condition \eqref{eq2.39} fails. 
   After applying suitable Schatten–von Neumann class estimates (see, for example, Chapter 4 of \cite{simon}), the issue reduces to the fact that $r^{-d} \notin L^1(\mathbb{R}^d)$. 
    
\end{rem}

\begin{definition}
    Let $B$ satisfy Hypothesis \ref{hypothesis}. 
    For $z \in \mathbb{C}$, set 
    $w  := \text{dist}(z, (- \infty, 0])$
    and for $C \geq 0$ define
    $\mathcal{D}_C(B) := \{ z \in \mathbb{C} \: 
    \vert \: w(z) \geq c \},$
    with $c \geq 1$ being the minimal 
    choice such that
    \begin{equation}
         \left(\rho_1^{-2N}(B)
+ \rho_1^{2N}(B) + 1\right)
      \left(  \rho_z^{-2N}(B) + \rho_z^{2N}(B)\right) 
      \leq \frac{1}{2}, \quad z \in \mathcal{D}_C(B). 
    \end{equation}
\end{definition}

\begin{definition}
    Let $s \in \mathbb{R}$.
    We denote by $Sym_s$ all smooth functions
    $f \: : \: \mathbb{H}^d \rightarrow 
    \mathbb{C}$, such that for all 
    $\gamma \in \mathbb{N}_0^d$,
    \begin{equation}
        \langle X \rangle^{\vert \gamma \vert -s} \partial^{\gamma} f \in 
        L^{\infty}(\mathbb{H}^d).
    \end{equation}
\end{definition}

\section{Schatten-von-Neumann membership of the partial trace}\label{section3}

\begin{lemma}\label{uisb}
    Let $B$ satisfy Hypothesis \ref{hypothesis} 
    and let $t \in [0,2]$, $s \in [-1,1]$,
    $k \in \mathbb{N}$, $k \leq N$,
    and $r \in [-2N,2N]$
    with $r -2(k-1) \in [-2N,2N]$.
    Then there exists a universal constant 
    $C$ such that the operators
\begin{equation}
 R :=   \langle M_0 \rangle^r R_z(T)^k 
    \langle M_0 \rangle^{2k-r-t} \langle c\nabla \rangle^t, 
    \quad T \in \{H_B, DD^*, D^*D \},
\end{equation}
admit extensions in $B(L^2(S,H))$, with
operator norms bounded by $2^k$ 
for $z \in \mathcal{D}_C(B)$.
\end{lemma}
\begin{proof}
    As in the proof of Lemma 2.16. in \cite{oliver23}, we start by writing 
    $R = \bigg( \prod_{i=1}^{k-1} R_i \bigg) S$,
    where 
\begin{equation}
    \begin{split}
        S &:= \langle M_0 \rangle^{-2(k-1) + r}
        R_z(T)
         \langle M_0 \rangle^{2k - r -t}
         \langle c \nabla \rangle^t, \\
         R^{i}_z &:= \langle M_0 \rangle^{-2(i-1) + r}
         R_z(T) \langle M_0 \rangle^{2i-r}, 
         \qquad i = \{1, \dots, k-1 \}.
    \end{split}
\end{equation}
Expanding $R_z(T)$ as a Neumann series around $H_0$
yields then,
\begin{equation}\label{neumann}
    \begin{split}
        S &:= \left[  \sum_{l =0}^{\infty} 
        (-1)^l \big(R_z(H_0) C_k\big)^l 
        \right] \langle M_0 \rangle^{2-t} \langle c \nabla \rangle^t R_z(H_0), \\
        R^{i}_z &:= \left[  \sum_{l =0}^{\infty} 
        (-1)^l \big(R_z(H_0) C_i \big)^l 
        \right] \langle M_0 \rangle^2 R_z(H_0),
    \end{split}
\end{equation}
with 
\begin{equation*}
    C_i := \langle M_0 \rangle^{-2(i-1) + r} 
    \left(A^2 - M_0^2 + \eta i(c \Id B) \right)\langle M_0
    \rangle^{2(i - 1) - r},
\end{equation*}
and where $\eta = -1$ for $T = D^*D$, $\eta = 0$
for $T = H_B$ and $\epsilon = 1$ for $T = DD^*$.
Now, 
\begin{equation*}
    \begin{split}
        &\norm{\overline{R_z(A_0) C_i}}_{B(L^2(S,H))}
         \leq \norm{\langle M_0 \rangle 
        \langle M_0 \rangle_z R_z(H_0)}_{B(L^2(S,H))} \\
        & \cdot \norm{\overline{\langle M_0 \rangle_z^{-1}
        \langle M_0 \rangle^{-2(i-1) + r} 
    \left(B^2 + BM_0 + M_0B + \eta i(c \Id B) \right)\langle M_0
    \rangle^{2(i - 1) - r}}}_{B(L^2(S,H))},
    \end{split}
\end{equation*}
where the first term on the right hand side can be 
bounded by $1$ given that
\begin{equation*}
    \sup_{\lambda \in \mathbb{R}} 
    \frac{(\lambda^2 + 1)^{1/2}(\lambda^2 +z)^{1/2}}{\lambda^2 +z} \leq 1,
\end{equation*}
which leads to 
\begin{equation*}
\begin{split}
   \norm{\overline{R_z(A_0) C_i}}_{B(L^2(S,H))} &\leq \left(\rho_1^{2(i-1) - r}(B) +1  \right) \rho_z^{2(i-1) -r}(B) \\
   &\leq \frac{1}{2},
   \end{split}
\end{equation*}
for $z \in \mathcal{D}_z(B)$ and by Remark \ref{trianglerho}.
Hence, the Neumann series of $R^{i}$ and $S$
converge in operator norm with bound $2$ 
for 
$z \in \mathcal{D}_z(B)$, which lets us conclude that
\begin{equation*}
    \norm{\overline{R}}_{B(L^2(S,H))} \leq 2^k, \qquad z \in \mathcal{D}_z(B).
\end{equation*}
    
\end{proof}

\begin{lemma}\label{fourier}
Provided that
$\langle X \rangle^{\frac{2d-2}{p}}f \in L^p(S^p(H))$ with $s> d/p$, then
$$M_f \cdot (\Delta^S + 1)^{-s/2} \in S^p(L^2(\mathbb{H}^d)).$$   
\end{lemma}
\begin{proof}
 Recall that if $f$ is a function on $\mathbb{H}^d$,
    its Fourier-Helgason transform $\mathcal{F}(f)$ is 
    given by 
    \begin{equation}
        \mathcal{F}(f) = \int_{G/K} f(x) 
        e^{(-i\lambda + \mathcal{\rho}) A(x,b)} 
        \Id x,
    \end{equation}
    where $\lambda \in \mathfrak{a}^*_{\mathbb{C}} \simeq \mathbb{C}$, $b \in B = K/M = 
 \text{Spin}(d)/\text{Spin}(d-1)   
    \simeq \mathbb{S}^{d-1}$,
    the $d-1$-dimensional sphere in $\mathbb{R}^d$,
    and 
$A(x,b)$ is defined via the Iwasawa decomposition, 
i.e. in $G = NAK = KAN$, one has 
$g \in N \exp A(g) K$ and $g \in K \exp H(g) N$, 
with unique $A(g),H(g) \in \mathfrak{a}$, 
such that $A(g) = - H(g^{-1})$ and 
\begin{equation}
 A(x,b)=   A(gK,kM) = A(k^{-1}g).
\end{equation}
Moreover (see e.g. \cite{helgason},  Theorem 1.5, p. 227), the Fourier transform extends to an isometry of $L^2(G/K)$ onto 
$L^2(X) = L^2(\mathfrak{a}^*_+ \times B; \vert c(\lambda)\vert^{-2} \Id \lambda \Id b)$, 
where $\mathfrak{a}^*_+ \simeq \mathbb{R}_+$ 
and $c$ is the Harish-Chandra $c$-function. \\
Let now $g(x) := (x^2+1)^{-s/2}$.
 Then $\mathcal{F} g(c\nabla) \mathcal{F}^{-1} = M_g \in L^2(X)$.
We now argue that 
$M_f \cdot g(c\nabla) \in S^2(L^2(
\mathbb{H}^d))$ by showing that
$$\int_{(\mathbb{H}^d)^2} \vert f(x) g(c\nabla)(x,y) 
\vert^2 \Id y \Id x < \infty,$$
where the growth condition on $f$ is to be 
determined.
Since we can separate the $x$ and $y$ integration, we
obtain  an expression involving
$$\int_{\mathbb{H}^d} g(c\nabla)(x,y) \overline{g}(c\nabla)(y,x) \Id y,$$
which reduces to 
the study of $\vert g \vert^2(c\nabla)(x,x)$ using   the fact that $\overline{g(c\nabla)(x,y)}
= \overline{g}(c \nabla)(y,x)$.
To that end, write
\begin{equation*}
\begin{split}
    g(c \nabla)(x,y) &= \int_X h(x, \xi, \lambda) 
    g\left( \sqrt{\lambda^2 + \rho^2}\right) \tilde{h}( \xi, \lambda, y) \: \Id \xi \Id \lambda, \\
    \end{split}
\end{equation*}
with 
\begin{equation*}
    \begin{split}
        h(x, (\lambda, \xi)) &:= e^{(-i \lambda - \rho) B(x,\xi)}, \\
        \tilde{h}(x, (\lambda, \xi)) &:= 
        e^{(i \lambda - \rho) B(x,\xi)}
        \vert c(\lambda) \vert^{-2},
    \end{split}
\end{equation*}
and the fact that the Fourier transform diagonalizes the spinor Laplacian.
Thus,
\begin{equation*}
   \begin{split}
        \vert g \vert^2(c\nabla)(x,x) 
    &= \int_X h(x, (\lambda, \xi) )
    \tilde{h}((\lambda, \xi), x) \vert g \vert^2 \left( 
  \sqrt{ \lambda^2 + \rho^2}
  \right) \Id \xi \Id \lambda \\
  &= \int_X e^{-2 \rho B(x, \xi)}
  \vert g \vert^2 \left( 
  \sqrt{ \lambda^2 + \rho^2}
  \right) \vert c(\lambda) \vert^{-2}
  \Id \xi \Id \lambda \\
  &= \int_{\mathbb{S}^{d-1}} e^{-2 \rho B(x, \xi)} 
  \Id \xi \cdot \int_{\mathbb{R}}
  \vert g \vert^2 (\langle \lambda \rangle_{\rho} )
  \vert c(\lambda) \vert^{-2} 
  \Id \lambda.
    \end{split}
\end{equation*}
It is known (\cite{helgason2}, Proposition 7.2.) that $\vert c(\lambda)\vert^{-2}$
can be estimated via
\begin{equation}
    \vert c(\lambda)\vert^{-2} \leq C(1 + \vert \lambda \vert)^{\beta}, \quad \lambda \in \mathfrak{a}^*,
\end{equation}
with $\beta = \text{dim}(N) = d-1$
and that it is regular at zero. 
On the other hand, 
in the Poincaré Ball model,
the Busemann function $B$ takes on the form 
$B(x, \xi) = \log \left( \frac{1- \vert x \vert^2}{
\vert x - \xi \vert^2
}
\right),$
so that 
$$e^{-2 \rho B(x, \xi)} =
\left( \frac{\vert x - \xi \vert^2}{1- \vert x \vert^2} \right)^{d-1},$$
and therefore $\int_{\mathbb{S}^{d-1}} e^{-2 \rho B(x, \xi)} 
  \Id \xi$ behaves  like $\langle x \rangle^{2(d-1)}$
  for large $x$.
Hence, at infinity, we need $-2s + d-1 < -1$, i.e. $ s > d/2$.
It follows that $M_f \cdot (\Delta^S +1)^{-s/2} \in S^2(L^2(\mathbb{H}^d))$, provided that 
$\langle X \rangle^{d-1}f \in L^2(S^2(H))$
and $s > d/2$.
Complex interpolation then implies 
that given 
$\langle X \rangle^{\frac{2d-2}{p}}f \in L^p(S^p(H))$ with $s> d/p$, then
$M_f \cdot (\Delta^S + 1)^{-s/2} \in S^p(L^2(\mathbb{H}^d))$.
    
\end{proof}

\begin{lemma}\label{viss}
    Assume that $B$ satisfies Hypothesis \ref{hypothesis}.
    Let $\beta \in \mathbb{R}$
    such that $\vert \beta - \alpha \vert 
    \leq 2N$ and $\vert \beta -1 \vert \leq 2N$.
 Then, for all $\delta >0$, 
\begin{equation}
    \langle c \nabla \rangle^{-1-\delta}
    \langle M_0 \rangle^{\beta -\alpha} 
   (ic \Id B) \langle M_0 \rangle^{-\beta},
\end{equation}
possesses an $S^d(L^2(S,H))$-extension. Moreover, 
\begin{equation}\label{eqdeltazero}
    \langle c \nabla \rangle^{-1}
    \langle M_0 \rangle^{\beta -\alpha} 
   (ic \Id B) \langle M_0 \rangle^{-\beta},
\end{equation}
possesses an $S^{d + \delta}(L^2(S,H))$-extension.
\end{lemma}
\begin{proof}
Since 
\begin{equation}
    x \mapsto \norm{ \langle M_0 \rangle^{\beta - \alpha} (c \Id B) \langle M_0 \rangle^{- \beta}}_{S^d(L^2(S,H))}
\end{equation}
behaves like $e^{-\text{dist}_{\mathbb{H}^d}(x,o)}$ as
$x \rightarrow \infty$ and 
\begin{equation*}
    \langle c\nabla
\rangle^{-1- \delta} \simeq 
\langle \Delta^S \rangle^{-\frac{1 + \delta}{2}} \simeq (\Delta^S + 1 )^{-\frac{\delta +1}{2}},
\end{equation*}
we need that $\delta + 1 >1$ in order to apply 
Lemma \ref{fourier}, which is 
guaranteed by $\delta >0$.
Therefore
\begin{equation*}
    \langle X \rangle^{\frac{2d-2}{d}}
    e^{-X} \in L^d(S^d(H)),
\end{equation*}
and we can conclude. Equation \eqref{eqdeltazero}
then follows by a similar argument as in the proof of Lemma \ref{fourier}, where the $S^{d + \delta}$-membership is deduced by standard integrability results
in $\mathbb{R}$.
\end{proof}

Analogous to Lemma 2.25. in \cite{oliver23},
we show the following:
\begin{lemma}\label{decomposition1}
    Let $z \in \mathbb{C} \backslash (-\infty, 0]$
    and assume Hypothesis \ref{hypothesis}.
    Then
\begin{equation}
  \begin{split}
      &
      \textnormal{tr}_{S}\big( R_z(D^*_B D_B)^N 
      - R_z(D_BD_B^*)^N\big)  \\
      &=2 \textnormal{tr}_{S} 
     \sum_{\substack{\gamma \in \mathbb{N}^{d+1} \\ \vert \gamma \vert = N + d}}
      \begin{pmatrix}
          N-1 \\
          \gamma - \one 
      \end{pmatrix}^{-1}
      R_z(H_B)^{\gamma} \shuffle (ic \Id B)^{\otimes d} \\
      & + \textnormal{tr}_{S} 
      \sum_{\substack{\gamma \in \mathbb{N}^{d+3} \\ \vert \gamma \vert = N + d + 2}}
       \begin{pmatrix}
          N-1 \\
          \gamma - \one 
      \end{pmatrix}^{-1} \bigg[ \left(R_z(H_B)^{\otimes(d+2)}, R_z(D^*_B D_B) \right)^{\gamma} \shuffle
      (ic \Id B)^{\otimes (d+2)} \\
      &+ \left(R_z(H_B)^{\otimes(d+2)}, R_z(D_B D_B^*)\right)^{\gamma}  \shuffle (ic \Id B)^{\otimes(d+2)} \bigg].
  \end{split}
\end{equation}

\end{lemma}

\begin{definition}\label{uandv}
    Let $k,l \in \mathbb{N}_0$, $\beta \in \mathbb{R}^{l+1}$, $\gamma \in \mathbb{N}^{l+1}$, 
     $\alpha \geq 0$
    and $\delta > 0$. Let $B$ be a family of operators that satisfies Hypothesis
    \ref{hypothesis}.
Define for 
$a,b,c \in \mathbb{R}$,  operators
\begin{equation*}
    \begin{split}
        u_z^{a,b,c,k}(B) &:= \langle M_0 \rangle^a
         R_z(H_B)^k 
        \langle M_0 \rangle^b
        \langle c \nabla \rangle^c, \\
        \big(U_z^{\alpha, \beta, \gamma, \delta}(B) \big)_i &:= 
        \begin{cases}
            u_z^{\beta_i, \alpha 
            - \beta_{i+1}, 1 + \delta, \gamma_i}(B), \quad 0 \leq i \leq l-1 \\
            u_z^{\beta_l, 0,0, \gamma_l}(B), 
             \quad  i =l,
        \end{cases} \\
        \big(V^{\alpha, \beta, \delta}(B) \big)_i &:= \langle c \nabla 
        \rangle^{-1 - \delta} \langle M_0 \rangle^{\beta_i - \alpha}  c(\Id B)
        \langle M_0 \rangle^{- \beta_i}, \quad 
        i \in \{1, \dots, l\}. \\
    \end{split}
\end{equation*}  
    \end{definition}

\begin{proposition}\label{decomposition2}
Under the assumptions made in Hypothesis \ref{hypothesis}, there exist, independently of the choices of
$B$ and $z \in \mathcal{D}_C(B)$, 
multi-indices $\beta \in [\alpha -2N, 2N + \alpha] \cap [1-2N, -1 +2N]$, a constant $\delta \in [0,1]$
and $C^1_{\gamma}, C^2_{\gamma} < \infty$
for $\gamma \in \mathbb{N}^{d+1}$, such that
\begin{equation}
    \begin{split}
         &
      \textnormal{tr}_{S}\big( R_z(D^*_B D_B)^N 
      - R_z(D_BD_B^*)^N\big)  \\ 
      &= 
      \textnormal{tr}_{S} 
      \sum_{\substack{\gamma \in \mathbb{N}^{d+1} \\ \vert \gamma \vert = N + d}}
     C_{\gamma}^1 \cdot 
      U_z^{\alpha, \beta, \gamma, \delta}(B)
    \shuffle V^{\alpha, \beta, \delta}(B) \\
    &+ 
    \textnormal{tr}_{S}
    \sum_{\substack{\gamma = 
    (\tilde{\gamma}, \gamma_{d+2})\in \mathbb{N}^{d+2}
    \times  \mathbb{N}\\ \vert \gamma \vert = N + d
    +2}} 
    C^2_{\gamma} \cdot \big(U_z^{\alpha, \beta, 
    \tilde{\gamma},0}(B) \shuffle V^{\alpha, \beta,0}(B) \big) \\
  &\cdot i(c\nabla B)\big(R_z(D^*_B D_B)^{\gamma_{d+2}}
  + R_z(D_B D_B^*)^{\gamma_{d+2}} \big),
    \end{split}
\end{equation}
for $z \in \mathcal{D}_C(B)$.
Then each summand is an operator in $S^1(L^2(S,H))$,
with $S^1$-norm bound given by a $z$-independent 
polynomial in $\tau^{\alpha, \beta, \nabla, d ,1}(B)$
for $\beta \in I$, where 
$I \subset [\alpha -2N, 2N + \alpha] \cap [1-2N, -1 +2N]$ are $z$-independent finite subsets.
    
\end{proposition}
\begin{proof}
    We consider the decomposition from 
    Lemma \ref{decomposition1} and start with the 
    first sum, which contains terms 
    of the form 
    \begin{equation*}
        R_z(H_B)^{\gamma} \shuffle (ic \Id B)^{\otimes d},
    \end{equation*}
    for $\gamma \in \mathbb{N}^{d+1}$, 
    $\vert \gamma \vert = N + d$.
We then rewrite those terms with the aid of the operators introduced in Definition \ref{uandv},
namely
\begin{equation}
    U_z^{\alpha, \beta, \gamma, \delta}(B)
    \shuffle V^{\alpha, \beta, \delta}(B),
\end{equation}
with $\beta_0 := 0$, $\beta_{i+1} := \beta_i + \alpha  \delta +1 - 2 \gamma_i$ and $\gamma \in \mathbb{N}^{d+1}$, 
$\vert \gamma \vert = N +d$.
We then have $\overline{(V^{\alpha, \beta, \delta}(B))_i}
\in S^d(L^2(S,H))$ by  Lemma \ref{viss}
for all $i$, $
\overline{(U_z^{\alpha, \beta, \gamma, \delta}(B))_i}
\in B(L^2(S,H))$
by Lemma \ref{uisb} for $i \in \{0, \dots, d-1 \}$
and 
$
\overline{(U_z^{\alpha, \beta, \gamma, \delta}(B))_d}
\in B(L^2(S,H))$ since
$(\alpha +1)d < 2(N+d)$, as assumed under Hypothesis 
\ref{hypothesis}, so that we can choose $\delta >0$
small enough such that $\beta_d < 2 \gamma_d$,
which lets us conclude with another application of
Lemma \ref{uisb}. \\
Let us now consider the second sum in the decomposition
of Lemma \ref{decomposition1}, indexed over 
$\gamma \in \mathbb{N}^{d+3}$ with 
$\vert \gamma \vert = N +d +2$. We then rewrite
every term of the form
\begin{equation*}
  \left(  R_z(H_B)^{\otimes(d+2)},R_z(T) \right)^{\gamma} \shuffle (ic \Id B)^{\otimes (d+2)}, 
\end{equation*}
where $T \in \{D^*D, DD^* \}$, as 
\begin{equation*}
   \left( U_z^{\alpha, \beta, \tilde{\gamma}, 0}(B)
    \shuffle V^{\alpha, \beta, 0}(B)  \right)
    ic(\Id B) \langle M_0 \rangle^{-1} 
    \langle M_0 \rangle R_z(T)^{\gamma_{d+2}},
\end{equation*}
with $\tilde{\gamma} := (\gamma_0, \dots, \gamma_{d +1})$, $\beta_0 := 0$ and $\beta_{i+1} := 
\beta_i + \alpha +1 - \gamma_i$.
We then have 
$\overline{\big(U_z^{\alpha, \beta, \tilde{\gamma}, 0}(B)\big)_i} \in B(L^2(S,H))$ by Lemma \ref{uisb} 
and $ic(\Id B)\langle M_0 \rangle^{-1} 
\in B(L^2(S,H))$ by the assumptions made in Hypothesis
\ref{hypothesis} and the fact that for a measurable
operator $W$, $\norm{W}_{B(L^2(S,H))} = 
\norm{W}_{L^{\infty}(\mathbb{H}^d, B(S \otimes H))}$.
Moreover, $\langle M_0 \rangle R_z(T)^{\gamma_{d+2}}
\in B(L^2(S,H))$
again by Lemma \ref{uisb} and 
$\overline{\big(V^{\alpha, \beta, 0}(B) \big)_i} \in 
S^{d+ 1}(L^2(S,H))$ by Lemma \ref{viss}.
We then have $z$-independent finite 
sets $I$ with 
$I \subset [\alpha -2N, 2N + \alpha] \cap [1-2N, -1 +2N]$ and $z$-independent constants $C$ such that
the $S^d$-, respectively $S^{d+1}$-norms 
of the numerous $\overline{V^-(B)}$ operators 
can be bounded uniformly by
\begin{equation*}
    C \sum_{\beta \in I} \tau^{\alpha, \beta, \nabla, d ,1}(B)
    < \infty.
\end{equation*}

\end{proof}

\begin{proposition}
    For $t>0$ and $B$ satisfying Hypothesis
    \ref{hypothesis}, we have
    \begin{equation}
        \textnormal{tr}_{S_x}\left(e^{-tD^*D} -
        e^{-tDD^*} \right) \in S^1(L^2(\mathbb{H}^d,H)).
    \end{equation}
\end{proposition}
\begin{proof}
    The proof follows by holomorphic functional 
    calculus and is completely analogous 
    to the proof of Proposition 2.28. 
    in \cite{oliver23}.
\end{proof}

\section{Approximation scheme for the
operator valued function $B$}\label{section4}
We construct an approximation scheme for the operator-valued function 
$B$ through finite-rank operators
with improved regularity while simultaneously ensuring that 
$\Id B$ is both more regular and trace-class valued. 
These properties are essential for computing the iterated trace in the sections that follow.
Since, in contrast to the Euclidean case, there is no need to distinguish a radial direction, the resulting approximation scheme is simpler than that employed in the Euclidean setting considered in \cite{oliver23}. \\
We start by establishing five technical lemmas addressing convergence of certain operators in the Schatten–von Neumann norm. These lemmas play a key role in the proof of the central result of this section, namely Lemma \ref{approximationscheme}.

\begin{lemma}[Lemma 3.4. \cite{gesztesy11}]\label{sotalemma}
    Let $R_n, S_n, T_n$, $n \in \mathbb{N}$
    be bounded operators in a 
    separable Hilbert space $X$ 
    and for $1 \leq p <  \infty$,
let $S_n \in S^p(X)$, such that
$\lim_{n \rightarrow \infty} S_n = S$
in $S^p(X)$-norm.
If $\lim_{n \rightarrow \infty} R_n = R$
and $\lim_{n \rightarrow \infty} T_n^* = T^*$ 
strongly, then
\begin{equation}
    \lim_{n \rightarrow \infty} R_n S_n T_n 
    = RST,
\end{equation}
in $S^p(X)$-norm.
\end{lemma}

\begin{definition}
    For a bounded operator $T$, and 
    $(T_n)_{n \in \mathbb{N}}$ a sequence of
    bounded operators in a Hilbert space $X$,
    we say that $T_n$ converges in SOTA to
    $T$, if $T_n$ and $T_n^*$ both converge in
    the strong operator topology to
    $T$, respectively $T^*$ as 
    $n \rightarrow \infty$.
\end{definition}

\begin{lemma}\label{convergence1}
    Assume that $B$ satisfies Hypothesis
    \ref{hypothesis} and let $(B_n)_{n\in \mathbb{N}}$ be a sequence of operators 
    each of which also satisfy 
    Hypothesis \ref{hypothesis}. 
    Suppose that, for $\beta \in [-2N+1, 2N-1]$
    and a.e. $x \in \mathbb{H}^d$, 
\begin{equation}\label{eqconvergence11}
    \langle A_0 \rangle^{\beta} (B_n(x) + A_0)^2 \langle A_0 \rangle^{-2 - \beta} 
\xrightarrow[]{n \rightarrow\infty}
 \langle A_0 \rangle^{\beta} (B(x) + A_0)^2 \langle A_0 \rangle^{-2 - \beta},
\end{equation}
converges in the strong operator topology of
$H$, and, for $w := \textnormal{dist}(z, (-\infty, 0])$,  
\begin{equation}\label{eqconvergence12}
    \sup_{n \in \mathbb{N}} \left(\rho_z^{-2N}(B_n) + \rho_z^{2N}(B_n) \right) 
    \xrightarrow[]{w \rightarrow \infty} 0.
\end{equation}
In that case, 
\begin{equation}
\begin{split}
    &\overline{ \langle M_0 \rangle^r  R_z(H_{B_n})^k  \langle M_0 \rangle^{2k-r-t} 
    \langle c\nabla \rangle^t} \\
& \xrightarrow[]{n \rightarrow\infty}    
 \overline{ \langle M_0 \rangle^r  R_z(H_{B})^k  \langle M_0 \rangle^{2k-r-t} 
    \langle c\nabla \rangle^t},
    \end{split}
\end{equation}
converges in SOTA of $L^2(S,H)$
for $t \in [0,2]$, 
$k \in \mathbb{N}$, $k \leq N$, 
and $r \in [-2N,2N]$ with $r-2(k-1) \in [-2N,2N]$.
Furthermore, if for a.e. $x \in \mathbb{H}^d$
we also have strong operator convergence
in $H$ of 
\begin{equation}
    \langle A_0 \rangle^{\beta} \partial_i B_n(x) \langle A_0 \rangle^{-2- \beta}
 \xrightarrow[]{n \rightarrow\infty}    
  \langle A_0 \rangle^{\beta} \partial_i B(x) \langle A_0 \rangle^{-2- \beta}, 
  \quad i \in \{1, \dots, d \},
\end{equation}
then 
\begin{equation}
    \begin{split}
    &\overline{ \langle M_0 \rangle^r  R_z(D_{B_n}D_{B_n}^*)^k  \langle M_0 \rangle^{2k-r-t} 
    \langle c\nabla \rangle^t} \\
& \xrightarrow[]{n \rightarrow\infty}    
 \overline{ \langle M_0 \rangle^r  R_z(D_B D_B^*)^k  \langle M_0 \rangle^{2k-r-t} 
    \langle c\nabla \rangle^t},
    \end{split}
\end{equation}
and 
\begin{equation}
    \begin{split}
    &\overline{ \langle M_0 \rangle^r  R_z(D_{B_n}^*D_{B_n})^k  \langle M_0 \rangle^{2k-r-t} 
    \langle c\nabla \rangle^t} \\
& \xrightarrow[]{n \rightarrow\infty}    
 \overline{ \langle M_0 \rangle^r  R_z(D_B^* D_B)^k  \langle M_0 \rangle^{2k-r-t} 
    \langle c\nabla \rangle^t},
    \end{split}
\end{equation}
converge in SOTA of $L^2(S,H)$.
    
\end{lemma}
\begin{proof}
We proceed similarly as in the proof of Lemma 
\ref{uisb}. To show the statement for $H_B$, we 
replace each occurrence of $B$ by $B_n$ and note
that the Neumann series \eqref{neumann} converge 
for $z \in \mathcal{D}_C(B_n)$. 
Since we are assuming \eqref{eqconvergence12},
there is an $n$-independent  $w_0 \geq 1$, such that
\begin{equation*}
    \{z \in \mathbb{C} \: : \:w(z) \geq w_0 \} 
    \subseteq \mathcal{D}_C(B) \cap \bigcap_{n \in \mathbb{N}} \mathcal{D}_C(B_n).
    \end{equation*}
Hence, after repeating the argument in the proof
of Lemma \ref{uisb}, we see that 
$$S_n(t) := \langle M_0 \rangle^{-2(k-1) + r}
        R_z(H_{B_n})
         \langle M_0 \rangle^{2k - r -t}
         \langle c \nabla \rangle^t,$$
and $$S(t) := \langle M_0 \rangle^{-2(k-1) + r}
        R_z(H_{B})
         \langle M_0 \rangle^{2k - r -t}
         \langle c \nabla \rangle^t,$$ have bounded extensions for $w(z) \geq w_0$ and that their
operator norms are uniformly bounded by $2$.
The resolvent identity yields then 
\begin{equation*}
    \overline{S_n(t)} - \overline{S(t)}
    = \overline{S_n(0)} \overline{
    \left( \langle M_0 \rangle^{-2k + r} \left((B+M_0)^2-(B_n+M_0)^2  \right) \langle M_0 \rangle^{2(k-1) -r}
    \right)} \overline{S(t)},
\end{equation*}
where $(B_n+M_0)^2$ converges strongly to
$(B+M_0)^2$ in $L^2(S,H)$
by assumptions \eqref{eqconvergence11}
and \eqref{eqconvergence12}, and 
dominated convergence.
Since the operator norm of $\overline{S_n(0)}$ is uniformly bounded
 by $2$, we conclude that we have strong operator convergence of $\overline{S_n(t)}$
 to $\overline{S(t)}$.
 The terms $R^{i}_{z,n} := \langle M_0 \rangle^{-2(i-1) + r}
         R_z(H_{B_n}) \langle M_0 \rangle^{2i-r}$
for $i = \{1, \dots, k-1 \}$ can then be treated in a similar manner, as well as the respective statement for the adjoints, where one uses 
that the functions $B_n$ and $B$ are point-wise symmetric and considers
\begin{equation*}
   S_n(t)^* - S(t)^*
    =S_n(t)^* \overline{
    \left( \langle M_0 \rangle^{2(k-1) -r} \left((B+M_0)^2-(B_n+M_0)^2  \right) \langle M_0 \rangle^{
    -2k + r}
    \right)} S(0)^*.
\end{equation*}
The statements for $D_BD_B^*$ and $D_B^*D_B$ are then completely 
analogous.

\end{proof}
Similar to the proof of Lemma \ref{viss}, one shows the following.
\begin{lemma}\label{extension2}
If $B \: : \: \mathbb{H}^d \rightarrow B(H)$
is a point-wise self-adjoint $C^1$-function
with respect to the strong operator topology and
\begin{equation*}
    \tau^{0,0,\nabla,d,1} < \infty,
\end{equation*}
then for all $\delta >0$,
    \begin{equation*}
        \langle c\nabla \rangle^{-1-\delta} (ic\Id B),
    \end{equation*}
    admits and $S^d(L^2(S,H))$-extension
    and 
    \begin{equation*}
        \langle c\nabla \rangle^{-1} (ic \Id B),
    \end{equation*}
    an $S^{\delta +1}(L^2(S,H))$-extension.
\end{lemma}

\begin{lemma}\label{lemmabound}
    Let $B \: : \: \mathbb{H}^d \rightarrow B(H)$
    be a point-wise self-adjoint $C^1$-function
with respect to the strong operator topology 
such that 
\begin{equation*}
    \norm{B}_{C^1_b(\mathbb{H}^d,B(H))}
    = \norm{B}_{L^{\infty}(\mathbb{H}^d,B(H))} 
    + \norm{\Id B}_{L^{\infty}(\mathbb{H}^d,B(H)^d)} < \infty,
\end{equation*}
and $B \in C^1_b(\mathbb{H}^d,A_0,1)$
with 
\begin{equation*}
    \lim_{w \rightarrow \infty} \rho_z(B) = 0,
\end{equation*}
where $w := \text{dist}(z, (- \infty,0])$.
Then for $t \in [0,2]$, $k \in \mathbb{N}$,
the operators
\begin{equation*}
    R_z(H_B)^k \langle c \nabla \rangle^t,
\end{equation*}
admit extensions in $B(L^2(S,H))$ with norm
bounded by $C^k(\norm{B}_{C^1_b(\mathbb{H}^d,B(H))}+1)^{2k}$
for $w \geq C(\norm{B}_{C^1_b(\mathbb{H}^d,B(H))}+1)^4$
and some universal constant $C$. \\
Moreover, 
for $T_B \in \{H_B, D_B^*D_B,D_BD_B^* \}$ 
and $w \geq C(\norm{B}_{C^1_b(\mathbb{H}^d,B(H))}+1)^4$, the operators
\begin{equation*}
    (H_0 +1)^{1/2} R_z(T_B) (H_0 +1)^{1/2},
\end{equation*}
admit extensions in 
$B(L^2(S,H))$ with their norms bounded by 
$C^k(\norm{B}_{C^1_b(\mathbb{H}^d,B(H))}+1)^{2k}$.
\end{lemma}
\begin{proof}
The proof proceeds via algebraic manipulations involving the Neumann series representation of the resolvents and the resolvent identity, and is completely analogous to the proof of Lemma 2.17. in \cite{oliver23}.
\end{proof}

\begin{lemma}\label{sota2}
    Let $B \: : \: \mathbb{H}^d \rightarrow B(H)$
    and $B_n \: : \: \mathbb{H}^d \rightarrow B(H)$,
    $n \in \mathbb{N}$,
be  point-wise self-adjoint $C^1$-functions
with respect to the strong operator topology 
and suppose that $\norm{B_{}}_{C^1_b(\mathbb{H}^d,B(H))}$
and $\norm{B_n}_{C^1_b(\mathbb{H}^d,B(H))}$ are 
uniformly bounded. 
Assume furthermore that $B,B_n \in C^1_b(\mathbb{H}^d, A_0,1)$ and that for all $n \in \mathbb{N}$,
\begin{equation*}
    \lim_{w \rightarrow \infty} \rho_z(B) 
    = \lim_{w \rightarrow \infty} \rho_z(B_n) =0,
\end{equation*}
where $w := \textnormal{dist}(z, (-\infty, 0])$, as well as that 
$B_n(x) \rightarrow B(x)$, $\Id B_n(x) \rightarrow \Id B(x)$ as $n \rightarrow \infty$
for a.e. $x \in \mathbb{H}^d$ in the strong operator 
topology of $H$.
Then there exist constants $C,C'$ such that for 
$w := \textnormal{dist}(z, (-\infty, 0]) \geq C'$, 
\begin{equation*}
    \overline{R_z(H_{B_n})^k \langle c\nabla \rangle^t} \longrightarrow \overline{R_z(H_{B})^k \langle c\nabla \rangle^t}, \quad n \rightarrow \infty,
\end{equation*}
converges in the strong operator topology of 
$L^2(S,H)$ for $t \in [0,2]$ and $k \in \mathbb{N}$, and even in SOTA in
$L^2(S,H)$ for $t =0$. 
For $w \geq 1$, we also have SOTA convergence in
$L^2(S,H)$ of 
\begin{equation*}
    R_z(D^*_{B_n}D_{B_n})^k \longrightarrow R_z(D^*_{B}D_{B})^k, \quad 
     R_z(D_{B_n}D^*_{B_n})^k \longrightarrow R_z(D_{B}D^*_{B})^k,
\end{equation*}
as $n \rightarrow \infty$.

\end{lemma}
\begin{proof}
By Lemma   \ref{lemmabound} and the uniform
boundedness principle, there are constants 
$C,C'$ such that for $w \geq C'$, we have 
$\norm{\overline{R_z(H_{B_n})^k\langle c\nabla \rangle^t}}_{B(L^2(S,H))} \leq C$, 
respectively the statement when $B_n$ is
replaced by $B$. 
After factorizing the operators, we can take $k=1$
and using the uniform norm bounds from 
Lemma \ref{lemmabound} and the fact
that $C^{\infty}_c(\mathbb{H}^d)$ is dense
in $L^2(\mathbb{H}^d)$, even $t=0$. The 
resolvent identity then yields
\begin{equation*}
    \overline{R_z(H_{B_n})}- \overline{R_z(H_{B}) }
    = \overline{R_z(H_{B_n})} \left(
    (B + M_0)^2 - (B_n + M_0)^2
    \right) \overline{R_z(H_{B}) },
\end{equation*}
which lets us deduce that 
$\overline{R_z(H_{B_n})}$ converges in the strong
operator topology to $\overline{R_z(H_{B}) }$,
given that $\overline{\overline{R_z(H_{B_n})}M_0}$ and 
$M_0 \overline{R_z(H_{B})}$ are uniformly 
bounded by Lemma \ref{lemmabound}. \\
The equivalent statements for the adjoints 
follow similarly from the uniform bounds in
Lemma \ref{lemmabound} and the 
resolvent identity.
\end{proof}

\begin{definition}\label{delta}
Let $\psi\in C_c^\infty([0,\infty))$ be nonnegative and supported in $[0,1)$.
Assume that $v\mapsto\psi(|v|)$ is smooth on $\mathbb R^d$ and 
for $l\in\mathbb N$, set
\[
 Z_l
 :=\int_{\mathbb H^d}
       l^d\psi\bigl(l\,d(o,y)\bigr)\,
       \mathrm d\operatorname{vol}(y),
 \qquad
 k_l(x,y)
 :=\frac{l^d}{Z_l}\psi\bigl(l\,d(x,y)\bigr),
\]
where, by homogeneity of $\mathbb H^d$, the normalization
constant is independent of the choice of base point.
In particular,
\[
 \int_{\mathbb H^d}k_l(x,y)\,
       \mathrm d\operatorname{vol}(y)=1
 \qquad\text{for every }x\in\mathbb H^d.
\]
For a Banach space $E$ and $f\in C(\mathbb H^d;E)$,
define
\[
 (\delta_l f)(x)
 :=\int_{\mathbb H^d}k_l(x,y)f(y)\,
       \mathrm d\operatorname{vol}(y),
\]
where the integral is understood in the Bochner sense.
\end{definition}

\begin{definition}
Let $E$ be a Banach space and $k>0$. Define
\[
 \mathcal G^k(\mathbb H^d;E)
 :=
 \left\{
 f\in C^1(\mathbb H^d;E):
 \sup_{x\in\mathbb H^d}
 e^{k d(x,o)}\|f(x)\|_E<\infty
 \right\},
\]
and write
\[
 \|f\|_{\mathcal G^k}
 :=
 \sup_{x\in\mathbb H^d}
 e^{k d(x,o)}\|f(x)\|_E.
\]
When the domain and target space are clear, we
simply write $\mathcal G^k$.
\end{definition}

\begin{lemma}\label{lemma4.3}
Let $k>0$ and $f\in\mathcal G^k(\mathbb H^d;E)$.
Then $\delta_l f\in\mathcal G^k(\mathbb H^d;E)$
for every $l\in\mathbb N$, and
\[
 \|\delta_l f\|_{\mathcal G^k}
 \leq e^{k/l}\|f\|_{\mathcal G^k}.
\]
In particular, the operators $\delta_l$ are
uniformly bounded with respect to this norm.
\end{lemma}

\begin{proof}
The smoothness and local support of $k_l$ imply
that $\delta_l f$ is smooth.
Moreover, if $k_l(x,y)\neq0$, then
$d(x,y)\leq l^{-1}$. The triangle inequality then  gives
$ d(x,o)\leq d(y,o)+l^{-1}$, so that,
using positivity and the normalization of the kernel,
we  obtain
\begin{align*}
 e^{k d(x,o)}\|\delta_l f(x)\|_E
 &\leq
 \int_{\mathbb H^d}k_l(x,y)
 e^{k d(x,o)}\|f(y)\|_E\,
          \mathrm d\operatorname{vol}(y)\\
 &\leq
 e^{k/l}\|f\|_{\mathcal G^k}
 \int_{\mathbb H^d}k_l(x,y)\,
          \mathrm d\operatorname{vol}(y)\\
 &=e^{k/l}\|f\|_{\mathcal G^k}, 
\end{align*}
so that taking the supremum over $x$ proves the assertion.
\end{proof}

\begin{definition}\label{approx}
    Let $P_n := \one_{[-n,n]} (A_0)$ for $n \in \mathbb{N}$ and denote by 
    $(\phi_n)_{n\in \mathbb{N}} \subset \mathcal{D}$ an orthonormal basis of $H$ with $\mathcal{D} \subset H$ being a dense subset of $H$. Write $K_n$ for the finite rank
    projections onto $\text{span}\big\{(\phi_j)_{j=1}^n\big\} \subset H$
    and identify the range of $K_n$ with
    $\mathbb{C}^n$. For a fixed choice of
     delta sequences, 
   denote by $\delta_l$ the family of functionals from Definition \ref{delta}.
   For $B$ satisfying Hypothesis \ref{hypothesis}, define
   \begin{equation}
       B_n(x) := P_n B(x) P_n, \quad 
       B_{n,m}(x) := K_m B_n(x) K_m, 
       \quad B_{n,m,l} := \delta_l(B_{n,m}).
   \end{equation}

\end{definition}

\begin{proposition}\label{approximationscheme}
    With the notations from Definition 
    \ref{approx}, there exist three
    constants $C, C_n$ and $C_{n,m}$, where 
    $C$ is uniform, $C_n$ is independent of $m,l$ and $C_{m,n}$ is independent of $l$, such that
    \begin{equation}\label{limit1}
    \begin{split}
 & \lim_{n \rightarrow \infty} \textnormal{tr}_{S} \bigg( 
  R_z\big(D^*_{B_n}D_{B_n}\big)^N - R_z\big(D_{B_n}D^*_{B_n}\big)^N \bigg) \\
  &= 
  \textnormal{tr}_{S} \bigg( 
  R_z(D^*_{B}D_{B})^N - R_z(D_{B}D^*_{B})^N \bigg),
  \end{split}
    \end{equation}
for $z \in \mathcal{D}_C(B)$, 
 \begin{equation}\label{limit2}
 \begin{split}
  &\lim_{m \rightarrow \infty} \textnormal{tr}_{S} \bigg( 
  R_z\bigg(D^*_{B_{n,m}}D_{B_{n,m}}\bigg)^N - R_z\bigg(D_{B_{n,m}}D^*_{B_{n,m}}\bigg)^N \bigg) \\
  &= 
  \textnormal{tr}_{S} \bigg( 
  R_z\big(D^*_{B_n}D_{B_n}\big)^N - R_z(D_{B_n}D^*_{B_n})^N \bigg),
  \end{split}
    \end{equation}
for $\textnormal{dist}(z, (-\infty, 0]) \geq C_n$, and 
 \begin{equation}\label{limit3}
 \begin{split}
  &\lim_{l \rightarrow \infty} \textnormal{tr}_{S} \bigg( 
  R_z\bigg(D^*_{B_{n,m,l}}D_{B_{n,m,l}}\bigg)^N - R_z\bigg(D_{B_{n,m,l}}D^*_{B_{n,m,l}}\bigg)^N \bigg) \\
  &= 
  \textnormal{tr}_{S} \bigg( 
  R_z\bigg(D^*_{B_{n,m}}D_{B_{n,m}}\bigg)^N - R_z\bigg(D_{B_{n,m}}D^*_{B_{n,m}}\bigg)^N \bigg),
  \end{split}
    \end{equation}
for $\textnormal{dist}(z, (-\infty, 0]) \geq C_{n,m}$, are limits in $S^1(L^2(S,H))$.
\end{proposition}
\begin{proof}
We start by showing Equation \eqref{limit1}.
After decomposing both sides as 
described in Proposition \ref{decomposition2}, we 
obtain a product of operators featuring 
$U_z^-(B_n)$, $V^-(B_n)$, $ic( \Id B_n)
\langle M_0 \rangle^{-1}$ and $\langle M_0 \rangle \big(
R_z(D^*_{B_n}D_{B_n})^k + R_z(D_{B_n}D^*_{B_n})^k \big)$,
where the various upper indices are independent of 
$n$ and $B$, just as for the counterparts 
of each factor that involves $B$ instead of $B_n$.
Now, given that $P_n$ commutes with $A_0$ and has norm
$1$, we have 
\begin{equation*}
    \sup_{n \in \mathbb{N}} (\rho_z^{-2N}(B_n) + \rho_z^{2N}(B_n))
    \leq \rho_z^{-2N}(B) + \rho_z^{2N}(B) 
    \rightarrow 0, \quad \text{dist}(z, (-\infty, 0]) 
    \rightarrow \infty.
\end{equation*}
On the other hand, $P_n \rightarrow \one_H$, $n \rightarrow \infty$ strongly and for $\beta \in [-2N, 2N]$, and a.e. $x \in \mathbb{H}^d$,
\begin{equation*}
    \langle A_0 \rangle^{\beta} (B_n(x) + A_0)^2 \langle A_0 \rangle^{-2 - \beta} 
\xrightarrow[]{n \rightarrow\infty}
 \langle A_0 \rangle^{\beta} (B(x) + A_0)^2 \langle A_0 \rangle^{-2 - \beta},
\end{equation*}
and 
\begin{equation*}
     \langle A_0 \rangle^{\beta} \partial_i B_n(x) \langle A_0 \rangle^{-2- \beta}
 \xrightarrow[]{n \rightarrow\infty}    
  \langle A_0 \rangle^{\beta} \partial_i B(x) \langle A_0 \rangle^{-2- \beta}, 
  \quad i \in \{1, \dots, d \},
\end{equation*}
converge
in the strong operator topology of $H$.
Thus we can apply Lemma \ref{convergence1} 
to get SOTA convergence of $\overline{U_z^-(B_n)}$
to $\overline{U_z^-(B)}$, as well as 
$\langle M_0 \rangle \big(
R_z(D^*_{B_n}D_{B_n})^k + R_z(D_{B_n}D^*_{B_n})^k \big)$ to $\langle M_0 \rangle \big(
R_z(D^*_{B}D_{B})^k + R_z(D_{B}D^*_{B})^k \big)$.
The same argument then also implies SOTA convergence of
$i(c \Id B_n) \langle M_0 \rangle^{-1}$ 
to $i(c \Id B) \langle M_0 \rangle^{-1}$.
Lastly, we discuss the terms 
$\overline{V^-(B_n)}$ which are, contrary to 
$\overline{U_z^-(B)}$, $S^d(L^2(S,H))$- or 
$S^{d+1}(L^2(S,H))$-valued.
Since $A_0$ does not depend on the spatial variable,
$P_n$ is constant on $\mathbb{H}^d$ and given that it
commutes with $A_0$, we can write 
$V^-(B_n) = \tilde{P_n} V^-(B) \tilde{P_n}$ 
where $\tilde{P_n} := \one_{L^2(\mathbb{H}^d,S)} 
\otimes P_n \rightarrow \one_{L^2(S,H)}$
in SOTA as $n \rightarrow \infty$.
Hence, Lemma \ref{sotalemma} guarantees that 
$\overline{V^-(B_n)}$ converges in 
$S^d(L^2(S,H))$, respectively in 
$S^{d+1}(L^2(S,H))$, to $\overline{V^-(B)}$.
Another application of Lemma \ref{sotalemma} 
lets us then conclude that \eqref{limit1}
converges in trace-class for $z \in \mathcal{D}_C(B)$ and some global constant $C$. \\
Let us now fix $n \in \mathbb{N}$ and study
Equation \eqref{limit2}. Since $K_m$, as 
opposed to $P_n$, does not necessarily 
commute with $A_0$, we can not repeat the 
argument we just used to justify trace-class convergence 
of \eqref{limit1}.
However, given that for all $s \in \mathbb{R}$ and  
$n \in \mathbb{N}$, $\langle M_0 \rangle^s P_n$
is a bounded operator, we have
\begin{equation*}
    \tau^{0,0,\nabla,d,1}(B_n) < \infty.
\end{equation*}
Hence we can set $\alpha = \beta = 0$ in the 
decomposition of Proposition \ref{decomposition2},
which also implies that no factors involving 
$\langle M_0 \rangle$ appear anymore.
We will denote this choice of indices by $"\circ"$.
Since
\begin{equation*}
    \tau^{0,0,\nabla,d,1}(B_{n,m}) \leq C < \infty,
\end{equation*}
for an $m$-independent constant $C$, the
decomposition as discussed above can also be 
employed for $B_{n,m}$.
Now, given that $K_m$ is constant on $\mathbb{H}^d$, 
we have that $V^{\circ}(B_{n,m}) = \widetilde{K_m} V^\circ(B_{n,m})\widetilde{K_m}$, where
$\widetilde{K_m} := \one_{L^2(\mathbb{H}^d,S)} \otimes K_m$ is self-adjoint and converges in 
the strong operator topology to $\one_{L^2(S,H)}$.
Hence by Lemma \ref{sotalemma} and Lemma
\ref{extension2}, $V^{\circ}(B_{n,m})$ converges in the $S^d(L^2(S,H))$-, respectively 
$S^{d+1}(L^2(S,H))$-norm, to $V^{\circ}(B_{n})$
as $m \rightarrow \infty$. 
On the other hand, $\norm{B_{n,m}}_{C^1_b(\mathbb{H}^d,B(H))}$
and $\norm{B_n}_{C^1_b(\mathbb{H}^d,B(H))}$
are uniformly bounded in $m$ so that by Lemma
\ref{sota2}, $\overline{U^\circ_z(B_{n,m})}$
converges in strong operator topology 
to $\overline{U^\circ_z(B_{n})}$, where the factors 
of $\overline{U^\circ_z(B_{n,m})}$ not involving $c \nabla$ even converge in SOTA.
By a similar argument, $R_z(D^*_{B_{n,m}}D_{B_{n,m}})^k + R_z(D_{B_{n,m}}D^*_{B_{n,m}})^k$ converges in SOTA
to $R_z(D^*_{B_n}D_{B_n})^k + R_z(D_{B_n}D^*_{B_n})^k$.
Since $\norm{B_{n,m}}_{C^1_b(\mathbb{H}^d,B(H))}$ is
uniformly bounded in $m$, Lemma \ref{sota2}
guarantees that convergence holds whenever
$w := \text{dist}(z, (-\infty,0]) \geq C$, where 
$C$ is an $m$-independent constant.
Note moreover that also $i(c \Id B_{n,m})$
converges in SOTA to $i(c\Id B_n)$ so that we can conclude by Lemma \ref{sotalemma} that Equation \eqref{limit2} 
converges in trace-class whenever $w \geq C$. \\
Finally, we show Equation \eqref{limit3}.
Since $B_{n,m}(x) = K_m B_{n,m}(x) K_m$, we have
$B_{n,m,l}(x) = K_m B_{n,m,l}(x) K_m$ and 
$\Id B_{n,m,l}(x) = K_m \Id B_{n,m,l}(x) K_m$.
So in particular, with $\text{rang}(K_m) \simeq 
\mathbb{C}^m$, we have
\begin{equation*}
    \tau^{0,0,\nabla,d,1} \leq \norm{B_{n,m,l}}_{\mathcal{G}^k} \leq C < \infty.
\end{equation*}
By Lemma \ref{lemma4.3}, we know that 
$B_{n,m,l}$ and $\Id B_{n,m,l}$ converge point-wise to $B_{n,m}$, respectively $\Id B_{n,m}$.
Posing $\alpha = \beta = 0$ in the decomposition 
from Lemma \ref{decomposition2} as in the previous
step, Lemma \ref{sota2} guarantees that 
$\overline{U_z^\circ(B_{n,m,l} )}$ converges
in strong operator topology to 
$\overline{U_z^\circ(B_{n,m} )}$, where the factor
of $\overline{U_z^\circ(B_{n,m,l} )}$ not involving 
$c \nabla$ even converges in SOTA.
By a similar argument, 
$R_z(D^*_{B_{n,m,l}}D_{B_{n,m,l}})^k + R_z(D_{B_{n,m,l}}D^*_{B_{n,m,l}})^k$ and 
$i(c \Id B_{n,m,l})$ also converge in SOTA to 
$R_z(D^*_{B_{n,m}}D_{B_{n,m}})^k + R_z(D_{B_{n,m}}D^*_{B_{n,m}})^k$, respectively 
$i(c \Id B_{n,m})$.
Since by Lemma \ref{lemma4.3}, 
$\norm{B_{n,m,l}}_{C^1_b(\mathbb{H}^d,B(H))} 
\leq \norm{B_{n,m,l}}_{\mathcal{G}^k}$
is uniformly bounded in $l$, Lemma \ref{lemmabound}
guarantees that the domain of convergence is
$\text{dist}(z,(-\infty,0]) \geq C$,
for some $l$-independent constant $C$.
Given that $K_m$ is finite-rank and that $\Id B_{n,m,l}$ can be bounded uniformly in $l$ in $S^d$-norm
by multiples of $e^{-X}$, the dominated convergence 
theorem and Lemma \ref{fourier} 
yield that the terms 
$V^\circ(B_{n,m,l})$ converge in 
$S^d(L^2(S,H))$-, respectively $S^{d+1}(L^2(S,H))$-norm 
to $V^\circ(B_{n,m})$. 
Hence we conclude with Lemma \ref{sotalemma}
that \eqref{limit3} 
converges in trace-class for 
$w \geq C$ and some $l$-independent constant $C$.

\end{proof}

\begin{corollary}\label{cor:semigrouplimit}
    For $B$ satisfying Hypothesis \ref{hypothesis}
    and $t >0$, we have that
\begin{equation*}
    \lim_{n \rightarrow \infty} \lim_{m \rightarrow \infty} \lim_{l \rightarrow \infty} 
    \textnormal{tr}_{S}
    \left(
    e^{-tD^*_{B_{n,m,l}} D_{B_{n,m,l}}}
    - e^{-tD_{B_{n,m,l}}D^*_{B_{n,m,l}}}
    \right) = \textnormal{tr}_{S}
    \left(
    e^{-tD^*_{B} D_{B}}
    - e^{-tD_{B}D^*_{B}}
    \right),
\end{equation*}
are limits in $S^1(L^2(\mathbb{H}^d,H))$.
\end{corollary}
\begin{proof}
This follows from holomorphic functional calculus
and the fact that  the integration 
curve can be chosen independently of the parameter over which  the limit is being taken, due to the
nature of the domains of definition for the parameter $z$.
\end{proof}

\section{Construction of the integral kernels}\label{section5}
\subsection{Auxiliary kernels and regularity}\label{subsection5.1}
Roughly speaking, the overarching
strategy for computing $$\textnormal{tr}_{L^2(\mathbb{H}^d,H)}
    \textnormal{tr}_{S} \left(
    e^{-tD^*D} - e^{-tDD^*}\right),$$
will be to find an integral kernel for 
$\textnormal{tr}_{S} \left(
    e^{-tD^*D} - e^{-tDD^*}\right)$ and 
    then use a Mercer-type theorem for 
    manifolds, c.f. Theorem \ref{mercer}. 
    To achieve the first, we 
    decompose the internal trace using Volterra series arguments and construct secondary integral 
    kernels
    (Subsection \ref{subsection5.1}), which 
combined will yield the main integral kernel
needed for the calculations performed in 
Section \ref{section6} and, in particular, 
Theorem \ref{maintheorem}.\\
We now work within the approximation framework developed in Section \ref{section4}, and thus consider the family $B_{n,m,l}$.
In order to deduce the corresponding statements for the limiting object $B$,
 Theorem \ref{maintheorem} requires passing to the  iterated limit $n,m,l \rightarrow \infty$.
Throughout the intermediate arguments, and in order to avoid cumbersome notation, we, however, suppress the explicit dependence on the parameters and simply write 
$B$
in place of 
$B_{n,m,l}$.	
The indices  $n,m,l$
 will be reinstated later, whenever their role becomes relevant to the analysis or to the limiting procedure.

\begin{definition}
    Let $l,s \in \mathbb{N}$ and $p \in [1, \infty]$. For a Banach space $X$, we denote by
    $C^{l,p}(\mathbb{H}^s,X)$ the space of 
    $l$-times continuously differentiable 
    functions with derivatives in  Bochner-$L^p$ and 
    values in $X$. Moreover, 
    for $l,m,n,s \in \mathbb{N}_0$
    with $m \geq n$, 
    let 
    \begin{equation}
    C^l_{m,n}(\mathbb{H}^s) :=
    C^{l,\infty}(\mathbb{H}^s, B(S_x) \otimes 
    \langle A_0 \rangle^{-m} B(H) \langle A_0 
    \rangle^n),
    \end{equation}
    be the Banach space with norm
\begin{equation}
    \norm{f}_{l,m,n} := \sum_{\alpha \in \mathbb{N}_0^s, \vert \alpha \vert \leq l} 
    \norm{f}_{L^{\infty}(\mathbb{H}^s,
    B(S_x) \otimes 
    \langle A_0 \rangle^{-m} B(H) \langle A_0 
    \rangle^n)}, \quad 
    f \in C^l_{m,n}(\mathbb{H}^s).
\end{equation}
Similarly, denote by 
\begin{equation}
C^{l,S^1}_{m,n}(\mathbb{H}^s)
:= C^{l,p}(\mathbb{H}^s, B(S_x) \otimes 
    \langle A_0 \rangle^{-m} S^1(H) \langle A_0 
    \rangle^n),
    \end{equation}
the Banach space with norm
\begin{equation}
 \norm{f}_{l,m,n}^{S^1} := \sum_{\alpha \in \mathbb{N}_0^s, \vert \alpha \vert \leq l} 
    \norm{f}_{L^{\infty}(\mathbb{H}^s,
    B(S_x) \otimes 
    \langle A_0 \rangle^{-m} S^1(H) \langle A_0 
    \rangle^n)}, \quad 
    f \in C^{l,S^1}_{m,n}(\mathbb{H}^s).   
\end{equation}

\end{definition}
\begin{rem}\label{rem:hypothesisapprox}
    Let us assume that $B$ satisfies Hypothesis
    \ref{hypothesis}. Then, for instance, 
    Equation \eqref{hyp2} says that
    for $\phi \in \text{Dom} \: A_0$, a.e. $y \in S_1(o)$,
a.e. $x \in \mathbb{H}^d$, and $\gamma \in 
\mathbb{N}_0^d$ with 
$\vert \gamma \vert \leq 1$, we have
\begin{equation}
    \lim_{R \rightarrow \infty} e^{R + \sqrt{R}}\norm{(
    \partial^{\gamma} B(w) -
    \tau^{*,LC}(w,w')
    \partial^{\gamma}B(w')) \phi}_H = 0,
\end{equation} 
where  $w = \exp_{w'}(y)$ and $
w' =  \exp_o(Ry)$.
It then follows, for  a.e. $y \in S_1(o)$,
a.e. $x \in \mathbb{H}^d$, and $\gamma \in 
\mathbb{N}_0^d$, that
\begin{equation}\label{hyp2approx}
    \lim_{R \rightarrow \infty} e^{R + \sqrt{R}}\norm{
    \partial^{\gamma} B_{n,m,l}(w) -
    \tau^{*,LC}(w,w')
    \partial^{\gamma}B_{n,m,l}(w')}_{\langle A_0 \rangle^{-i} S^1(H) \langle A_0 \rangle^{j}} = 0,
\end{equation} 
holds as well, 
for all $n,m,l \in \mathbb{N}$, $i,j \in \mathbb{N}_0$.
One can see this by noting that $P_n$ commutes with
$A_0$, and $K_m$ is finite rank, hence mapping into
any $\text{Dom}A_0^{i}$, $i \in \mathbb{N}_0$,
making \eqref{hyp2approx} hold
for $B_{n,m}$. That it holds 
for  $B_{n,m,l}$  follows
from $B_{n,m,l} = \delta_l(B_{n,m})$,
$\norm{B_{n,m}}_{\mathcal{G}^k} < \infty$
for some $k >0$, and the dominated convergence
Theorem. One then argues similarly 
for the remaining statements of 
Hypothesis \ref{hypothesis}. See also Remark 3.6 in \cite{oliver23}. \\
Note that we also have, for $n,m,l \in \mathbb{N}$ and
$i,j \in \mathbb{N}_0$,
\begin{equation}
    (B_{n,m,l} + M_0)^2 -M_0^2, \quad 
    ic(\Id B_{n,m,l}) \in C_{i,j}^{\infty, S^1}(\mathbb{H}^d).
\end{equation}
This also follows from  $P_n$ commuting with
$A_0$ and $K_m$ being finite rank, hence mapping into
any $\text{Dom}A_0^{i}$, $i \in \mathbb{N}_0$,
and the smoothness of 
$\delta_l(f)$ for any $f \in \mathcal{G}^k$.

\end{rem}

\begin{definition}
The Euclidean $l$-simplex
$\Delta_l \subset [0,1]^l$ is the set 
\begin{equation}
    \Delta_l := \{ w \in [0,1]^l \: \vert \: w_l \leq w_{l-1} \leq \dots \leq w_1 \}.
\end{equation}
A frequently used parametrization is 
realized with $l +1$ coordinates 
by setting 
\begin{equation}
    s_0 := 1-w_1, \: s_1 := w_1 - w_2, \dots, s_{l-1} := w_{l-1} -w_l, \: s_l := w_l,
\end{equation}
such that $s_j \in [0,1]$ with 
$\sum_{j =0}^l s_j = 1$.
\end{definition}

\begin{rem}
    The volume of $\Delta_l$ is 
    $\text{vol}(\Delta_l) = \frac{1}{l!}$.
\end{rem}

\begin{definition}
Let $Y$ be a Hilbert space.
    The convolution shuffle 
    product $q \shuffle_* M 
    \in L^0( \mathbb{H}^d \times \mathbb{H}^d, B(Y))$ 
    pairing $q = (q_0, \dots, q_l)
    \in (L^0( \mathbb{H}^d \times \mathbb{H}^d, B(Y))^{l+1}$ with 
    $M = (M_1, \dots, M_l) \in 
    (L^0(\mathbb{H}^d, B(Y))^{l}$ 
    is defined by 
    \begin{equation}
(q \shuffle_* M)(x,y) :=  
\int_{u \in (\mathbb{H}^d)^l}
q_0(x,u_1)M_1(u_1)q_1(u_1,u_2) \dots 
M_l(u_l) q_l(u_l,y) \Id u,
    \end{equation}
whenever the integral converges absolutely 
in $B(Y)$-norm for a.e. $(x,y) \in \mathbb{H}^d \times \mathbb{H}^d$.
\end{definition}

\begin{theorem}[\cite{mercermfd}]\label{mercer}
    Let $H$ be a separable 
    Hilbert space and
    $T \in S^1(L^2(\mathbb{H}^d,H))$. 
   Then $T$ admits a unique, 
    $S^1(H)$-valued integral kernel 
    $t$, such that 
    \begin{equation}
        (Tf)(x) = \int_{\mathbb{H}^d} t(x,y) f(y) 
        \Id y,
    \end{equation}
    for $f \in L^2(\mathbb{H}^d, H)$.
    Moreover, there is a bounded, continuous map
    \begin{equation}
        \begin{split}
            \mathbb{H}^d &\longrightarrow L^1(\mathbb{H}^d, S^1(H)) \\
    h &\longmapsto   (x \mapsto t(x,L_hx)),     
        \end{split}
    \end{equation}
    where $L_h$ denotes 
    translation by $h \in \mathbb{H}^d$,
    and 
    \begin{equation}
    \begin{split}
        \int_{\mathbb{H}^d} \norm{t(x,x)}_{S^1(H)} 
        &\leq \norm{T}_{S^1(L^2(\mathbb{H}^d,H))}, \\
    \textnormal{tr}_{L^2(\mathbb{H}^d,H)}(T) &= 
    \int_{\mathbb{H}^d} \textnormal{tr}_H(t(x,x)) \Id x.
 \end{split}       
    \end{equation}
\end{theorem}

\begin{rem}\label{remmercer}
    In the context of the stated uniqueness
    in
    Theorem \ref{mercer},
    note that 
    the separability of the Hilbert space
    $H$ and the variational principle guarantee that any jointly continuous,
    $B(H)$-valued integral kernel 
    $\tilde{t}$ of $T$ already satisfies
    $t = \tilde{t}$ as an element
    of $C_b(\mathbb{H}^d, L^1\left(\mathbb{H}^d,S^1(H))\right)$.
\end{rem}

\begin{lemma}[\cite{oliver23}, Lemma 4.4] 
Let $T,S$ be self-adjoint operators 
in a Hilbert space $H$, such that 
$S$ is bounded and $T \geq 0$. Then for $t \geq 0$, 
\begin{equation}
    e^{-t(T+S)} = e^{-tT} + \sum_{l =1}^{\infty}
    (-1)^l \int_{s \in t \Delta_l} e^{-sT} \shuffle 
    S^{\otimes l} \Id s.
\end{equation} 
\end{lemma}

\begin{definition}
    For $t >0$ and $k \in \mathbb{N}$, 
    define 
    \begin{equation}
    R^k_t :=
        \begin{dcases*}
e^{-t H_B}
   & if  $k = 0$\,, \\[1ex]
\int_{s \in t \Delta_k} e^{-s H_B} \shuffle 
\left(ic \Id B \right)^{\otimes k} \Id s
   & if $k \geq 1$\,.
\end{dcases*}
    \end{equation}
\end{definition}

\begin{lemma}
    For $t > 0$, 
    \begin{equation}
        e^{-t D D^*} = \sum_{l = 0}^{\infty}
        (-1)^l R^l_t 
        \quad \text{and} \quad 
e^{-t D^* D} = \sum_{l = 0}^{\infty}
         R^l_t         
    \end{equation}
    hold, with absolute convergence in 
    operator norm.
\end{lemma}

\begin{definition}
    Let $t >0$, $l \in \mathbb{N}$. 
    Set
    \begin{align*}
K^{-,l}_t&:= \sum_{j = \frac{l-1}{2}}^{\infty} R^{2j+1}_t, \: l \: \text{odd},   &  K^+_t &:= \sum_{j =\frac{l}{2}}^{\infty} R^{2j}_t, \: l \: \text{even},               \\
k^{-,l}_t&:= \sum_{j =0}^{\frac{l-3}{2}}   R^{2j+1}_t, \: l \: \text{odd},        &  k^+_t&:= \sum_{j = 1}^{\frac{l-2}{2}}   R^{2j}_t, \: l \: \text{even}.  
\end{align*}
\end{definition}

\begin{lemma}\label{commutator}
Let $t >0$. Then, on $\textnormal{dom} H_0$,
\begin{equation}
    \begin{split}
   & [D,D^*e^{-sDD^*}]
        -[D^*,D e^{-sD^*D}] \\
        & \quad + [D, e^{-sD^*D}D^*] - [D^*, e^{-sDD^*}D] \\
& =  [ic\nabla, ic\nabla (e^{-sD^*D} - e^{-sDD^*})]
+ [ic\nabla, (e^{-sDD^*} - e^{-sD^*D}) ic\nabla] \\
&\quad - [A, ic\nabla (e^{-sD^*D }+ e^{-sDD^*})] 
- [A, (e^{-sDD^*} + e^{-sD^*D})ic\nabla] \\
&= [ic\nabla, [ic\nabla, e^{-sD^*D}-e^{-sDD^*}]] 
- [A, \{ic\nabla, e^{-sD^*D} + e^{-sDD^*} \}].
    \end{split}
\end{equation}
\end{lemma}
\begin{proof}
    Recall that $D = ic \nabla + A$ and $D^* = ic\nabla - A$, then expand the commutators.
\end{proof}

\begin{corollary}
      Let $t>0$. Then 
    \begin{equation}\label{eq5.15}
    \begin{split}
        e^{-tD^*D} - e^{-tDD^*} 
        &= \int_0^t
        \overline{[ic\nabla, ic \nabla(
        K^{-,d-2}_s + k^{-,d-2}_s)]} 
        - \overline{[A, A(K^{-,d}_s +k^{-,d}_s)]} \\
 &+       \overline{[ic\nabla, A(K^{+,d-1}_s + k^{+,d-1}_s +e^{-sH_B})]} \\
 &- \overline{[A, ic\nabla(K^{+,d-1}_s+ k^{+,d-1}_s + e^{-sH_B})]} \Id s,
        \end{split}
    \end{equation}
    with the integrals converging strongly on 
    $\text{dom}(H_0)$.
\end{corollary}
  
\begin{proposition}[\cite{camporesi92}, Theorem 3.1.]\label{prop:spinorheatkernel}
    For $d$ odd, the spinor
      heat kernel $P^{CP}_t$ on $\mathbb{H}^d$ 
      is given by, for $x,y \in \mathbb{H}^d$,
\begin{equation*}
    P^{CP}_t(x,y) := p^{CP}_t(\textnormal{dist}(x,y)) \cdot \tau(x,y),
\end{equation*}
where  
\begin{equation}\label{eq:scalerheatkernelpart}
    \begin{split}
&p^{CP}_t(s) \\
&= 2^{\frac{5-3d}{2}} \sqrt{\pi} \: \Gamma\left(\frac{d}{2}\right)^{-1}
(4 \pi t)^{-1/2} \cosh \left(\frac{s}{2}\right)
\left(-\frac{d}{d(\cosh{s})}\right)^{\frac{d-1}{2}} \left(\cosh \left(\frac{s}{2}\right) \right)^{-1} e^{-\frac{s^2}{4t}}, \\ 
\end{split}
    \end{equation}
    for all $s \in \mathbb{R}$
    and $\tau(x,y) \: : \: S_y \rightarrow S_x$
    denotes spinor parallel transport.
\end{proposition}

\begin{rem}
    We write $P^{CP}_t$ for the hyperbolic-spinor heat kernel with the Camporesi-Pedon normalization, see Remark \ref{rem:haarnormalization}, and 
    $P_t$ for the rescaled spinor heat kernel giving the usual small-time asymptotics, see Remark \ref{rem:smalltimespinorkernel}.
\end{rem}

\begin{rem}\label{taulc}
    Let $x,y \in \mathbb{H}^d$ and consider two
    sections $\xi_y \in T_y^* \mathbb{H}^d$
    and $s_y \in S_y$. We want to argue that
    \begin{equation*}
        \tau(x,y) \big(c(\xi_y) s_y \big)
        = c \big( \tau^{*,LC}(x,y) \xi_y \big)
        \tau(x,y)s_y,
    \end{equation*}
by showing that both sides of the Equation satisfy the same differential equation, 
where $\tau^{LC}(x,y) \: : \: T_y \mathbb{H}^d 
\rightarrow T_x \mathbb{H}^d$ denotes the Levi-Civita parallel transport and $\tau^{*,LC}$
the induced parallel transport on the 
cotangent bundle.
Let $f(t) := \tau(\gamma(t),y) \big(c(\xi_y) s_y \big) \in S_{\gamma(t)}$, for the unique  geodesic $\gamma \: : \: [0,1] \rightarrow \mathbb{H}^d$ such that $\gamma(0) =y$, $\gamma(1)= x$. We then have
$$\nabla_{\dot \gamma(t)}^S f(t) = 0,
\quad  f(y) = c(\xi_y) s_y.$$
On the other hand,
\begin{equation*}
\begin{split}
    &\nabla_{\dot \gamma(t)}^S \big[ 
    c \big( \tau^{*,LC}(\gamma(t),y) \xi_y \big)
        \tau(\gamma(t),y)s_y\big]\\
        &= c \big(\nabla_{\dot \gamma(t)}^{*,LC}
        \tau^{*,LC}(\gamma(t),y) \xi_y \big) \cdot\tau(\gamma(t),y) s_y 
        + c \big( \tau^{*,LC}(\gamma(t),y) \xi_y \big) 
        \cdot \nabla_{\dot \gamma(t)}^S\tau(\gamma(t),y) s_y \\
        &= 0,
        \end{split} 
\end{equation*}
satisfying the same initial condition
as $f$, which
lets us conclude.

\end{rem}

\begin{lemma}
    Let $t >0$. Then
 \begin{equation}\label{eq5.16}
    \begin{split}
  \textnormal{tr}_{S}\left(     e^{-tD^*D} - e^{-tDD^*}\right)
        &= \int_0^t  
        \textnormal{tr}_{S} \left( \left[    
        ic \nabla, ic \nabla K_s^{-,d-2}
        + AK_s^{+,d-1}
        \right] \right)
        \\
        &- \left[A, \textnormal{tr}_{S} \left(
        AK_s^{-,d} +ic \nabla K_s^{+,d-1}
         \right)\right] \Id s. 
        \end{split}
    \end{equation}    
\end{lemma}
\begin{proof}
Note that $\text{tr}_{\mathbb{C}^r}$ is  continuous 
from $B(L^2(S,H))$ to $B(L^2(\mathbb{H}^d,H))$,
and that regardless of whether the spaces 
are equipped with the weak operator, the strong operator or the norm topology. 
In Equation \eqref{eq5.15}, restrict the 
integrand to $\text{Dom} H_0 = S_x \otimes 
(H^2(\mathbb{H}^d,H) \cap L^2(\mathbb{H}^d, 
\text{Dom} A_0^2))$ which allows us, given that the integral
converges strongly 
on $\text{Dom} H_0$, to commute $\text{tr}_{\mathbb{C}^r}$ with the integration. \\
    We have to argue that the lower terms vanish.
    For that, observe that pointwise,  
\begin{equation*}
\begin{split}
    & \text{tr}_{S} \big(
    P_{s_0}(x,y_1) ic (\Id B)(y_1) 
    \dots ic(\Id B)(y_{d-1})
    P_{s_{d-1}}(y_{d-1},x) \big) \\
    =& 
    p_{s_0}(\text{dist}(x,y_1)) \dots 
    p_{s_{d-1}}(\text{dist}(y_{d-1},x)) \\
    &  \text{tr}_{S}\bigg(  
    \tau(x,y_1) c(\Id B)(y_1) 
    \dots \tau(y_{d-2}, y_{d-1}) 
    c(\Id B)(y_{d-1}) \tau(y_{d-1},x) \bigg)  \\
=& 
    p_{s_0}(\text{dist}(x,y_1)) \dots 
    p_{s_{d-1}}(\text{dist}(y_{d-1},x)) \\
    &  \text{tr}_{S}\bigg( 
   c\big( \tau^{*,LC}(x,y_1) \Id B \big)(y_1) 
    \dots c\big( \tau^{*,LC}(x,y_1) \dots
    \tau^{*,LC
    }(y_{d-2}, y_{d-1}) 
    \Id B\big)(y_{d-1})\\
    &\tau(x,y_1) \dots \tau(y_{d-2},y_{d-1})\tau(y_{d-1},x) \bigg),\\   
    \end{split}
\end{equation*}
where the second equality follows from Remark
\ref{taulc}.
Hence Proposition \ref{cliffordid} implies 
that
\eqref{eq5.16} holds on 
$H^2(\mathbb{H}^d,H) \cap L^2(\mathbb{H}^d, 
\text{Dom} A_0^2)$.
On the other hand, for $k \geq d-2 \geq 1$, $U,W \in \{ic \nabla, A \}$, 
$c := \norm{ic( \Id B)}_{B(L^2(S,H))}$
and some constant $C$, the estimates 
$\norm{\overline{UWR^k_s}}_{B(L^2(S,H))} \leq C \frac{c^{k-1}s^{k-1}}{(k-1)!}$ and $\norm{\overline{UR^k_sW}}_{B(L^2(S,H))} \leq C \frac{c^{k-1}s^{k-1}}{(k-1)!}$,
which are obtained by exploiting the fact that 
$i c (\Id B)$ is a $C^1_b$-family of finite rank
operators with bounded derivatives and range in 
$\text{Dom}A_0$, imply 
\begin{equation*}
    \norm{\overline{UWK^k_s}}_{B(L^2(S,H))}, \norm{\overline{UK^k_sW}}_{B(L^2(S,H))}
    \leq C \sum_{k = j-1}^{\infty} \frac{c^k s^k}{k!}.
\end{equation*}
Since this is $O(s^{j-1})$ for $s \rightarrow 0$,
we have that $\overline{UWK^k_s}$ and 
$\overline{UK^k_sW}$ are Bochner integrable
in operator norm for $s \in (0,t)$ and any
$t>0$. 
Thus, by boundedness of $e^{-tD^*D} - e^{-tDD^*}$,
continuity of $\text{tr}_{\mathbb{C}^r}$
and density of $H^2(\mathbb{H}^d,H) \cap L^2(\mathbb{H}^d, 
\text{Dom} A_0^2)$ in $L^2(\mathbb{H}^d,H)$, 
we conclude that Equation \eqref{eq5.16} holds in 
$L^2(\mathbb{H}^d,H)$ and that the integral 
converges in operator norm.

\end{proof}

\begin{proposition}[\cite{grig}, Theorem 1.1.]
    For $d = 2m +1$ odd, the scalar heat kernel $q_t$
    on $\mathbb{H}^d$ is given by, 
    for $s = \text{dist}_{\mathbb{H}^d}(x,y)$, 
   \begin{equation*}
       q_t(s) = 
       \frac{(-1)^m}{2^m \pi^m}
       \frac{1}{(4 \pi t)^{1/2}} 
       \left( \frac{1}{\sinh s} \frac{\partial}{\partial s}
       \right)^m e^{-m^2t - \frac{s^2}{4t}}.
   \end{equation*}
\end{proposition}

\begin{rem}\label{heatestimate}
By Proposition \ref{casimir}, the spinor Laplacian satisfies 
$\slashed{\Delta} = \nabla^* \nabla + \frac{\text{scal}(\mathbb{H}^d)}{4}$, where $\nabla^* \nabla$ denotes the
Bochner Laplacian and $\text{scal}(\mathbb{H}^d) = -d(d-1)$. Hence 
$$e^{-t\slashed{\Delta}} = e^{-t \nabla^* \nabla} 
e^{-t \frac{\text{scal}(\mathbb{H}^d)}{4}}.$$
Since the Clifford representation is unitary,
the Kato-Simon inequality (see 
\cite{batu}, Theorem VII.8) implies 
that the heat kernel of the connection Laplacian 
on spinors is dominated by the scalar heat kernel 
$q_t(x,y)$ on $\mathbb{H}^d$,
$$\norm{{K_{\nabla \nabla^*}(x,y)}}_{\text{op}}
\leq  
q_t(x,y), \quad t>0,$$
where $\norm{\cdot}_{\text{op}}$ is the operator
norm on the spinor fibers.
Therefore, 
\begin{equation*}
\norm{P_t(x,y)}_{\text{op}} 
\leq e^{t \frac{d(d-1)}{4}} q_t(x,y).
\end{equation*}

\end{rem}

\begin{proposition}
    On the diagonal, $\text{dist}_{\mathbb{H}^d}(x,y) = s = 0$, we have, for $d= 2m +1 \geq 3$,
    \begin{equation}
        \begin{split}
            p^{CP}_t(0)
=
&\frac{1}{2^{4m}\Gamma\left(m+\frac12\right)^2} \\
& \cdot\sum_{k=0}^{m}
\Gamma\left(k+\frac12\right)
e_{m-k}\left(
\frac14,\frac94,\ldots,\frac{(2m-1)^2}{4}
\right)
t^{-k-\frac12},
        \end{split}
    \end{equation}
    where $e_j$ denotes the $j$-th elementary symmetric polynomial,
with $e_0=1$.
In particular, the leading small-time behavior is
\begin{equation}\label{eq:pt-zero-asymptotic}
p^{CP}_t(0)
\sim
\frac{1}{2^{2d-2}\Gamma\left(\frac d2\right)}
t^{-d/2},
\end{equation}
for $ t\downarrow0$.
\end{proposition}
\begin{proof}
By combining 
Equation (6.10) with Equation (8.7) and the
normalization of spherical functions via 
(5.24)--(5.25) from Camporesi-Pedon \cite{camporesi}, we obtain 
\begin{equation}\label{eq:spectral-pt-zero}
p^{CP}_t(0)
=
\frac{2^{1-4m}}{\Gamma\left(m+\frac12\right)^2}
\int_0^\infty
e^{-t\lambda^2}
\prod_{j=0}^{m-1}
\left(
\lambda^2+\left(j+\frac12\right)^2
\right)
\,d\lambda,
\end{equation}
after rewriting  the coefficients 
appearing in Equation (6.10) in \cite{camporesi}
in terms of the $\Gamma$-function.
Now, the polynomial in~\eqref{eq:spectral-pt-zero} can be expanded in terms of
elementary symmetric polynomials such that
\begin{equation}\label{eq:polynomial-expansion}
\prod_{j=0}^{m-1}
\left(
\lambda^2+\left(j+\frac12\right)^2
\right)
=
\sum_{k=0}^{m}
e_{m-k}\left(
\frac14,\frac94,\ldots,\frac{(2m-1)^2}{4}
\right)
\lambda^{2k},
\end{equation}
which, together with the Gaussian integral,
\begin{equation}
\int_0^\infty
\lambda^{2k}e^{-t\lambda^2}\,d\lambda
=
\frac12
\Gamma\left(k+\frac12\right)
t^{-k-\frac12},
\end{equation}
lets us conclude
\end{proof}

\begin{table}[H]
    \centering
    \begin{tabular}{|c|c|}\hline
         $d$&  $p_t(0)$\\\hline
         $3$&  $\frac{1}{16\sqrt{\pi}}t^{-1/2} + \frac{1}{8\sqrt{\pi}}t^{-3/2}$\\\hline
         $5$&  $\frac{1}{256\sqrt{\pi}}t^{-1/2} + \frac{5}{576\sqrt{\pi}}t^{-3/2} + \frac{1}{192\sqrt{\pi}}t^{-5/2}$\\ \hline
 $7$&$\frac{1}{4096\sqrt{\pi}}t^{-1/2}+ \frac{259}{460800\sqrt{\pi}}t^{-3/2} + \frac{7}{15360\sqrt{\pi}}t^{-5/2} + \frac{1}{7680\sqrt{\pi}}t^{-7/2}$\\\hline
    \end{tabular}
    \caption{First three values of 
    $p^{CP}_t(0)$.}
    \label{tab:valuespt}
\end{table}

\begin{rem}\label{rem:smalltimespinorkernel}
    Employing the rescaling as described in 
    Remark \ref{rem:haarnormalization} and Equation \eqref{eq:pt-zero-asymptotic}, we recover the usual small-time 
    asymptotics 
    \begin{equation}
        P_t(x,x) \sim (4 \pi t)^{-d/2} \one_S, \qquad t \downarrow 0. 
    \end{equation}
\end{rem}

\begin{rem}\label{vanishingder}
Since the hyperbolic spinor heat kernel
only depends on the signed geodesic distance and
that dependence is an even smooth function, its derivatives with respect to the first 
and second variables vanish on the diagonal,
\begin{equation}
    \Id_1 P_t(x,y) \bigg\vert_{y=x} =
    \Id_2 P_t(x,y) \bigg\vert_{y=x} = 0.
\end{equation}

\end{rem}

\begin{definition}\label{kernels1}
    Let $t>0$ and $j \geq 1$. Define
     \begin{align*} 
P_t &:= p_t( \text{dist}(x,y)) \cdot \tau(x,y),          &              \\ 
Q_t &:= e^{-tA_0^2} \otimes P_t,      &  \theta_t^j&:= \int_{s \in t\Delta_j} Q_s \shuffle_* (A^2 - M_0^2)^{\otimes j} \Id s, \\
\theta_t &:= Q_t + \sum_{j = 1}^{\infty} (-1)^j\theta^j_t,  &  r^j_t&:= \int_{s \in t\Delta_j}  \theta_s \shuffle_* (ic \Id B)^{\otimes j} \Id s,      
\end{align*}
where $\tau(x,y)$ denotes parallel transport 
\end{definition}

We next examine the regularity of the kernels from Definition \ref{kernels1} and derive the associated estimates. The proofs follow classical arguments and reduce to inserting the bound 
   $\norm{\partial^{\beta}P_t   }_{L^1(\mathbb{H}^d)} \leq C_{\beta} e^{\frac{d(d-1)}{4} t} t^{-\frac{\vert \beta \vert}{2}}$
  into the corresponding steps of Lemmas 4.12, 4.13, 4.15, 4.18, and Corollaries 4.14, 4.16 in \cite{oliver23}, so we 
  do not repeat them here. 
  Note that the assumption that
  $ic (\Id B)$ is finite rank is crucial in what follows.

\begin{lemma}
    Let $l,l',m,m',n,n' \in \mathbb{N}_0$
    with $m' \geq m \geq n \geq n'$, and 
    $l' \geq l$. 
    For $t>0$, define two linear 
    operators via 
    \begin{equation}
        (Q^L_t f)(x) := \int_{\mathbb{H}^d}
        Q_t(x,y) f(y) \Id y, \quad 
  (Q^R_t f)(x) := \int_{\mathbb{H}^d}
          f(y) Q_t(x,y) \Id y,
    \end{equation}
for $f \in C^{l,(S^1)}_{m,n}(\mathbb{H}^d)$.
Then $Q^L_t$ and $Q^R_t$ are continuous 
as maps to $C^{l',(S^1)}_{m',n}$, respectively $C^{l',(S^1)}_{m,n'}$, and 
there is a $t$-independent constant 
$C$ such that 
\begin{equation}
    \norm{Q^L_t}_{B\left(C^{l,(S^1)}_{m,n},C^{l',(S^1)}_{m',n}\right)} \leq C e^{\frac{d(d-1)}{4}t}t^{-\frac{l'-l+m'-m}{2}}, 
\end{equation}
and 
\begin{equation}
    \norm{Q^R_t}_{B\left(C^{l,(S^1)}_{m,n},C^{l',(S^1)}_{m,n'}\right)} \leq C e^{\frac{d(d-1)}{4}t}t^{-\frac{l'-l+n-n'}{2}}.
\end{equation}
\end{lemma}

\begin{lemma}
    For $n \in \mathbb{N}_0$,
 $t>0$ and $j \in \mathbb{N}$, we have
    \begin{equation}
\int_{s \in t \Delta_j} \norm{ \left(Q_s \shuffle_* (A^2-M_0^2)^{\otimes j} \right)(x,y) }_{B(S \otimes \textnormal{Dom}(A_0^n) )} \Id s \leq \frac{(ct)^j}{j!} P_t(x,y),
    \end{equation}
    where $x,y \in \mathbb{H}^d$ and $c := \sup_{x \in \mathbb{H}^d} 
    \norm{\overline{A(x)^2-A_0^2}}$.
    Moreover, let $l,l',m,m',n,n' \in \mathbb{N}_0$
    with $m' \geq m \geq n \geq n'$, and 
    $l' \geq l$. 
    For $t>0$, define two linear 
    operators via 
    \begin{equation}
        \left(\Theta^{j,L}_t f\right)(x) := \int_{\mathbb{H}^d}
        \theta^j_t(x,y) f(y) \Id y, \quad 
  \left(\Theta^{j,R}_t f\right)(x) := \int_{\mathbb{H}^d}
          f(y) \theta^j_t(x,y) \Id y,
    \end{equation}
for $f \in C^{l,(S^1)}_{m,n}(\mathbb{H}^d)$.
Then $Q^L_t$ and $Q^R_t$ are continuous 
as maps to $C^{l',(S^1)}_{m',n}$, respectively $C^{l',(S^1)}_{m,n'}$, and 
there is a $t,j$-independent constant 
$C$ such that 
\begin{equation}
    \norm{\Theta^{j,L}_t}_{B\left(C^{l,(S^1)}_{m,n},C^{l',(S^1)}_{m',n}\right)}, \norm{\Theta^{j,R}_t}_{B\left(C^{l,(S^1)}_{m,n},C^{l',(S^1)}_{m,n'}\right)} \leq \frac{ e^{\frac{d(d-1)}{4}t(1+j)} 
    C^j t^{j-\frac{l'-l}{2}}}{\prod_{i=1}^j \max(1,i-\frac{l'-l}{2})}.
\end{equation}
\end{lemma}

\begin{corollary}\label{cor1}
    For $n \in \mathbb{N}_0$, there 
    is a constant $c$, such that 
    for $t>0$,
    \begin{equation}
        \sum_{j=0}^{\infty} \norm{\theta^j_t(x,y)}_{B(S \otimes \textnormal{Dom}(A_0^n) )} 
        \leq e^{ct} P_t(x,y),
    \end{equation}
    where $x,y \in \mathbb{H}^d$.
 Moreover, let $l,l',m,m',n,n' \in \mathbb{N}_0$
    with $m' \geq m \geq n \geq n'$, and 
    $l' \geq l$. 
    For $t>0$, define two linear 
    operators via 
    \begin{equation}
        \left(\Theta^{L}_t f\right)(x) := \int_{\mathbb{H}^d}
        \theta_t(x,y) f(y) \Id y, \quad 
  \left(\Theta^{R}_t f\right)(x) := \int_{\mathbb{H}^d}
          f(y) \theta_t(x,y) \Id y,
    \end{equation}
for $f \in C^{l,(S^1)}_{m,n}(\mathbb{H}^d)$.
Then $Q^L_t$ and $Q^R_t$ are continuous 
as maps to $C^{l',(S^1)}_{m',n}$, respectively $C^{l',(S^1)}_{m,n'}$, and 
there is a continuous function $h$ such that 
\begin{equation}
    \norm{\theta^L_t}_{B\left(C^{l,(S^1)}_{m,n},C^{l',(S^1)}_{m',n}\right)} \leq t^{-\frac{l'-l+m'-m}{2}} h(t), 
\end{equation}
and 
\begin{equation}
    \norm{\theta^R_t}_{B\left(C^{l,(S^1)}_{m,n},C^{l',(S^1)}_{m,n'}\right)} \leq t^{-\frac{l'-l+n-n'}{2}} h(t).
\end{equation}
\end{corollary}

\begin{lemma}
Let $l,m,n \in \mathbb{N}_0$ and $m \geq n$.
Then there exists a constant $C>0$, such that 
for $t>0$ and $j \in \mathbb{N}_0$, 
\begin{equation}
    \norm{\theta_t^j}_{C^l_{m,n}(\mathbb{H}^{2d})}
    \leq \begin{dcases*}
Ce^{\frac{d(d-1)}{4}t}t^{-\frac{l+m-n+d}{2}} & if $\: j =0$\\
\frac{ e^{\frac{d(d-1)}{4}t(1+j)} 
    C^j t^{j-\frac{l+d}{2}}}{\prod_{i=1}^j \max(1,i-\frac{l+d}{2})} & if $\: j \geq 1$
\end{dcases*}.
\end{equation}
\end{lemma}

\begin{corollary}
   Let $l,m,n \in \mathbb{N}_0$ and $m \geq n$.
   Then 
\begin{equation}
    \theta_t \in C^l_{m,n}\big(\mathbb{H}^{2d}\big), 
    \quad t >0,
\end{equation}
and there exists a continuous function $h$,
such that
   \begin{equation}
       \norm{\theta_t}_{C^l_{m,n}(\mathbb{H}^{2d})} 
       \leq t^{-\frac{l+m-n+d}{2}} h(t), \quad t>0.
   \end{equation}
\end{corollary}

\begin{lemma}
Let $n \in \mathbb{N}_0$. Then there exists a constant $c$, such that for $t>0$, $j \in \mathbb{N}$, and $x,y \in \mathbb{H}^d$,
\begin{equation}
\begin{split}
    &\int_{s \in t \Delta_j} \norm{(\theta_s 
    \shuffle_* (ic\Id B)^{\otimes j})(x,y)}_{B(S) \otimes \langle A_0 \rangle^{-n}S^1(H) \langle A_0 \rangle^n} \Id s \\
   & \leq \frac{(ct)^j}{j!} e^{ct} P_t(x,y).
    \end{split}
\end{equation}
 Moreover, let $l,l',m,m',n,n' \in \mathbb{N}_0$
    with $m' \geq m \geq n \geq n'$, and 
    $l' \geq l$. 
    For $t>0$ and $j \in \mathbb{N}$, define a linear 
    operator via 
\begin{equation}
  \big(R^{j,L}_tf\big)(x) := \int_{\mathbb{H}^d} r^j_t(x,y)f(y) \Id y, 
\end{equation}
where $f \in C^l_{m,n}(\mathbb{H}^d)$. 
Then $R^{j,L}_t$ is continuous as a map 
to $C^{l',S^1}_{m',n}$ and for $0 < t \leq t_0$,
there exists a $j,t$-independent constant $C$ 
such that 
\begin{equation}
\norm{R^{j,L}_t}_{B\left(C^l_{m,n},C^{l',S^1}_{m',n}\right) } \leq \frac{ e^{\frac{d(d-1)}{4}t(1+j)} 
    C^j t^{j-\frac{l'-l}{2}}}{\prod_{i=1}^j \max(1,i-\frac{l'-l}{2})}.
\end{equation}
\end{lemma}

\begin{corollary}\label{restimate}
    Let $n \in \mathbb{N}_0$. Then there
    is a constant $c$, such that for $t>0$, 
 \begin{equation}
        \sum_{j=1}^{\infty} \norm{r^j_t(x,y)}_{B(S) \otimes \langle A_0 \rangle^{-n}S^1(H) \langle A_0 \rangle^n} 
        \leq e^{ct} P_t(x,y), \quad x,y \in \mathbb{H}^d.
    \end{equation}    
    Moreover, let $l,m,n \in \mathbb{N}_0$, 
    with $m \geq n $, $1 \geq l$, and $j \in \mathbb{N}$. Then, for
    $0< t \leq t_0$, there exists a 
    $t,j$-independent constant $C$, such that
\begin{equation}
    \norm{r^j_t}_{C^{l,S^1}_{m,n}(\mathbb{H}^{2d})}
    \leq 
    \frac{ e^{\frac{d(d-1)}{4}t(1+j)} 
    C^j t^{j-\frac{l+d}{2}}}{\prod_{i=1}^j \max(1,i-\frac{l+d}{2})}.
\end{equation}

\end{corollary}

\subsection{From auxiliary to main integral kernels}

\begin{definition}
     For $t>0$ and $k \in \mathbb{N}$, let 
    \begin{equation}
        \kappa^{-,k}_ := \sum_{j=\frac{k-1}{2}}^{\infty} r_t^{2j+1}, \: \text{ $k$ odd}, 
        \quad \text{and} \quad \kappa^{+,k}_t := \sum_{j=\frac{k}{2}}^{\infty} r_t^{2j}, \: \text{ $k$ even}.
    \end{equation}
\end{definition}

\begin{definition}
    For $\epsilon, r > 0$, define
    \begin{equation}
        \begin{split}
            \rho_t^{\epsilon,r}(x,y) &:= 
           \mathbbm{1}_{B_r(o)}(x) 
           \mathbbm{1}_{B_r(o)}(y) 
 \:   \text{tr}_{S}        (ic \nabla_x + ic\nabla_y)\left((ic\nabla
            P_{\epsilon}\otimes P_{\epsilon}
            ) * \kappa^{-,d-2}_t \right. \\
&+ \left. (P_{\epsilon}\otimes P_{\epsilon}) * 
\left(A \kappa^{+,d-1}_t\right)  \right) (x,y), \\
   \rho_t^{r}(x,y) &:= 
           \mathbbm{1}_{B_r(o)}(x) \mathbbm{1}_{B_r(o)}(y) 
 \:   \text{tr}_{S}        (ic \nabla_x + ic\nabla_y)\left(ic\nabla_1
            \kappa^{-,d-2}_t \right. 
+ \left. 
A \kappa^{+,d-1}_t \right) (x,y), \\
\sigma_t(x,y) &:= - 
 \left[ A, \text{tr}_{S} \left( A \kappa^{-,d}_t + ic\nabla_1 \kappa_t^{+,d-1}\right)\right] (x,y).
        \end{split}
    \end{equation}
\end{definition}
The next Corollary will allow us to use the dominated
convergence Theorem in the computations 
performed in Section \ref{section6} and its proof works analogously to the proof of Corollary 5.2. in \cite{oliver23}.
\begin{corollary}\label{corollary29}
    Let $l,m,n \in \mathbb{N}_0$ with $m \geq n$,
    $t>0$ and $k \in \mathbb{N}$.
    Then 
    \begin{equation}
        \kappa_t^{\pm,k} \in C^l_{m,n}\big(\mathbb{H}^{2d}\big),
    \end{equation}
and there exist continuous functions $h_k$, such 
that
\begin{equation}
    \norm{\kappa_t^{-,k}}_{C_{m,n}^{l,S^1}(\mathbb{H}^{2d})}
    + \norm{\kappa_t^{+,k}}_{C_{m,n}^{l,S^1}(\mathbb{H}^{2d})} \leq h_k(t) t^{k-\frac{l+d}{2}}.
\end{equation}
Moreover, let $\epsilon, r >0$ and $\psi \in C_c^{\infty}(\mathbb{H}^d)$ with 
$\text{supp}(\psi) \in B_r(o)$. In that case, 
\begin{equation}
    \rho_t^{\epsilon, r}, \rho_t^{ r},
    \sigma_t \in C_{0,0}^{0,S^1}\big(\mathbb{H}^{2d}\big),
\end{equation}
and there exist continuous functions $f,g$ such that
\begin{equation}
    \begin{split}
        \norm{(\psi \otimes \psi) \rho_t^{\epsilon,r}}_{C_{0,0}^{0,S^1}\big(\mathbb{H}^{2d}\big)} &\leq f(t) t^{-\frac{1}{2}}, \\
 \norm{(\psi \otimes \psi) {\sigma_t}}_{C_{0,0}^{0,S^1}\big(\mathbb{H}^{2d}\big)} &\leq    g(t).  
    \end{split}
\end{equation}

\end{corollary}

\begin{lemma}[\cite{oliver23}, Lemma 5.3.]
    Let H be a separable Hilbert space and 
    let $\epsilon, r >0$. 
    Denote, for $f \in L^2(\mathbb{H}^d,H)$
    and $T \in B(L^2(\mathbb{H}^d,H))$
    with $x \in \mathbb{H}^d$,
\begin{align*}
(\chi_rf)(x)&:=\one_{B_r(o)}(x)f(x),           &   (\chi_r^LT)(f)&:= \chi_r Tf,              &  (\chi_r^RT)(f)&:= T \chi_r f,\\
\mathcal{P}_{\epsilon}f&:= P_{\epsilon}*f,         &  (\mathcal{P}^L_{\epsilon}T)(f)&:= 
\mathcal{P}_{\epsilon}Tf, &  
(\mathcal{P}^R_{\epsilon}T)(f)&:= T\mathcal{P}_{\epsilon}f.
\end{align*}
If, 
for  $p \in [1, \infty)$,
$T \in S^p(L^2(\mathbb{H}^d,H))$, 
then 
\begin{equation}
    \chi^L_r \chi^R_r \mathcal{P}^L_{\epsilon} 
    \mathcal{P}^R_{\epsilon} T \longrightarrow \chi^L_r \chi^R_r T, 
    \quad \epsilon \rightarrow 0^+,
\end{equation}
uniformly in $r$ in the 
$S^p(L^2(\mathbb{H}^d,H))$-norm, so that
the double limit
\begin{equation}
    \lim_{(\epsilon, r) \rightarrow (0, \infty)} 
     \chi^L_r \chi^R_r \mathcal{P}^L_{\epsilon} 
    \mathcal{P}^R_{\epsilon} T = T,
\end{equation}
exists and equals the iterated limits.

\end{lemma}

\begin{proposition}\label{prop1}
    Let $\epsilon, r >0$ and 
    $\mathscr{P}_{\epsilon}= \one_{S_x} \otimes \mathcal{P}_{\epsilon}$.
    Then for $t>0$, we have 
\begin{equation}
   \begin{split}
 & \textnormal{tr}_{S}\left(     e^{-tD^*D} - e^{-tDD^*}\right) \\
        &= 
 \lim_{r \rightarrow  \infty}         
        \int_0^t  
    \chi_r    \textnormal{tr}_{S} \left( \left[    
        ic \nabla, ic \nabla 
        K_s^{-,d-2} 
        +
        AK_s^{+,d-1} 
        \right] \right) \chi_r
        \\
        &- 
 \chi_r        \left[A, \textnormal{tr}_{S} \left(
        AK_s^{-,d} +ic \nabla K_s^{+,d-1}
         \right)\right]
  \chi_r        
         \Id s,
        \end{split}
\end{equation}    
with the integral converging in $S^1(L^2(\mathbb{H}^d,H))$-norm.
\end{proposition}

\begin{proposition}
    Let $t>0$, $j \in \mathbb{N}_0$ and 
    $k \in \mathbb{N}$.
    For $f \in L^2(S \otimes H)$ 
    and a.e. $x \in \mathbb{H}^d$,
\begin{equation}
    \begin{split}
 \int_{\mathbb{H}^d} \theta^j_t(x,y) f(y) 
 \Id y &= \left(  \int_{s \in \Delta_j}  
 e^{-sH_0} \shuffle (A^2 - M_0^2)^{\otimes j} \Id s f\right)(x),
 \\
  \int_{\mathbb{H}^d} \theta_t(x,y) f(y) 
  \Id y &= \left( e^{-tH_B}f \right)(x), \\
   \int_{\mathbb{H}^d} r^j_t(x,y) f(y) \Id y
   &=  \left( R^j_t f \right)(x),\\
 \int_{\mathbb{H}^d} \kappa^{\pm,k}_t(x,y) 
 f(y) \Id y &= \left( K^{\pm,k}_t f \right)(x). 
    \end{split}
\end{equation}
\end{proposition}

\begin{proposition}\label{kernels2}
    Let $t, \epsilon, r >0$ and 
    $f \in L^2(\mathbb{H}^d, H)$.
    Then, for a.e. $x \in \mathbb{H}^d$,
\begin{equation*}
    \begin{split}
        \int_{\mathbb{H}^d} \rho_t^{\epsilon,r}(x,y) f(y) \Id y
        & = \left(
\chi_r   
\overline{\textnormal{tr}_{S} \left( \left[    
        ic \nabla, ic \nabla 
        \mathscr{P}_{\epsilon}
        K_s^{-,d-2} \mathscr{P}_{\epsilon}
        + \mathscr{P}_{\epsilon}
        AK_s^{+,d-1} \mathscr{P}_{\epsilon}
        \right] \right)} \chi_r f \right)(x)      
        \\
  \int_{\mathbb{H}^d} \chi_r(x) \sigma_t(x,y) \chi_r(y) f(y) \Id y &= 
   - 
\left( \chi_r       \left[A, \textnormal{tr}_{S} \left(
        AK_s^{-,d} +ic \nabla K_s^{+,d-1}
         \right)\right]
  \chi_r f\right)(x).
    \end{split}
\end{equation*}   
\end{proposition}
\begin{proof}
For the first equality, let 
$W_t := (ic\nabla P_{\epsilon}\otimes P_{\epsilon}
            ) * \kappa^{-,d-2}_t 
+ (P_{\epsilon}\otimes P_{\epsilon}) * 
\left(A \kappa^{+,d-1}_t\right)$. 
Since the spinor heat kernel $P_t$ is a Schwartz 
function, $W_t$ is smooth. 
Integration by parts, for 
$f \in C^{\infty}_c(\mathbb{H}^d,S)$, yields
\begin{equation*} 
\begin{split}
    &\int_{\mathbb{H}^d} \left( \left(ic\nabla_1 + ic\nabla_2 \right) W_t\right)(x,y) f(y) \Id y \\
   &  = ic \nabla \left(ic \nabla 
        \mathscr{P}_{\epsilon}
        K_s^{-,d-2} \mathscr{P}_{\epsilon}
        + \mathscr{P}_{\epsilon}
        AK_s^{+,d-1} \mathscr{P}_{\epsilon}\right)f\\
    & \quad \quad - 
    i \sum_{j=1}^d c^j 
    \left(ic \nabla 
        \mathscr{P}_{\epsilon}
        K_s^{-,d-2} \mathscr{P}_{\epsilon}
        + \mathscr{P}_{\epsilon}
        AK_s^{+,d-1} \mathscr{P}_{\epsilon}\right)
        \nabla_{e_j}f,
 \end{split}       
\end{equation*}
which, after applying 
the trace and commuting the $c^j$'s,
results in
\begin{equation*}
    \begin{split}
&\int_{\mathbb{H}^d} \text{tr}_{S}\left( \left(ic\nabla_1 + ic\nabla_2 \right) W_t\right)(x,y) f(y) \Id y \\
&= \text{tr}_{S} \left( \left[    
        ic \nabla, ic \nabla 
        \mathscr{P}_{\epsilon}
        K_s^{-,d-2} \mathscr{P}_{\epsilon}
        + \mathscr{P}_{\epsilon}
        AK_s^{+,d-1} \mathscr{P}_{\epsilon}
        \right]f \right)(x).
    \end{split}
\end{equation*}
Due to $\mathscr{P}_{\epsilon}$, the right-hand side
admits a bounded extension in $L^2(\mathbb{H}^d,H)$ and since the latter must equal its closure, we can conclude. \\
The claim for $\sigma_t$ follow by Fubini's theorem, 
$\kappa_t^{\pm,k} \in C^{1,S^1}_{m,n}(\mathbb{H}^{2d})$
and the fact that $\chi_r$ is compactly supported.

\end{proof}

\section{The trace formula}\label{section6}
In this section, we start by showing that the 
interior trace  reduces to a boundary integral at infinity (Proposition \ref{prop6.5}), before 
computing the iterated trace in Theorem 
\ref{maintheorem}.
The argument proceeds by 
\begin{enumerate}
    \item first discarding 
    the commutator term with $\sigma_t$ from 
    Proposition 
    \ref{kernels2} (Lemma \ref{lemmarholimit}), and then the factor appearing in the expression of 
    $\rho_t^{\epsilon,r}$
    involving 
    $ic \nabla \kappa_t^{-,d-2}$
    (Lemma \ref{lemma6.2}) by employing a symmetrization procedure as 
    described in Lemma \ref{opproddec}, as well as 
   estimates involving the spinor heat kernel (Subsection \ref{subsectionspinorestimates});
    \item freezing coefficients in the $Q_t$ factors (Lemma \ref{lemma3}) and reducing the 
    $\kappa^{+,d-1}$ term to an expression 
    that depends only on the diagonal of the spinor heat kernel, as well as the frozen operators $e^{-uA^2(x_0)}$ and $ic( \Id B)(x_0)$, for some $x_0 \in \mathbb{H}^d$;
    \item showing that only the term of order $d-1$
   survives in the rescaled limit (Lemma \ref{highervanish}).
\end{enumerate}

\begin{lemma}\label{lemmarholimit}
    For $t>0$, we have 
    \begin{equation}
 \textnormal{tr}_{L^2(\mathbb{H}^d,H)}
    \textnormal{tr}_{S} \left(
    e^{-tD^*D} - e^{-tDD^*}\right)       = \textnormal{tr}_{S \otimes H} \left(
     \lim_{R \rightarrow \infty} 
     \int_{\mathbb{H}^d} \int_0^t 
     \rho_s^R(x,x) \Id s \Id x
        \right),
    \end{equation}
    where the right-hand side converges in  $S^1(S \otimes H)$ and the first integral 
    is a Bochner integral in $L^1(\mathbb{H}^d, S^1(S \otimes H))$.
\end{lemma}
\begin{proof}
    The proof is similar to the proof of 
    Lemma 5.5. in \cite{oliver23}. 
    We know that 
    $$-\int_0^t \left[ A,  A K^{-,d}_s + ic\nabla K_s^{+,d-1}\right] \Id s,$$
    admits 
    $$\Sigma_t(x,y) := \int_0^t \tilde \sigma_s(x,y) \Id s, 
    \quad 
  \tilde \sigma_s(x,y)   := - \left[ A,  A \kappa^{-,d}_s + ic\nabla_1 \kappa_s^{+,d-1}\right],$$
    as an operator-valued integral kernel, given that $s \mapsto \tilde \sigma_s$ is integrable at $0$ in 
    $C^{0,S^1}_{m,n}(\mathbb{H}^{2d})$-norm since, by Corollary \ref{corollary29}, 
    $$ \tilde \sigma_s \in C_{m,n}^{0,S^1}(\mathbb{H}^{2d}), \quad \text{and} \quad 
  \norm{\tilde \sigma_s}_{C_{m,n}^{0,S^1}(\mathbb{H}^{2d})} \leq h(s) s^{d-1 - \frac{d+1}{2}} = h(s) s^{\frac{d-3}{2}}, $$
  for all $m \geq n$, $m,n \in \mathbb{N}_0$,
  and
  for a function $h$ which is continuous at $0$. Moreover, $\tilde \sigma_s$ is the operator-valued integral
  kernel of $-\left[ A,  A K^{-,d}_s + ic\nabla K_s^{+,d-1}\right]$ by 
  Proposition \ref{kernels2}. Since, for $R>0$, 
  $\chi_R \textnormal{tr}_{S} \left(
    e^{-tD^*D} - e^{-tDD^*}\right)  \chi_R \in S^1(L^2(\mathbb{H}^d, H))$, we know by Theorem \ref{mercer} that it admits an operator valued integral kernel $k_{t,R}$ such that
    \begin{equation} 
        \begin{split}
            \mathbb{H}^d &\longrightarrow L^1(\mathbb{H}^d, S^1(H)) \\
    h &\longmapsto   (x \mapsto k_{t,R}(x,L_hx)),     
        \end{split}
    \end{equation}
    is a bounded continuous map, and given that 
    $\chi_R \otimes \chi_R \cdot \textnormal{tr}_{S} \Sigma_t \in  S^1(L^2(\mathbb{H}^d, H))$ as well, Proposition \ref{prop1} 
    implies that $\int_0^t \rho_s^R \Id s \in S^1(L^2(\mathbb{H}^d, H))$.
    Hence, by Theorem \ref{mercer},
    \begin{equation}\label{eq:threesummands}
        \begin{split}
 \textnormal{tr}_{L^2(\mathbb{H}^d,H)}  \left(  \chi_R \textnormal{tr}_{S} \left(
    e^{-tD^*D} - e^{-tDD^*}\right)  \chi_R \right) &= \int_{\mathbb{H}^d} \chi_R(x) \textnormal{tr}_{S \otimes H} \Sigma_t(x,x) \Id x     \\
    &+ \int_{\mathbb{H}^d} \textnormal{tr}_{H}  \left( \int_0^t \rho_s^R(x,x) \Id s \right) \Id x. 
        \end{split}
    \end{equation}
Now, by Proposition \ref{prop1}, $\tau_s(x,y) := A \kappa_s^{-,d} + ic \nabla_1 \kappa_s^{+,d-1} \in C_{m,n}^{0,S^1}(\mathbb{H}^{2d})$ for all $m,n \in \mathbb{N}_0$, and so 
\begin{equation*}
    \textnormal{tr}_{S \otimes H} A(x) \tau_s(x,x) = \textnormal{tr}_{S \otimes H} \langle A_0 \rangle \tau_s(x,x) A(x) \langle A_0 \rangle^{-1} = \textnormal{tr}_{S \otimes H} \overline{\tau_s(x,x) A(x)},
\end{equation*}
where we used that $A(x) \langle A_0\rangle^{-1}$ is a bounded operator. Therefore, 
\begin{equation*}
    \textnormal{tr}_{S \otimes H} \sigma_s(x,x) = \textnormal{tr}_{S \otimes H} [A(x), \tau(x,x)] =0,
\end{equation*}
and thus 
\begin{equation*}
    \textnormal{tr}_{S \otimes H} \Sigma_t(x,x) = \int_0^t \textnormal{tr}_{S \otimes H} \sigma_s(x,x) \Id s = 0,
\end{equation*}
which allows us to conclude with Equation \eqref{eq:threesummands} and Proposition \ref{prop1}.
    
\end{proof}

\begin{lemma}[\cite{oliver23}, Lemma 5.7.]\label{opproddec}
    With the notations from Definition 
    \ref{deltadef} and an elementary 
    tensor of operator functions 
    $\textbf{T} := T_1 \otimes ... \otimes T_n$,
    \begin{equation}
    \begin{split}    &\textnormal{OpProd}\left(\tau_{\big(y/2\big)^{-1}}
        - \tau_{\big(y/2\big)}\right)^2  
        \textbf{T} \\
    &= \sum_{j =1}^n \bigg( \prod_{i = 1}^{j-1} \tau_{y_i} T_i \bigg) \mathcal{W}^2_{y_j}
    (T_j) \bigg( \prod_{i =j +1}^n \tau_{y_i^{-1}}
    T_i  \bigg) \\
    & \quad + 2 \sum_{1 \leq k < l \leq n} 
    \bigg( \prod_{i = 1}^{k-1} T_i \bigg) 
    \tau_{\frac{y_k}{2}} \delta_{\frac{y_k}{2}} T_k \bigg( \prod_{j = k+1}^{l-1} 
    \tau_{y_j} T_j \bigg) 
    \tau_{\big(\frac{y_l}{2}\big)^{-1}} \delta_{\frac{y_l}{2}} T_l
    \bigg( \prod_{m = l+1}^n \tau_{y_m^{-1}} T_m 
    \bigg).
        \end{split}
    \end{equation}
\end{lemma}

\subsection{Further estimates concerning the spinor heat kernel}\label{subsectionspinorestimates}

\begin{lemma}
Let \(d\geq3\) be odd and $m=\frac{d-1}{2}$.
For \(x,y\in\mathbb H^d\), set 
$r = d(x,y)$ and denote by $p_t(r)$
the scalar part of the spinor heat kernel
as in Equation \eqref{eq:scalerheatkernelpart}.
Then, there exists $c >4$  and a constant $C=C(c,d)>0$ 
such that, for all \(t>0\) and \(r\geq0\),
\begin{equation}\label{eq:scalarheatkernelbound}
    |p_t(r)|
\leq
C
\left(
t^{-d/2}+t^{-1/2}
\right)
\exp\left(-\frac{r^2}{ct}\right)
\exp\left(-\frac{d-1}{2}r\right).
\end{equation}
Moreover, for every \(\eta\geq0\), there exists $C=C(c,d,\eta)>0$ 
such that
\begin{equation}\label{eq:scalarheatkernelbound2}
    \left(
\frac{r}{t}
+
\frac{1}{\sqrt t}
+
1
\right)^\eta
|p_t(r)|
\leq
C
\left(
t^{-d/2}+t^{-1/2}
\right)
\left(
1+t^{-\eta/2}
\right)
\exp\left(-\frac{r^2}{ct}\right)
\exp\left(-\frac{d-1}{2}r\right).
\end{equation}
\end{lemma}
\begin{proof}
Let $f_t(r):=
\left(\cosh\frac r2\right)^{-1}e^{-r^2/(4t)}$
and introduce the  differential operator
$L:=-\frac{1}{\sinh r}\frac{d}{dr}$.
We start by proving that, for every $c_0>4$, there exists a constant
$C=C(c_0,d)$ such that
\begin{equation}
\left|
\cosh\left(\frac r2\right)L^m f_t(r)
\right|
\leq
C(1+t^{-m})
\exp\left(-\frac{r^2}{c_0t}\right)e^{-mr},
\label{eq:auxiliary-spinor-estimate}
\end{equation}
for all $t>0$ and $r\geq0$. 
For that matter, we
distinguish the regions $0\leq r\leq1$ and $r\geq1$.
For $0\leq r\leq1$, set $s=r^2$. Since
$\frac{d}{dr}=2r\frac{d}{ds}$,
we may write $L
=
-\frac{2r}{\sinh r}\frac{d}{ds}
=
-a(s)\frac{d}{ds}$, where 
$a(s):=\frac{2\sqrt{s}}{\sinh\sqrt{s}}$.
Notice that the function $a$ is smooth on $[0,1]$, with $a(0)=2$.
Consequently, for every fixed $m$,
\[
L^m=\sum_{j=1}^m a_j(s)\frac{d^j}{ds^j},
\]
where the coefficients $a_j$ are smooth  on $[0,1]$.
Moreover, 
$f_t(r)=b(s)e^{-s/(4t)}$, 
with $b(s):=\left(\cosh\frac{\sqrt{s}}2\right)^{-1}$
and where 
 $b$ is smooth on $[0,1]$ with $b(0)=1$. Hence, by the Leibniz rule and for $1\leq j\leq m$,
\[
\left|
\frac{d^j}{ds^j}
\left(
b(s)e^{-s/(4t)}
\right)
\right|
\leq
C_j(1+t^{-j})e^{-s/(4t)}
\leq
C(1+t^{-m})e^{-s/(4t)},
\]
so that $|L^m f_t(r)|
\leq
C(1+t^{-m})e^{-r^2/(4t)}$.
Since $0\leq r\leq1$, the functions
$\cosh(r/2)$ and $e^{-mr}$ are bounded above and below by
positive constants. Thus, after multiplying by one and absorbing constants into $C$,
\begin{equation}
\left|
\cosh\left(\frac r2\right)L^m f_t(r)
\right|
\leq
C(1+t^{-m})
\exp\left(-\frac{r^2}{4t}\right)e^{-mr},
\qquad 0\leq r\leq1.
\label{eq:small-r}
\end{equation}
We next consider $r\geq1$. In this region,
$\left|
\frac{d^k}{dr^k}\frac{1}{\sinh r}
\right|
\leq
C_k e^{-r}$, $k\geq0$, so that
repeated applications of $L=-(\sinh r)^{-1}d/dr$ and Leibniz' rule
 give
\[
L^m f
=
\sum_{j=1}^m A_{m,j}(r)f^{(j)}(r),
\]
where $|A_{m,j}(r)|\leq C e^{-mr}$ for 
$r\geq1$.
On the other hand, $f_t(r)
=
\sech\left(\frac r2\right)e^{-r^2/(4t)}$
and, 
for every $k\geq0$,
\[
\left|
\frac{d^k}{dr^k}\sech\left(\frac r2\right)
\right|
\leq
C_k e^{-r/2},
\qquad r\geq1.
\]
Furthermore, for every $\varepsilon>0$ and every $1\leq j\leq m$,
the standard Gaussian derivative estimate gives
\[
\left|
\frac{d^j}{dr^j}e^{-r^2/(4t)}
\right|
\leq
C_{j,\varepsilon}
(1+t^{-j/2})
\exp\left(
-\frac{r^2}{4(1+\epsilon)t}
\right).
\]
By Leibniz' rule, it follows that
\[
|f_t^{(j)}(r)|
\leq
C_\varepsilon
(1+t^{-m/2})
\exp\left(
-\frac{r^2}{4(1+\epsilon)t}
\right)e^{-r/2}.
\]
Since $1+t^{-m/2}\leq C(1+t^{-m})$,
we obtain, together with $\cosh\left(\frac r2\right)\leq Ce^{r/2}$,
\[
\left|
\cosh\left(\frac r2\right)L^m f_t(r)
\right|
\leq
C(1+t^{-m})
\exp\left(
-\frac{r^2}{4(1+\epsilon)t}
\right)e^{-mr}, \qquad r \geq 1.
\]
Choosing $\epsilon>0$ sufficiently small so that
$4(1+\epsilon)<c_0$, 
proves \eqref{eq:auxiliary-spinor-estimate} for $r\geq1$ and
together with \eqref{eq:small-r}, this proves
\eqref{eq:auxiliary-spinor-estimate} for all $r\geq0$.
Take now $c_0<c$. From
\eqref{eq:auxiliary-spinor-estimate}
and $t^{-1/2}(1+t^{-m})
=
t^{-1/2}+t^{-d/2}$, we obtain 
Equation \eqref{eq:scalarheatkernelbound}.
For the weighted estimate, fix $\eta\geq0$. For every
$\delta>0$, the elementary inequalities
\[
\left(
\frac rt+\frac1{\sqrt t}+1
\right)^\eta
\leq
C_\alpha
(1+t^{-\eta/2})
\left(1+\frac r{\sqrt t}\right)^\eta,
\]
and
\[
\left(1+\frac r{\sqrt t}\right)^\alpha
\leq
C_{\alpha,\delta}
\exp\left(\frac{\delta r^2}{t}\right),
\]
hold for all $r,t>0$. Choose $c_0<c$ and $\delta>0$ sufficiently
small so that $\frac{1}{c_0}-\delta>\frac{1}{c}$.
Multiplying \eqref{eq:auxiliary-spinor-estimate} by the polynomial
factor above, the resulting inequality 
 by $C_dt^{-1/2}$ and again using 
 $t^{-1/2}(1+t^{-m})
=
t^{-1/2}+t^{-d/2}$, 
we conclude with Equation \eqref{eq:scalarheatkernelbound2}.

\end{proof}

\begin{lemma}\label{scalerheatkerneldecreasing}
    Let $d=2m+1$ be odd and let $p_t$ denote the scalar part of the
spinor heat kernel on $\mathbb{H}^d$. 
Then
\[
    p_t'(s)<0
    \qquad\text{for all }s>0.
\]
In particular, $p_t$ is strictly decreasing as a function of $s$.
\end{lemma}
\begin{proof}
    Let $a:=\frac{1}{4t}$ and $
    x:=\cosh s$.
Then $x>1$ for $s>0$, and
$ \cosh\frac{s}{2}
    =
    \sqrt{\frac{x+1}{2}}$ by the half angle formula. 
Define now
\[
    F_a(x)
    :=
    \sqrt{\frac{2}{x+1}}\,
    e^{-a(\operatorname{arcosh}x)^2}
\]
and
\[
    G_m(x):=(-1)^mF_a^{(m)}(x).
\]
We start by showing that $F_a$ is completely monotone on $(1,\infty)$, so, in particular, $(-1)^k F_a^{(k)}(x)>0$ for 
$x >1$ and $k \geq 0$.
To this end, let $\phi(x):=(\operatorname{arcosh}x)^2$ for $x>1$.
We claim that $\phi$ is a Bernstein function on $(1,\infty)$ \footnote{Recall that $f$ is a Bernstein function if and only if it is nonnegative and $f'$ is a completely monotone function. }. 
One computes 
\[
    \phi'(x)
    =
    \frac{2\operatorname{arcosh}x}{\sqrt{x^2-1}}.
\]
For $x>1$, we have the identity
\begin{equation}
    \frac{\operatorname{arcosh}x}{\sqrt{x^2-1}}
    =
    \int_0^\infty
    \frac{du}{x+\cosh u},
    \label{eq:arcosh-integral}
\end{equation}
obtained by performing the 
change of variables $v=\tanh(u/2)$, 
thus $\cosh u=\frac{1+v^2}{1-v^2}$ and
$du=\frac{2\,dv}{1-v^2}$, so that
\[
\begin{aligned}
    \int_0^\infty\frac{du}{x+\cosh u}
    &=
    2\int_0^1
    \frac{dv}{(x+1)-(x-1)v^2} 
    =
    \frac{2}{\sqrt{x^2-1}}
    \operatorname{arctanh}
    \sqrt{\frac{x-1}{x+1}} \\
    &=
    \frac{\operatorname{arcosh}x}{\sqrt{x^2-1}}.
\end{aligned}
\]
Since
\[
    \frac{1}{x+\cosh u}
    =
    \int_0^\infty
    e^{-\lambda(x+\cosh u)}\,d\lambda,
\]
Fubini's Theorem gives
\begin{equation}
    \phi'(x)
    =
    2\int_0^\infty
    e^{-\lambda x}
    \left(
        \int_0^\infty
        e^{-\lambda\cosh u}\,du
    \right)
    d\lambda.
    \label{eq:phi-complete-monotone}
\end{equation}
Hence $\phi'$ is completely monotone, 
since it is the Laplace transform of a 
positive measure,
and therefore $\phi$ is a Bernstein function.
By the Bernstein representation Theorem, for every $a>0$, $e^{-a\phi(x)}$
is completely monotone. 
Moreover,
$x\longmapsto (x+1)^{-1/2}$
is completely monotone, since
\[
    (-1)^k\frac{d^k}{dx^k}(x+1)^{-1/2}
    =
    \frac{(2k)!}{4^k k!}
    (x+1)^{-k-1/2}>0.
\]
Finally, the product of two completely monotone functions is completely monotone
and hence
\[
    F_a(x)
    =
    \sqrt{2}\,(x+1)^{-1/2}
    e^{-a\phi(x)},
\]
is completely monotone, so that, 
by what we have just shown, $G_m(x)>0$ for
$x > 1$.
Now, Equation \eqref{eq:scalerheatkernelpart} can be written as
\begin{equation}
    p_t(s)
    =
    C_{d,t}
    \sqrt{\frac{x+1}{2}}\,
    G_m(x),
    \qquad x=\cosh s,
    \label{eq:pt-Gm}
\end{equation}
where $C_{d,t}>0$ is independent of $s$.
Since $F_a$ is completely monotone, it admits a Laplace-transform
representation
\[
    F_a(x)
    =
    \int_0^{\infty}
    e^{-(x-1)\lambda}\,d\mu_a(\lambda)
\]
for some positive Borel measure $\mu_a$. Consequently,
\[
    G_m(x)
    =
    \int_0^{\infty}
    \lambda^m e^{-(x-1)\lambda}\,d\mu_a(\lambda),
\]
and
\[
    -G_m'(x)
    =
    \int_0^{\infty}
    \lambda^{m+1}e^{-(x-1)\lambda}\,d\mu_a(\lambda).
\]
Writing $ M_k(x)
    :=
    \int_0^{\infty}
    \lambda^k e^{-(x-1)\lambda}\,d\mu_a(\lambda)$,
we have $G_m=M_m$ and $-G_m'=M_{m+1}$.
By the Cauchy--Schwarz inequality,
\[
    M_{k+1}^2(x)\leq M_k(x)M_{k+2}(x),
\]
so that the ratios $M_{k+1}/M_k$ are nondecreasing in $k$. Hence
\begin{equation}
    -\frac{G_m'(x)}{G_m(x)}
    =
    \frac{M_{m+1}(x)}{M_m(x)}
    \geq
    \frac{M_1(x)}{M_0(x)}
    =
    -\frac{F_a'(x)}{F_a(x)},
    \label{eq:ratio-comparison}
\end{equation}
with
a direct differentiation giving
\[
    -\frac{F_a'(x)}{F_a(x)}
    =
    \frac{1}{2(x+1)}
    +
    \frac{2a\operatorname{arcosh}x}{\sqrt{x^2-1}}
    >
    \frac{1}{2(x+1)}.
\]
Therefore,
\[
    G_m'(x)
    +
    \frac{G_m(x)}{2(x+1)}
    <0.
    \label{eq:Gm-negative}
\]
Differentiating \eqref{eq:pt-Gm} with respect to $x$ now yields
\[
\begin{aligned}
    \frac{d}{dx}p_t(\operatorname{arcosh}x)
    &=
    C_{d,t}
    \sqrt{\frac{x+1}{2}}
    \left(
        G_m'(x)
        +
        \frac{G_m(x)}{2(x+1)}
    \right)\\
    &<0.
\end{aligned}
\]
Finally, $\frac{dx}{ds}=\sinh s>0,$ $
    s>0$, and hence $ p_t'(s)<0$, $s>0$.
Since $p_t$ is even, $p_t'(0)=0$, so $p_t$ is strictly decreasing on
$[0,\infty)$.

\end{proof}

\begin{lemma}\label{lemmageometricgaussian}
Let $m=(d-1)/2$ and $h := \frac{st}{s+t}$. For every $\eta>0$ there exists
$C=C(\eta,d)$ such that, for all $s,t>0$ and all
$x,y\in\mathbb H^d$,
\begin{equation}
\begin{aligned}
&
\int_{\mathbb H^d}
\exp\left(
-\eta\frac{d(x,z)^2}{s}
-\eta\frac{d(z,y)^2}{t}
\right)
\\
&\qquad\qquad\qquad\qquad
\cdot
\exp\left(
-m\bigl(d(x,z)+d(z,y)-d(x,y)\bigr)
\right)
\,dz
\\
&\qquad\leq
C_{\eta,d} \exp{\left(-\frac{c\eta d(x,y)^2}{s+t}\right)}
\sqrt{h}
\min \left\{\sqrt{h}, h^{\frac{d-1}{2}}\right\} .
\end{aligned}
\label{eq:hyperbolic-gaussian-lemma}
\end{equation}
Moreover, there exist
$C=C(a,\eta,d)$ and $c =c(d)$
such that, for every $R>0$,
\begin{equation}
\begin{aligned}
&
\int_{\{d(x,z)\geq R\}}
\exp\left(
-\eta\frac{d(x,z)^2}{s}
-\eta\frac{d(z,y)^2}{t}
\right)
\\
&\qquad\qquad\qquad\qquad
\cdot
\exp\left(
-m\bigl(d(x,z)+d(z,y)-d(x,y)\bigr)
\right)
\,dz
\\
&\qquad\leq
C_{\eta,d} \exp{\left(-\frac{c\eta d(x,y)^2}{s+t}\right)}
\sqrt{h}
\min \left\{\sqrt{h}, h^{\frac{d-1}{2}}\right\}
\exp\left(-\frac{c\eta R^2}{s+t}\right).
\end{aligned}
\label{eq:hyperbolic-gaussian-tail}
\end{equation}
\end{lemma}

\begin{proof}
Let us set $D:=d(x,y)$.
If $D=0$, then $x=y$, and
the proof is immediate from polar coordinates centered around
$x$. We may therefore assume that $D>0$. \\
Let $\gamma$ be the unique geodesic joining $x$ to $y$, so that 
$x=\gamma(0)$ and 
$y=\gamma(D)$.
We use Fermi
coordinates  around $\gamma$.
More precisely, for $u \in \mathbb{R}$, let
$\rho \geq 0$ denote the distance from 
$z$ to $\gamma$, and $\omega\in\mathbb S^{d-2}$ the direction of the normal geodesic 
from the nearest point $\gamma(u)$ to
$z$.
Thus, away from $\gamma$, we may write 
$ \mathbb{H}^d \ni z=(u,\rho,\omega)$ where 
$\rho >0$.
The hyperbolic metric in these coordinates is
\[
d\rho^2+\cosh^2\rho\,du^2+\sinh^2\rho\,d\omega^2,
\]
and consequently
\begin{equation}
dz=
(\sinh\rho)^{d-2}\cosh\rho
\,d\rho\,du\,d\omega.
\label{eq:fermi-volume}
\end{equation}
The geodesic $\gamma$ itself corresponds 
to $\rho =0$ and has measure zero, so it
does not effect any of the integrals below.
Set now $r := d(x,z)$ and $q:=d(z,y)$.
We first estimate the distance
$E:=r+q-D$
For $a\geq0$ and $\rho\geq0$, define
\[
F(a,\rho)
:=
\operatorname{arcosh}(\cosh a\cosh\rho)-a.
\]
Differentiating with respect to $\rho$ gives
\[
\partial_\rho F(a,\rho)
=
\frac{\sinh\rho\,\cosh a}
{\sqrt{\cosh^2a\cosh^2\rho-1}},
\]
which, together with $\sqrt{\cosh^2a\cosh^2\rho-1}
\leq
\cosh a\cosh\rho$, leaves us with
$\partial_\rho F(a,\rho)
\geq
\tanh\rho$.
As $F(a,0)=0$, it follows that
\begin{equation}
F(a,\rho)
\geq
\int_0^\rho\tanh v\,dv
=
\log\cosh\rho,
\qquad
a,\rho\geq0.
\label{eq:F-lower-bound}
\end{equation}
 By
the hyperbolic law of cosines, 
$\cosh r=\cosh\rho\,\cosh \vert u \vert$
and 
$\cosh q=\cosh\rho\,\cosh \vert D-u \vert$, so
\[
r= \vert u \vert +F( \vert u \vert,\rho),
\qquad
q= \vert  D-u \vert +F( \vert D-u \vert,\rho).
\]
Hence
\begin{equation}
E =r+q-D 
=
F(\vert u \vert ,\rho)+F(\vert D-u \vert ,\rho)
\geq
2\log\cosh\rho.
\label{eq:excess-distance-final}
\end{equation} 
Since $2m=d-1$, we therefore have
$e^{-m(r+q-D)}
\leq
(\cosh\rho)^{-(d-1)}$,
which,
combined with \eqref{eq:fermi-volume}, gives
\[
e^{-m(r+q-D)}
(\sinh\rho)^{d-2}\cosh\rho
\leq
\frac{\sinh^{d-2}\rho}{\cosh^{d-2}\rho}
=
\tanh^{d-2}\rho.
\]

Moreover, we have the estimates 
$r^2 \geq u^2 + \rho^2$ and 
$q^2 \geq (D-u)^2 + \rho^2$, where the first follows from 
$\cosh r = \cosh  u \cosh \rho$ and the inequality
$\cosh \sqrt{u^2 + \rho^2} \leq \cosh u \cosh \rho $. A similar argument 
is true for $q$.
Thus 
\begin{equation*}
\frac{r^2}{s} + \frac{q^2}{t} \geq \frac{u^2 + \rho^2}{s} + \frac{(D-u)^2 + \rho^2}{t} 
= \frac{D^2}{s+t}
+
\frac{s+t}{st}
\left(
u-\frac{sD}{s+t}
\right)^2 + \frac{\rho^2}{h^2}.
\end{equation*}
Using the  estimates above and integrating first in $\omega$, we obtain
\begin{equation}
\begin{aligned}
&
\int_{\mathbb H^d}
e^{-\eta r^2/s-\eta q^2/t}
e^{-m(r+q-D)}
\,dz
\\
&\qquad\leq
C
e^{-\eta D^2/(s+t)}
\int_{\mathbb R}
e^{-\eta\frac{s+t}{st}
\left(u-\frac{sD}{s+t}\right)^2}
\,du
\\
&\qquad\qquad\qquad\qquad
\cdot
\int_0^\infty
e^{-\eta\rho^2(\frac1s+\frac1t)}
\tanh^{d-2} \rho
\,d\rho.
\end{aligned}
\label{eq:separated-integrals}
\end{equation}
The $u$-integral is a one-dimensional Gaussian integral, and hence
\[
\int_{\mathbb R}
e^{-\eta\frac{s+t}{st}
\left(u-\frac{sD}{s+t}\right)^2}
\,du
\leq
C_\eta
\sqrt{h}.
\]
For the $\rho$-integration, we have, when
$0 \leq h \leq 1$,
\begin{equation*}
    \int_0^\infty
e^{-\eta\rho^2/h}
\tanh^{d-2} \rho
\,d\rho \leq  \int_0^\infty
e^{-\eta\rho^2/h}
\rho^{d-2}
\,d\rho \leq C_{\eta,d} h^{\frac{d-1}{2}},
\end{equation*}
and for $h \geq 1$, $\tanh^{d-2} \rho \leq 1$, so
\begin{equation*}
    \int_0^\infty
e^{-\eta\rho^2/h}
\tanh^{d-2} \rho
\,d\rho \leq C \sqrt{h},
\end{equation*}
 thus 
proving Equation \eqref{eq:hyperbolic-gaussian-lemma}. 
Remains to show Equation   \eqref{eq:hyperbolic-gaussian-tail}.
Fix $\theta\in(0,1)$. 
Then 
$e^{-\eta r^2/s} = e^{-\theta\eta r^2/s}e^{-(1- \theta)\eta r^2/s}$ and 
$e^{-(1- \theta)\eta r^2/s} \leq e^{-(1- \theta)\eta R^2/(s+t)}$ for $r = d(x,z)\geq R$, so that the left-hand side of Equation  \eqref{eq:hyperbolic-gaussian-tail} can be bounded by 
\begin{equation*}
    \begin{split}
        e^{-(1- \theta)\eta R^2/(s+t)} 
        \int_{\mathbb{H}^d} 
    e^{-\theta\eta r^2/s-\theta\eta q^2/t}
e^{-m(r+q-D)}
\,dz,    
    \end{split}
\end{equation*}
which, together with \eqref{eq:hyperbolic-gaussian-lemma}, lets us conclude.
\end{proof}

\begin{rem}\label{remupperboundmax}
    As this is going to be useful later, 
    note that, for $s,t >0$ and $h = \frac{st}{s+t}$, we have 
    \begin{equation}
        A_s A_t \sqrt{h}
\min \left\{\sqrt{h}, h^{\frac{d-1}{2}}\right\} \leq C_d(1+ A_{s+t}),
    \end{equation}
    with notation as in Lemma \ref{lemmageometricgaussian}.
\end{rem}

\begin{lemma}\label{lemmaspinorchainbound1}
Let $u_j>0$, $j\in\{0,\dots,k\}$, and set
\[
t_0:=\sum_{j=0}^k u_j,
\qquad
V_i:=u_0+\cdots+u_{i-1},
\qquad
T_i:=u_i+\cdots+u_k,
\]
so that $V_i+T_i=t_0$. Let $\delta\geq0$ and $r_0>0$. Then, for every
$i\in\{1,\dots,k\}$,
\begin{equation}
\begin{aligned}
&
\bigg|
\int_{(\mathbb H^d)^k}
\mathbf 1_{\{d(z_i,x)\geq r_0\}}
\left(
\frac{d(x,z_1)}{u_0}
+\frac{1}{\sqrt{u_0}}+1
\right)^\delta
\\
&\qquad\qquad\qquad\qquad
\cdot
P_{u_0}(x,z_1)
P_{u_1}(z_1,z_2)
\cdots
P_{u_k}(z_k,x)
\,dz
\bigg|
\\
&\qquad\leq
C_{c,d,\delta}
\left(1+u_0^{-\delta/2}\right)
\left(1+A_{V_i}\right)
\left(1+T_i^{-(d-1)/2}\right)
\exp\left(-c\frac{r_0^2}{V_i}\right),
\end{aligned}
\label{eq:epsilon-zero-bound}
\end{equation}
where
\[
A_t:=t^{-d/2}+t^{-1/2}.
\]
Note that for $i=1$, the right-hand side can be replaced by
\begin{equation}\label{eq:i=1casebound}
    C_{c,d,\delta}
\left(1+u_0^{-\delta/2}\right)
\left(1+A_{t_0}\right)
\exp\left(-c\frac{r_0^2}{u_0}\right).
\end{equation}
\end{lemma}
\begin{proof}
We write
\[
W_{u_0,\delta}(z)
:=
\left(
\frac{d(x,z)}{u_0}
+\frac{1}{\sqrt{u_0}}+1
\right)^\delta,
\]
and
 consider the case $i=1$ first. By the semigroup property of the spinor
heat kernel,
\[
\int_{(\mathbb H^d)^{k-1}}
P_{u_1}(z_1,z_2)
\cdots
P_{u_k}(z_k,x)
\,dz_2\cdots dz_k
=
P_{T_1}(z_1,x).
\]
Consequently, the integral in \eqref{eq:epsilon-zero-bound}, which we label by $I_1$ in what follows, is equal to
\[
\int_{\{d(x,z)\geq r_0\}}
W_{u_0,\delta}(z)
P_{u_0}(x,z)P_{T_1}(z,x)
\,dz.
\]
Using Equation 
\eqref{eq:scalarheatkernelbound2} and
Equation \eqref{eq:scalarheatkernelbound}, we obtain
\[
\begin{aligned}
|I_1|
&\leq
C
A_{u_0}A_{T_1}
\left(1+u_0^{-\delta/2}\right)
\\
&\qquad\qquad\cdot
\int_{\{d(x,z)\geq r_0\}}
\exp\left(
-\frac{d(x,z)^2}{cu_0}
-\frac{d(x,z)^2}{cT_1}
\right)
e^{-2m d(x,z)}
\,dz.
\end{aligned}
\]
Setting $H_1:=\frac{u_0T_1}{u_0+T_1}$,
we can now apply Lemma \ref{lemmageometricgaussian},
\begin{equation*}
    C
A_{u_0}A_{T_1}
\left(1+u_0^{-\delta/2}\right)
\sqrt{H_1}
\min \left\{\sqrt{H_1}, H_1^{\frac{d-1}{2}}\right\}
\exp\left(-\frac{c r_0^2}{H_1}\right),
\end{equation*}
which lets us conclude with Remark \ref{remupperboundmax}
and the fact that $H_1 \leq u_0$. \\
We now turn to the case $i>1$. By the semigroup property,
 the  integral in Equation \eqref{eq:epsilon-zero-bound}, denoted by $I_i$, reduces to
\[
\begin{aligned}
I_i
=
\int_{\mathbb H^d\times\mathbb H^d}
&
\mathbf 1_{\{d(z_i,x)\geq r_0\}}
W_{u_0,\delta}(z_1)
P_{u_0}(x,z_1)
\\
&\qquad\qquad\qquad\qquad
\cdot
P_{V_i-u_0}(z_1,z_i)
P_{T_i}(z_i,x)
\,dz_1\,dz_i.
\end{aligned}
\]
For fixed $z_i$, define
\[
F_i(z_i)
:=
\int_{\mathbb H^d}
W_{u_0,\delta}(z_1)
P_{u_0}(x,z_1)
P_{V_i-u_0}(z_1,z_i)
\,dz_1.
\]
Applying the same argument as above with total time $V_i$ gives
\begin{equation}
\begin{aligned}
|F_i(z_i)|
\leq
C
\left(1+u_0^{-\delta/2}\right)
\left(1+A_{V_i}\right)
\exp\left(-\frac{d(x,z_i)^2}{cV_i}\right)
e^{-m d(x,z_i)}, 
\end{aligned}
\label{eq:Fi-bound}
\end{equation}
which, together with Equation \eqref{eq:scalarheatkernelbound}, yields
\[
\begin{aligned}
|I_i|
&\leq
C
\left(1+u_0^{-\delta/2}\right)
\left(1+A_{V_i}\right)
A_{T_i}
\\
&\qquad\cdot
\int_{\{d(x,z_i)\geq r_0\}}
\exp\left(
-\frac{d(x,z_i)^2}{cV_i}
-\frac{d(x,z_i)^2}{cT_i}
\right)
e^{-2m d(x,z_i)}
\,dz_i.
\end{aligned}
\]
As before, set $H_i:=\frac{V_iT_i}{V_i+T_i}$.
Then
\[
\begin{aligned}
|I_i|
&\leq
C
\left(1+u_0^{-\delta/2}\right)
\left(1+A_{V_i}\right)
A_{T_i}
\\
&\qquad\cdot
\int_{\{d(x,z_i)\geq r_0\}}
\exp\left(-c\frac{d(x,z_i)^2}{H_i}\right)
e^{-2m d(x,z_i)}
\,dz_i \\
&\leq
C
\left(1+u_0^{-\delta/2}\right)
\left(1+A_{V_i}\right)
A_{T_i}
\int_{r_0}^{\infty}
\exp\left(-c\frac{r^2}{H_i}\right)\,dr \\
&\leq C \left(1+u_0^{-\delta/2}\right)
\left(1+A_{V_i}\right)
A_{T_i} \sqrt{H_i}
\exp\left(-c\frac{r_0^2}{V_i}\right),
\end{aligned}
\]
which, together with $A_{T_i}\sqrt{H_i}
\leq
A_{T_i}\sqrt{T_i}
=
T_i^{-(d-1)/2}+1$, lets us conclude.

\end{proof}
By a similar argument to the one that was used in the proof of Lemma
\ref{lemmaspinorchainbound1}, one obtains the following.

\begin{lemma}\label{lemmaspinorchainbound2}
    Let $u_j >0$, with $j \in \{0, \dots,k \}$, $t_0 = \sum_{j=0}^k u_j$, and 
    \begin{equation*}
    U_j :=
        \begin{dcases*}
            u_0 + \dots + u_{j-1}, & $j<i$, \\
            u_j + \dots + u_k, & $j \geq i$.
        \end{dcases*}
    \end{equation*}
Moreover, set $\epsilon >0$ and
$r_0 > 0$. 
Then, there exist constants $C,c >0$
such that for $i,j \in \{1, \dots, k\}$, and $x \in \mathbb{H}^d$,
\begin{equation}\label{eq:lemmaboundkernels2}
    \begin{split}
       & \bigg| \int_{(\mathbb{H}^d)^k} 
        \one_{ \{d(z_i,x) \geq r_0 \} }
        \left( \frac{d(z_j,x)}{U_j} + \frac{1}{\sqrt{U_j}} +1 \right)^{\epsilon} \\
        & \quad \quad \quad \cdot  P_{u_0}(x,z_1) 
     P_{u_1}(z_1,z_2)   \cdot ... \cdot  P_{u_k}(z_k,x)  
     \Id z \bigg| \\
     &\leq C_{c, \epsilon, d}
     \left(1+u_0^{-\epsilon/2}\right)
\left(1+A_{V_i}\right)
\left(1+T_i^{-(d-1)/2}\right)
\begin{dcases*}
            1, & $j<i$, \\
            1 + T_i^{-\epsilon/2}, & $j \geq i$,
        \end{dcases*} \cdot
\exp\left(-c\frac{r_0^2}{V_i}\right).
    \end{split}
\end{equation}
    
\end{lemma}

\subsection{Deriving the main result}

\begin{lemma}\label{shuffleheatkernellemma}
Let
$j \geq d-1$ and $u \in t \Delta_j$. Then
\begin{equation}\label{eq:shuffleheatkernellemma}
    \begin{split}
    & 
\int_{S_R(o)} Q^{x_0}_{u} \shuffle_* (ic(\Id B)(x_0))^{\otimes j}  \Id_{S_R(o)} x_0
  \\
          =&\int_{S_R(o)} P_{u}(x_0,x_0) \cdot \left[
  e^{-uA^2(x_0)} \shuffle (ic(\Id B)(x_0))^{\otimes j}  \right]\Id_{S_R(o)} x_0.
    \end{split}
\end{equation}

\end{lemma}
\begin{proof}
   We begin with the lowest order term from Equation \eqref{eq:shuffleheatkernellemma}  and 
    expand the integrand, for $j=d-1$, as 
\begin{equation}\label{parallelrep}
\begin{split}
   & \int_{(\mathbb{H}^d)^{d-1}}\prod_{i=0}^{d-1} p_{u_i}(y_i,y_{i+1})  \cdot   e^{-u_0 A^2(x_0)} \cdot \\
   &\prod_{j=1}^{d-1} \tau(y_{j-1},y_j)
    ic \left( (\Id B) (y_j)\right) 
    e^{-u_j A^2(x_0)} \cdot \tau(y_{d-1},x_0)
    e^{-u_{d-1}A^2(x_0)} \Id y,
    \end{split}
\end{equation}
with the convention that $y_d := y_0 := x_0$.  
In a first step, we argue that 
in Equation \eqref{parallelrep}, the term
$$\tau(x_0,y_1) ic(\Id B)(y_1)
\stackrel{\text{Remark  }\ref{taulc}}{=} ic(\tau^{*,LC}(x_0,y_1)\Id B(y_1))\tau(x_0,y_1), $$
can be replaced
by $ic(\Id B)(x_0) \tau(x_0,y_1)$, which we do by
estimating 
\begin{equation*}
    \begin{split}
       &\int_{S_R(o)} \int_{(\mathbb{H}^d)^{d-1}} 
        \prod_{i=0}^{d-1} p_{u_i}(y_i,y_{i+1}) 
     e^{-u_0 A^2(x_0)} \\
     &\bigg[ 
     \tau(x_0,y_1) ic(\Id B)(y_1) -
     ic(\Id B)(x_0) \tau(x_0,y_1) \bigg]e^{-u_1 A^2(x_0)} \tau(y_1,y_2) \\
     &\prod_{j=2}^{d-1} 
    ic \left( (\Id B) (y_j)\right)  
    e^{-u_j A^2(x_0)} \tau(y_{j},y_{j+1}) \cdot \tau(y_{d-1},x_0)
    e^{-u_{d-1}A^2(x_0)}  \Id y \Id_{S_R(o)} x_0.
    \end{split}
\end{equation*}
We now focus only on the integral over $y_1$
and partition
$\mathbb{H}^d$ via a dyadic decomposition, 
\begin{equation*}
    \mathbb{H}^d = \bigsqcup_{n=0}^{\infty}
\big(B_{r_n}(x_0)\backslash B_{r_{n-1}}(x_0)\big),
\end{equation*}
where $r_0(R)$ is a monotonously increasing 
function such that $\frac{r_0(R)}{\sqrt{R}} 
\rightarrow + \infty$ as $R \rightarrow \infty$
and 
$r_{n-1} =: r_n/2$.
Firstly, we treat the outer annuli.
This gives
\begin{equation}\label{equationdyadic}
    \begin{split}
        & \int_{S_R(o)}
        \sum_{n=1}^{\infty} 
        \int_{B_{r_n}(x_0)/B_{r_{n-1}}(x_0)}
   \bigg\Vert     p_{u_0}(x_0,y_1) p_{u_1}(y_1,y_2) \\
    &\bigg[ 
     \tau(x_0,y_1) ic(\Id B)(y_1) -
     ic(\Id B)(x_0) \tau(x_0,y_1) \bigg]
     \tau(y_1,y_2)  \bigg\Vert \Id y_1 
     \Id_{S_R(o)} x_0\\
   \leq  & \int_{S_R(o)} \sum_{n=1}^{\infty} 
        \exp{-C_1r_{n-1}^2 + C_2 r_{n-1}}
        \cdot C_3 \cdot \text{vol} B_{r_n}(x_0)
        \Id_{S_R(o)} x_0\\
    \leq & \int_{S_R(o)} C \cdot     \sum_{n=1}^{\infty}  
    \exp{-C_1' r_{n-1}^2 + d r_n}  \Id_{S_R(o)} x_0,
    \end{split}
\end{equation}
where, in the first inequality, we used 
the Kato-Simon inequality, see Remark \ref{heatestimate}, and the Gaussian bound
on the scalar heat kernel,
\begin{equation*}
    q_t(x,y) \leq C t^{-d/2} \exp{- \frac{\text{dist}(x,y)^2}{4t} + \frac{d-1}{2}
 \text{dist}(x,y)   }.
\end{equation*}
Since $r_{n-1} = r_n/2$, we can estimate the last 
line of Equation \eqref{equationdyadic} by
\begin{equation*}
\begin{split}
   &\int_{S_R(o)} C_3 e^{-C_4 r_{0}^2(R)} \Id_{S_R(o)} x_0 \\
   \leq & C_5 \text{vol}(S_R(o)) e^{-C_4 r_{0}^2(R)} \\
   \leq & C_6 \exp{(d-1)R - C_4 r_0^2(R)}  \quad 
   \xrightarrow{R\to\infty} 0,
   \end{split}
\end{equation*}
since, by assumption, $\frac{r_0(R)}{\sqrt{R}} 
\rightarrow + \infty$ as $R \rightarrow \infty$. \\
Remains to treat the inner ball of the dyadic 
decomposition. Again, separating the 
$y_1$-integration yields, together with
\begin{equation*}
    \sup_{x_0 \in S_R(o)} 
    \sup_{y_1 \in B_{r_0(R)}(x_0)} e^{R + \sqrt{R}} 
    \norm{\tau(x_0,y_1) ic(\Id B)(y_1) -
     ic(\Id B)(x_0) \tau(x_0,y_1)}
    \leq C < \infty,
\end{equation*}
for some $R$-independent constant $C$,
\begin{equation*}
    \begin{split}
        &\int_{S_R(o)}
        \int_{B_{r_0}(x_0)}  \bigg\Vert
        p_{u_0}(x_0,y_1)   p_{u_1}(y_1,y_2) \\
    &\bigg[ 
     \tau(x_0,y_1) ic(\Id B)(y_1) -
     ic(\Id B)(x_0) \tau(x_0,y_1) \bigg]
     \tau(y_1,y_2)  \bigg\Vert \Id y_1 
     \Id_{S_R(o)} x_0\\
     \leq & C \exp{- (R + \sqrt{R}) + dr_0 
     +(d-1)R}.
    \end{split}
\end{equation*}
Combining this with the remaining $d-2$-terms 
$ic(\Id B)(y_i)$ from Equation \eqref{parallelrep}
which all decay like $e^{-r}$ as 
$\vert y_i-o \vert = r \rightarrow \infty$,
as well as the assumption $\frac{r_0(R)}{\sqrt{R}} 
\rightarrow + \infty$ as $R \rightarrow \infty$, we are allowed to conclude.
Hence, Equation \eqref{parallelrep} becomes
\begin{equation*}
    \begin{split}
         &e^{-u_0 A^2(x_0)} 
    ic \left (\Id B)(x_0)\right) \tau(x_0,y_1)e^{-u_1 A^2(x_0)} \cdot \\     
         &\prod_{j=2}^{d-1} \tau(y_{j-1},y_j)
    ic \left( (\Id B) (y_j)\right) 
    e^{-u_{j} A^2(x_0)} \cdot \tau(y_{d-1},x_0)
    e^{-u_{d-1}A^2(x_0)},
    \end{split}
\end{equation*}
so that  iterating the procedure yields 
\begin{equation*}
\begin{split}
   &e^{-u_0 A^2(x_0)} 
    ic \left((\Id B)(x_0) \right) \tau(x_0,y_1)e^{-u_1 A^2(x_0)} \cdot \\     
         &\prod_{j=1}^{d-2} 
    ic \left( (\Id B) (y_j)\right) \tau(y_{j},y_{j+1})
    e^{-u_{j+1} A^2(x_0)} \cdot \tau(y_{d-1},x_0)
    e^{-u_{d-1}A^2(x_0)}.
    \end{split}
\end{equation*}
Note that there is now no term of 
$ic \left( (\Id B) (y_j)\right) $  that
depends on $y_{d-1}$ anymore, meaning that we can combine 
the parallel transports depending on $y_{d-1}$
with the corresponding scalar heat kernel parts 
$p_{u_{d-2}}(y_{d-2},y_{d-1})p_{u_{d-1}}(y_{d-1},x_0)$ 
as to produce $P_{u_{d-2} +u_{d-1}}(y_{d-2},x_0)$ 
via the semi-group property. 
Thus, we have integrated out $y_{d-1}$
and reduced the problem to $d-2$ variables.
Proceeding similarly by separating the $y_2$-integration, we eliminate the 
$y_{d-2}$  variable via
\begin{equation*}
\begin{split}
   &e^{-u_0 A^2(x_0)} 
    ic \left((\Id B)(x_0) \right) e^{-u_1 A^2(x_0)}ic \left((\Id B)(x_0) \right)
    \tau(x_0,y_1)
    e^{-u_2 A^2(x_0)} \cdot \\     
         &\prod_{j=1}^{d-3} 
    ic \left( (\Id B) (y_j)\right) \tau(y_{j},y_{j+1})
    e^{-u_{j+1} A^2(x_0)} \cdot P_{u_{d-2} +u_{d-1}}(y_{d-2},x_0)
    e^{-u_{d-1}A^2(x_0)},
    \end{split}
\end{equation*}
so that iteratively, we end up with
\begin{equation*}
    P_s(x_0,x_0) \cdot \left[ 
  e^{-uA^2(x_0)} \shuffle (ic(\Id B)(x_0))^{\otimes (d-1)} \right],
\end{equation*}
which allows us to conclude. 
The higher order terms 
$j >d-1$ then follow by the same argument.

\end{proof}

\begin{lemma}\label{lemma6.2}
    For $j \in \{1, \dots, d \}$ and $t >0$,
   \begin{equation}
        \begin{split}
     &\lim_{R \rightarrow \infty} 
        \sinh^{d-1} R \lim_{\epsilon \rightarrow 0}
        \int_0^t \int_{S_1(o)} \\
&\norm{ \left( \partial_{j,1} \left(P_{\epsilon
} \otimes P_{\epsilon} * r_s^{d-2}\right) \right)(\exp_o(Ry), 
        \exp_o(Ry) ) }_{S^1(H)}\Id S_1(y) \Id s =0.        
        \end{split}
\end{equation}
\end{lemma}
\begin{proof}
We proceed as in \cite{oliver23}, Lemma 5.8.
    For $x_0 \in \mathbb{H}^d$, and 
    $\gamma \in \{0,1\}^k$, define 
   $$ S^{x_0}_{\gamma} := (S_{\gamma_1}, \dots, S_{\gamma_k}), \quad \text{with}
   \quad S_0 := A^2 - A^2(x_0), \: S_1 := 
   ic \nabla B.$$
Then, for a.e.  $x_0, x,y \in \mathbb{H}^d$, 
we can rewrite $r_s^{j}$ as 
\begin{equation}
    r^{j}_s(x,y) = - \sum_{k = j}^{\infty}
    (-s)^k \sum_{\substack{\gamma \in \{0,1\}^k \\ \vert \gamma \vert = j}} 
    \int_{u \in \Delta_k} \big( Q_{su}^{x_0}\shuffle_* S_{\gamma}^{x_0} \big)(x,y) 
    \Id u, 
\end{equation}
where $j \geq d-2$ is odd, so that
\begin{equation}\label{eq:serieskappa}
    \kappa_s^{-,d-2}(x,y) = - \sum_{j = \frac{d-3}{2}}^{\infty}
    \sum_{k = 2j+1}^{\infty}
    (-s)^k \sum_{\substack{\gamma \in \{0,1\}^k \\ \vert \gamma \vert = 2j+1}} 
    \int_{u \in \Delta_k} \big( Q_{su}^{x_0}\shuffle_* S_{\gamma}^{x_0} \big)(x,y) 
    \Id u, 
\end{equation}
which converges absolutely in $C_{m,n}^{1,S^1}(\mathbb{H}^{2d})$ for $m,n \in \mathbb{N}_0$.
Here, $ Q_t^{x_0}(x,y) :=  e^{-tA^2(x_0)} 
    \otimes P_t(x,y)$, $x,y \in \mathbb{H}^d$,
denotes the smooth integral kernel associated to
the semi-group $e^{-t(\slashed{\Delta} + M_{x_0}^2)}$, with $M_{x_0}$ being multiplication by
$A(x_0)$.
We now use Equation \eqref{eq:serieskappa} to study 
$((\partial_j P_{\epsilon} \otimes P_{\epsilon}) * \kappa^{-,d-2}_s)(\exp_o(Ry), \exp_o(Ry))$, where
$R>0$, $y \in S_1(o)$ and $x_0 = \exp_o(Ry)$.
From 
$\norm{Q_t^{x_0}}_{B(S \otimes \text{Dom}(A_0^n))} \leq  P_t$, 
we obtain
\begin{equation}\label{eq:boundonderivativeshuffle}
    \begin{split}
        &\norm{\int_{u \in s \Delta_k} \partial_{j,1} ( P_{\epsilon} \otimes P_{\epsilon}
        * (Q_u^{x} \shuffle_* S_{\gamma}))(x,x) \Id u
        }_{S^1(H)} \\
        \leq& s^k \bigg\vert \int_{s \in \Delta_k} \int_{(\mathbb{H}^d)^k} 
        \left( \frac{\text{dist}(z_1,o)}{ su_0 + \epsilon} + 
        \frac{1}{\sqrt{su_0+ \epsilon}} +1 \right) 
         P_{su_0 + \epsilon}(z_1)  \\
        & \quad \cdot P_{su_1}(z_1, z_2) \cdot \dots \cdot 
         P_{su_k + \epsilon}(z_k) \prod_{i = 1}^k \norm{S_{\gamma_i}^x(x+z_i)}_{S^1(H)} \Id z \Id u \bigg\vert.
    \end{split}
\end{equation}
To estimate the integral, we restrict from 
$z \in \mathbb{H}^d$
to $z \in B_{r_0}(o)^k$ 
where $r_0 = \log d(x,o)$ and consider two cases. 
Using $\norm{S^x_{\gamma_i}}_{L^{\infty}(\mathbb{H}^d,S^1(H))} \leq c$, we can 
estimate the $2^k-1$ terms where at least one
$z_i \notin B_{r_0}(o)$ by Lemmas
\ref{lemmaspinorchainbound1} and 
\ref{lemmaspinorchainbound2} for 
$\epsilon \in (0,1)$, hence allowing the factor of $1+ T_i^{-\epsilon/2}$ to be absorbed 
into $1+ T_i^{-(d-1)/2}$ and thus eliminating the possible singularity for the integral over $s \in (0,1)$, as well as providing summability over $k$.
Now, for the integral over $z \in B_{r_0}(o)^k$, we combine Equation \eqref{hyp4.1} of Hypothesis 
\ref{hypothesis} as by Remark 
\ref{rem:hypothesisapprox}
and the decay 
$\norm{\Id B(x)}_{S^1(H)} \sim e^{-d(x,o)} $
as $d(x,o) \rightarrow \infty$, to show absolute summability of \eqref{eq:boundonderivativeshuffle}
with respect to both sums in Equation \eqref{eq:serieskappa}, as well
as integrability in $s $ at zero, uniformly
for $\epsilon \in (0,1)$.

Due to the limit $d(x,o) \rightarrow \infty$, we may neglect all terms in 
Equation \eqref{eq:serieskappa}
with $j \geq \frac{d-1}{2}$ and $j = \frac{d-3}{2}$, $k \geq d$, and 
it suffices, for the Lemma to hold,  to argue that
 for $\gamma \in \{0,1\}^{d-1}$,
$\vert \gamma \vert = d-2$,
\begin{equation}\label{eqd21}
\begin{split}
    &\lim_{R \rightarrow \infty} 
    \sinh^{d-1} R \lim_{\epsilon \rightarrow 0} \int_0^t s^{d-2} \int_{S_1(o)} \\
   &\norm{\int_{u \in 
    \Delta_{d-2}} 
    \partial_{j,1} 
    \big( P_{\epsilon} \otimes P_{\epsilon} *(Q_{su}^{\exp_o(Ry)}\shuffle_* (ic \Id B)^{\otimes (d-2)}\big)(\exp_o(Ry), 
    \exp_o(Ry)) \Id u }_{S^1(H)}
    \\
    &\quad \quad \quad \quad \quad \quad \quad \quad \quad\quad \quad \quad \quad \quad \quad \quad \quad \quad \quad \quad \quad \quad \quad \quad \quad \quad \quad \quad \quad \quad \Id S_1(y) =0,
    \end{split}
\end{equation}
and 
\begin{equation}\label{eqd22}
\begin{split}
    &\lim_{R \rightarrow \infty} 
    \sinh^{d-1} R \lim_{\epsilon \rightarrow 0} \int_0^t s^{d-1} \int_{S_1(o)} \\
   &\norm{\int_{u \in 
    \Delta_{d-1}} 
    \partial_{j,1} 
    \big( P_{\epsilon} \otimes P_{\epsilon} *(Q_{su}^{\exp_o(Ry)}\shuffle_* S_{\gamma})\big)(\exp_o(Ry), 
    \exp_o(Ry)) \Id u }_{S^1(H)}
    \\
    &\quad \quad \quad \quad \quad \quad \quad \quad \quad\quad \quad \quad \quad \quad \quad \quad \quad \quad \quad \quad \quad \quad \quad \quad \quad \quad \quad  \Id S_1(y) =0.
    \end{split}
\end{equation}
Let us start with Equation \eqref{eqd21}.
We have 
\begin{equation*}
    \begin{split}
 &\partial_{j,1} 
    \big( Q_{su}^{\exp_o(Ry)}\shuffle_* (ic \Id B)^{\otimes (d-2)}\big)(\exp_o(Ry), 
    \exp_o(Ry))     \\
    &=  \int_{(\mathbb{H}^d)^{d-2}} 
    \partial_{j,1} (P_{su_0}(\exp_o(Ry), w_1))
    p_{su_1}(w_1,w_2) \cdot ... \cdot  
    p_{su_{d-2}}(w_{d-2}, \exp_o(Ry)) \\
    &  e^{-su_0 A^2(\exp_o(Ry))} 
    \prod_{j=1}^{d-3} 
    (ic\Id B(w_j)) \tau(w_j, w_{j+1})
    e^{-su_j A^2(\exp_o(Ry))} \\
    & (ic\Id B(w_{d-2})) \tau(w_{d-2}, \exp_o(Ry))
    e^{-su_{d-2} A^2(\exp_o(Ry))} \Id w, \\
\end{split}
\end{equation*}
which yields, upon changing variables 
from $w$ to $L_g \exp_o(Ry)$,
\begin{equation*}
\begin{split}
    & 
    \int_{(\mathbb{H}^d)^{d-2}} 
    \partial_{j,1} (P_{su_0}(g_1,o))
    p_{su_1}(g_1,g_2) \cdot ... \cdot  
    p_{su_{d-2}}(g_{d-2}, o) \\
    & \cdot e^{-su_0 A^2( \exp_o(Ry))}
    \prod_{j=1}^{d-3} (ic\Id B(L_{g_j}\exp_o(Ry))) \\
& \cdot \tau(L_{g_j} \exp_o(Ry), L_{g_{j+1}} \exp_o(Ry))    e^{-su_j A^2(\exp_o(Ry))} \\
    & \cdot (ic\Id B(L_{g_{d-2}} \exp_o(Ry))) \tau(
    L_{g_{d-2}} \exp_o(Ry), \exp_o(Ry))
    e^{-su_{d-2} A^2(\exp_o(Ry))} \Id g\\
    \stackrel{L_{g} \mapsto L_{g^{-1}}}{=} &
    \int_{(\mathbb{H}^d)^{d-2}} 
    \partial_{j,1} (P_{su_0}(z_1,0))
    p_{su_1}(z_1,z_2) \cdot ... \cdot  
    p_{su_{d-2}}(z_{d-2}, 0) \\
    & e^{-su_0 A^2( \exp_o(Ry))}
    \prod_{j=1}^{d-3} (ic\Id B(L_{g_j^{-1}}\exp_o(Ry))) \\
& \tau(L_{g_j^{-1}} \exp_o(Ry), L_{g_{j+1}^{-1}} \exp_o(Ry))    e^{-su_j A^2(\exp_o(Ry))} \\
    &(ic\Id B(L_{g_{d-2}^{-1}} \exp_o(Ry))) \tau(
    L_{g_{d-2}^{-1}} \exp_o(Ry), \exp_o(Ry))
    e^{-su_{d-2} A^2(\exp_o(Ry))} \Id g
    \end{split}
\end{equation*}
and hence, with the notations from 
Definition \ref{deltadef},
\begin{equation*}
    \begin{split}
        &\sinh^{d-1}( R) \partial_{j,1} \bigg( P_{\epsilon} \otimes 
        P_{\epsilon} * \bigg(
         Q_{su}^{\exp_o(Ry)}\shuffle_* (ic \Id B)^{\otimes (d-2)}
         \\
        & \quad -   Q_{su}^{\exp_o(Ry)}\shuffle_* (ic \Id B(\exp_o(Ry))^{\otimes (d-2)} \bigg) \bigg)
        (\exp_o(Ry), \exp_o(Ry))
    \end{split}
\end{equation*}
\begin{equation*}
    \begin{split}
        &= \sinh^{d-1}( R)\int_{(\mathbb{H}^d)^{d-2}} 
        \partial_{j,1} (P_{su_0 + \epsilon}(z_1,0))
    p_{su_1}(z_1,z_2) \cdot ... \cdot  
    p_{su_{d-2} + \epsilon}(z_{d-2}, 0)    e^{-su_0A^2(\exp_o(Ry))} \\
    & \textnormal{OpProd}\left(\tau_{\big(z/2\big)^{-1}}
        - \tau_{\big(z/2\big)}\right)^2  
    \bigg( ic\Id B 
    \tau
    e^{-su_1A^2(\exp_o(Ry))},\\
   & \qquad \qquad \qquad \qquad \qquad \qquad \qquad \qquad..., ic \Id B \tau e^{-su_{d-2}A^2(\exp_o(Ry))} \bigg) (\exp_o(Ry)) \Id z,
    \end{split}
\end{equation*}
 which suffices to be studied by Lemma
\ref{shuffleheatkernellemma}.
Due to the small-time asymptotics of the 
hyperbolic spinor heat kernel, see Remark 
\ref{rem:smalltimespinorkernel},
 and restricting to $z \in B_{r_0}(o)^{d-2}$, as was done in the argument above, we 
obtain, similar to the Euclidean case 
\cite{oliver23}, that the left hand side 
of Equation \eqref{eqd21} is equivalent to
\begin{equation}
\begin{split}
     &- \lim_{R \rightarrow \infty} \sinh^{d-1} R
    \int_0^t s^{d-2} \int_{S_1(o)} 
    \int_{u \in \Delta_{d-2}} 
    \int_{z \in B_{r_0}(o)^{d-2}} \\
   & \partial_{j,1} (P_{su_0 }(z_1,0))
    p_{su_1}(z_1,z_2) \cdot ... \cdot  
    p_{su_{d-2} }(z_{d-2}, 0)  \\
    &\Vert e^{-su_0A^2(\exp_o(Ry))} \textnormal{OpProd}\left(\tau_{\big(z/2\big)^{-1}}
        - \tau_{\big(z/2\big)}\right)^2  
    \bigg( ic\Id B 
    \tau
    e^{-su_1A^2(\exp_o(Ry))},\\
   & \qquad \qquad \qquad \qquad \qquad \qquad \qquad \qquad..., ic \Id B \tau e^{-su_{d-2}A^2(\exp_o(Ry))} \bigg) (\exp_o(Ry)) \Vert_{S^1(H)} \\
   &\Id z \Id u \Id S_1(o) \Id s.
\end{split}
\end{equation}      
We now perform  the decomposition  of $\textnormal{OpProd}\left(\tau_{\big(z/2\big)^{-1}}
        - \tau_{\big(z/2\big)}\right)^2 $ in the sense of Lemma \ref{opproddec}.
   The products incorporating  a factor of 
   $\mathcal{W}^2_\cdot$ vanish as 
   $R \rightarrow \infty$ due to Equations \eqref{hyp4.3} and \eqref{hyp2}
   from Hypothesis 
   \ref{hypothesis}, after an application 
   of the dominated convergence Theorem.
   On the other hand, the products involving 
   two occurrences of $\delta_\cdot$ 
   go to zero as $R \rightarrow \infty$ 
   as a consequence of Equation \eqref{hyp4.2} from Hypothesis \ref{hypothesis} and, again, the dominated convergence Theorem. In three dimensions, it is then necassary to apply Hypothesis \eqref{hyp4.2} to the approximated family $B_{n,m,l}$. \\
   Let us now turn our attention to Equation
\eqref{eqd22}.
Take $S_{\gamma}:= S_{\gamma}^{\exp_o(Ry)}$ with a factor of 
$A^2 - A^2(\exp_o(Ry))$ in the 
$j$-th slot and denote by
$\tilde{S}_{\gamma}$ the same $S_{\gamma}^{\exp_o(Ry)}$, but
with all $ic \Id B$ replaced by 
$ic \Id B(\exp_o(Ry))$. 
Since the spinor heat kernel behaves 
like its Euclidean counterpart when
$s \rightarrow 0$, we have, uniformly in 
$R$ and $y$, 
\begin{equation*}
    s^{d-1} \norm{\int_{u \in \Delta_{d-1}}
 Q_{su}^{\exp_o(Ry)}\shuffle_* S_{\gamma}   \Id u }_{C_{0,0}^{1,S^1}(\mathbb{H}^{2d})}
 \leq C s^{\frac{d-3}{2}},
\end{equation*}
whenever $s$ is close to zero.
Hence we may commute the limit 
$\epsilon \rightarrow 0$ with the integrals.

After changing coordinates in the convolution,
employing the dominated convergence theorem
and making use of Equations \eqref{hyp4.1} and
\eqref{hyp4.2} from Hypothesis 
\ref{hypothesis}, we see that
\begin{equation}
    \begin{split}
        &\lim_{R \rightarrow \infty} 
        \sinh^{d-1}R \int_0^t s^{d-1} \int_{u \in \Delta_{d-1}} \bigg\vert \bigg\vert
        \partial_{j,1} \bigg(
           \big( Q_{su}^{\exp_o(Ry)}\shuffle_* S_{\gamma}\big)\\
        & \qquad \qquad \qquad -  
           \big( Q_{su}^{\exp_o(Ry)}\shuffle_* \tilde{S_{\gamma}}\big) \bigg)
        (\exp_o(Ry), \exp_o(Ry)) \bigg\vert 
        \bigg\vert_{S^1(H)}
        \Id u \Id s =0.
    \end{split}
\end{equation}
Let now $T$ denote the $(d-1)$-vector of operator
valued functions, containing 
$(d-2)$-times the constant operator 
$i c \Id B(\exp_o(Ry))$ and one occurrence of $A^2 - A^2(\exp_o(Ry))$ in the $j$-th
slot. Then, after successive applications 
of Lemma \ref{shuffleheatkernellemma},
\begin{equation}\label{helper}
    \begin{split}
 &\partial_{j,1}
           \big( Q_{su}^{\exp_o(Ry)}\shuffle_* \tilde{S_{\gamma}}\big)
        (\exp_o(Ry), \exp_o(Ry))    \\
    &= \int_{\mathbb{H}^d}  
   \big( \partial_{j,1} P_{s(u_0 + \dots + u_{j-1})} \big)(\exp_o(Ry), w_j) 
 P_{s(u_j + \dots + u_{d-1})}(w_j, \exp_o(Ry))\\
 & \qquad \qquad \cdot \bigg( 
 e^{-suA^2(\exp_o(Ry))} \shuffle T(w_j) \bigg) \Id w_j,
    \end{split}
\end{equation}
and, after
changing coordinates,
as in the argument above, to $w_j \mapsto L_{g_j} \exp_o(Ry)$
and then to $L_{g_j} \mapsto L_{g_j^{-1}} $, we obtain a factor of $\mathcal{W}^2_{g_j}(A^2)
(\exp_o(Ry))$ in the $j$-th slot of 
the shuffle product of Equation
\eqref{helper}. 
Consequently, Equation \eqref{eqd22} follows
from dominated convergence and Equation 
\eqref{hyp4.3} from Hypothesis \ref{hypothesis}.

\end{proof}
We are now ready to state a first trace formula. 
In what follows, for $v \in  T_x  \mathbb{H}^d$ 
and $x \in \mathbb{H}^d$, 
we set $\langle c,v \rangle := \sum_{j=1}^d c^j v_j$.
The proof is completely analogous to the proof of Proposition 5.9 in \cite{oliver23}
and works by exploiting the results established in Lemma \ref{lemmarholimit}, 
Theorem \ref{mercer}, Corollary 
\ref{corollary29}, and Lemma \ref{lemma6.2}.

\begin{proposition}\label{prop6.5}
    Let $t>0$. Then 
    \begin{equation}
    \begin{split}
&\textnormal{tr}_{L^2(\mathbb{H}^d,H)}
    \textnormal{tr}_{S} \left(
    e^{-tD^*D} - e^{-tDD^*}\right) \\
   & = \textnormal{tr}_H
   \lim_{ r \rightarrow  \infty} 
   \sinh^{d-1}{r} 
   \int_0^t 
   \int_{S_1(o)}  \\
& \qquad \qquad \qquad i \langle c, \nu_x \rangle \bigg(A(\exp_o(ry)) \kappa_s^{+,d-1}(
   \exp_o(ry),\exp_o(ry)) \bigg) \Id S_1(y) \Id s,\\
    \end{split}
\end{equation} 
where $\nu_x = \frac{d}{dr} \exp_o(ry) \bigg\vert_{r = \textnormal{dist}(o,x)}$ is the unit normal
vector at a point $x = \exp_o(ry)$, $y \in S^{d-1}(o)$ and  with the limit existing 
in $S^1(H)$-norm.
\end{proposition}

\begin{lemma}\label{lemma3}
For arbitrary but fixed $t>0$,
\begin{equation}\label{lemma3eq}
    \begin{split}
&\textnormal{tr}_{L^2(\mathbb{H}^d,H)}
    \textnormal{tr}_{S} \left(
    e^{-tD^*D} - e^{-tDD^*}\right) \\
   & = \textnormal{tr}_H
   \lim_{ r \rightarrow  \infty} 
   \sinh^{d-1}{r} 
   \int_0^t 
   \int_{S_1(o)}  \\
& i \langle c, \nu_x \rangle \bigg(A(\exp_o(ry)) J_s^{+,d-1}(
   \exp_o(ry),\exp_o(ry)) \bigg) \Id S_1(y) \Id s, \\
    \end{split}
\end{equation} 
where 
\begin{equation}
   J_s^{+,d-1} := \sum_{j=\frac{d-1}{2}}^{\infty} \Xi_s^{2j}, 
\end{equation}
and 
\begin{equation}\label{xiterm}
\begin{split}
    &\Xi^{j}_t(x_0) :=   
\int_{u \in t \Delta_{d-1}} 
 \left(   Q_u^{x_0} \shuffle_* (ic \Id B)^{\otimes j} 
 \right)(x_0,x_0) \Id u,
    \end{split}
\end{equation}
for $Q_t^{x_0}(x,y) :=  e^{-tA^2(x_0)} 
    \otimes P_t(x,y)$ and $x,y,x_0 \in \mathbb{H}^d$.
    
\end{lemma}
\begin{proof}
 Pick $x_0 \in \mathbb{H}^d$
and  denote by $M_{x_0}$ the multiplication operator
$(M_{x_0} f)(x) := A(x_0)f(x)$, for $x \in \mathbb{H}^d$
and $f \in \text{Dom}(M_0)$. 
Then $H_{x_0} := \slashed{\Delta} + M^2_{x_0}$ defines 
a self-adjoint operator on $\text{Dom}(H_0)$. 
The associated semi-group, $e^{-tH_{x_0}}$, $t>0$, 
admits a smooth integral kernel given by 
\begin{equation}
    Q_t^{x_0}(x,y) :=  e^{-tA^2(x_0)} 
    \otimes P_t(x,y), \quad 
    x,y \in \mathbb{H}^d.
\end{equation}
Denoting the norm $\norm{\cdot}_{\text{End}(S_x) \otimes\langle A_0 \rangle^{-n}B(H) \langle A_0 \rangle^{-m}}$ by $\norm{\cdot}_{\Box}$, we have 
\begin{equation}\label{eq6}
\begin{split}
   & \limsup_{R \rightarrow \infty}
   \int_0^t \int_{S_R(o)} \int_{s \in \Delta_{d-1}} 
    \int_{\mathbb{H}^d} \\
&\norm{(\theta_{ts_0}(x_0,y) - Q_{ts_0}^{x_0}(x_0,y))(ic\Id B)(y) (Q_{t(s_1, \dots,s_{d-1})} \shuffle (ic\Id B)^{\otimes(d-2)})(y,x_0)}_{\Box} \\
&\Id_{\mathbb{H}^d} y \Id s \Id_{S_R(o)} x_0 \Id u \\
 \leq & \: \limsup_{R \rightarrow \infty} \:
C \int_0^t \int_{S_1(o)} \int_{s \in \Delta_{d-1}} 
    \int_{\mathbb{H}^d} \\
&\norm{\theta_{ts_0}(\exp_o(R\xi),y) - Q_{ts_0}^{\exp_o(R\xi)}(\exp_o(R\xi),y)}_{\Box} \:
e^{-\frac{\text{scal}(\mathbb{H}^d)}{4}t(1-s_0)} \\
&q_{t(1-s_0)}(y,\exp_o(R\xi))  
\Id_{\mathbb{H}^d} y \Id s \Id_{S_1(o)} \xi \Id u,
\end{split}
\end{equation}
where the inequality follows 
from the assumption that 
$$\norm{ic(\Id B)(x)}_{\text{op}} \leq C' \: 
\text{Vol}_{\mathbb{H}^d}(S_R(o))^{-\frac{1}{d-1}}$$ 
as $\vert x-o \vert =R \rightarrow \infty$, which 
neutralizes the $\sinh^{d-1}(R)$ term that
appears after rescaling to the unit sphere,
and Remark \ref{heatestimate}.
In Equation \ref{eq6}, we  
use Duhamel's 
principle to obtain 
\begin{equation}
    e^{-tH_B} = e^{-tH_{x_0}} - \int_0^t 
    e^{-sH_{x_0}}(A^2-M^2_{x_0})e^{-(t-s)H_B} \Id s,
\end{equation}
and hence 
the decomposition
\begin{equation}
\begin{split}
 &\theta_{ts_0}(x,y) - Q_{ts_0}^{x_0}(x,y) \\
 &= - \int_0^{s_0} \int_{\mathbb{H}^d} 
 Q_{tu}^{x_0}(x,w)\left(A^2(w) - A^2(x_0)\right) 
 \theta_{t(s_0-u)}(w,y) \Id w \Id u.
 \end{split}
\end{equation}
The term $ Q_{tu}^{x_0}(\exp_o(R\xi),w)$ will then 
produce the factor $e^{tu\frac{d(d-1)}{4}} q_{tu}(\exp_o(R\xi),w)$, so that Equation \ref{eq6} becomes
\begin{equation}
\begin{split}
& \: \limsup_{R \rightarrow \infty} \: C e^{-\frac{\text{scal}(\mathbb{H}^d)}{4}t} \int_0^t
\int_{S_1(o)}  
 \int_{\mathbb{H}^d} \int_{s \in \Delta_d} \\
 &q_{ts_0}(\exp_o(R\xi),w) \norm{A^2(w) - A^2(\exp_o(R\xi))}_{B(H,\text{Dom}(A_0))} 
 q_{t(1-s_0)}(w,\exp_o(R\xi)) \\
 &\Id s \Id_{\mathbb{H}^d} w \Id_{S_1(o)} \xi \Id u. 
\end{split}
\end{equation}
We now change variables from 
$w$ to $L_{z} \cdot \exp_o(R\xi)$, where 
$L_z \in G$ denotes the element 
in $G$ that satisfies $L_z \cdot o = z$ and  
accounts for the translation by $z \in \mathbb{H}^d$,
where $z = \exp_o(s\xi)$ for some $s \geq 0$,
so that $q_{ts_0}(\exp_o(R\xi), g_z \cdot\exp_o(R\xi)) 
= q_{ts_0}(\vert z \vert)$, where 
$\vert z \vert := \text{dist}_{\mathbb{H}^d}(o,z)$,
\begin{equation}\label{eq6.3}
\begin{split}
 & \: \limsup_{R \rightarrow \infty} \: Ce^{-\frac{\text{scal}(\mathbb{H}^d)}{4}t} \int_0^t
 \int_{S_1(o)} 
    \int_{\mathbb{H}^d} \int_{s \in \Delta_{d}} \\
&q_{ts_0}(\exp_o(R\xi), L_z \cdot\exp_o(R\xi)) q_{t(1-s_0)}(\exp_o(R\xi), L_z \cdot\exp_o(R\xi)) \\
 &\norm{A^2(L_z \cdot \exp_o(R\xi)) - A^2(\exp_o(R\xi))}_{B(H,\text{Dom}(A_0))} \\
 &\Id s \Id_{\mathbb{H}^d} z \Id_{S_1(o)} \xi \Id u\\
=&  \: \limsup_{R \rightarrow \infty} \: Ce^{-\frac{\text{scal}(\mathbb{H}^d)}{4}t} \int_0^t
 \int_{S_1(o)} 
    \int_{\mathbb{H}^d} \int_{s \in \Delta_{d}} \\
&q_{ts_0}(\vert z \vert) q_{t(1-s_0)}(\vert z \vert) 
 \norm{A^2(L_z \cdot \exp_o(R\xi)) - A^2(\exp_o(R\xi))}_{B(H,\text{Dom}(A_0))} \\
 &\Id s \Id_{\mathbb{H}^d} z \Id_{S_1(o)} \xi \Id u.
    \end{split}
\end{equation}
Since the factors involving the heat
kernels are independent of $R$, the integral admits 
an integrable dominant given that 
\begin{equation}\label{eq5}
\norm{A^2(g_z \cdot \exp_o(R\xi)) - A^2(\exp_o(R\xi)}_{B(H,\text{Dom}(A_0))}
\end{equation}
is uniformly bounded, allowing us to take
the limit $R \rightarrow \infty$ inside of the
integral. Since Equation \eqref{eq5} 
goes to zero as $R$ goes to infinity, so does the 
integral in the last line of Equation \eqref{eq6.3}.
We have therefore shown that 
\begin{equation}
\begin{split}
    &\limsup_{R \rightarrow \infty} 
    \int_{S_R(o)} \int_{s \in t \Delta_{d-1}}   
   \left( Q_{s}^{x_0} \shuffle_* (ic\Id B)^{\otimes(d-1)}
   \right)
   (x_0,x_0) \Id s \Id x_0 \\
&=\limsup_{R \rightarrow \infty} 
    \int_{S_R(o)} \int_{s \in t\Delta_{d-1}}   \\
    &\left(\left(\theta_{s_0},Q^{x_0}_{(s_1, \dots, s_{d-1})}\right)
    \shuffle_* (ic\Id B)^{\otimes(d-1)}\right)(x_0,x_0) 
    \Id s \Id x_0.
   \end{split}
\end{equation}
Repeating the procedure yields then
\begin{equation*}
\begin{split}
    &\lim_{r\rightarrow \infty} \sinh^{d-1}(r)
    \int_0^t \int_{S_1(o)}
    \text{tr}_{S_x} \big(i \langle c, \nu_x \rangle A(\exp_o(ry)) \\
   & \quad \qquad \quad \quad \qquad \qquad \qquad
   r^{d-1}_s(
    \exp_o(ry), \exp_o(ry)) \big)  \Id S_1(y) \Id s\\
&= \lim_{r \rightarrow \infty} \sinh^{d-1}(r)
\int_0^t 
\int_{S_1(o)} \text{tr}_{S_x} 
\left(i \langle c, \nu_x \rangle A(\exp_o(ry)) \:
\Xi_s(\exp_o(ry)) \right) \Id S_1(y) \Id s,
    \end{split}
\end{equation*}
where for $t>0$ and $x \in \mathbb{H}^d$, 
\begin{equation}
    \Xi^{d-1}_t(x_0) := \int_{u \in t\Delta_{d-1}}
 \left(   Q^{x_0}_u \shuffle_*(ic\Id B)^{\otimes (d-1)}
 \right) (x_0,x_0)
    \Id u.
\end{equation}
\end{proof}
Combining the preceding Lemma with Lemma 
\ref{shuffleheatkernellemma}, we are lead to the following.
\begin{proposition}\label{eqbeforemain}
    For arbitrary but fixed $t>0$,
\begin{equation}\label{premaineq}
    \begin{split}
&\textnormal{tr}_{L^2(\mathbb{H}^d,H)}
    \textnormal{tr}_{S} \left(
    e^{-tD^*D} - e^{-tDD^*}\right) \\
   & = \textnormal{tr}_H
   \lim_{ r \rightarrow  \infty} 
   \sinh^{d-1}{r} 
   \int_0^t 
   \int_{S_1(o)}  \\
& i \langle c, \nu_x \rangle \bigg(A(\exp_o(ry)) \tilde{\kappa}_s^{+,d-1}(
   \exp_o(ry),\exp_o(ry)) \bigg) \Id S_1(y) \Id s, \\
    \end{split}
\end{equation} 
where 
$$\tilde{\kappa}_s^{+,d-1} := \sum_{j=\frac{d-1}{2}}^{\infty} \tilde{r}_s^{2j},$$
\begin{equation}
   \tilde{r}_s^{j}(x_0) := 
   P_s(x_0,x_0)
   \int_{u \in t \Delta_{d-1}} 
  e^{-uA^2(x_0)} \shuffle (ic(\Id B)(x_0))^{\otimes j} 
  \Id u.
\end{equation}
\end{proposition}

We now argue that the only factor that remains
of $J_s^{+,d-1}$ is the one of lowest order.
\begin{lemma}\label{highervanish}
Let $t>0$. Then 
\begin{equation}\label{prelformula}
    \begin{split}
&\textnormal{tr}_{L^2(\mathbb{H}^d,H)}
    \textnormal{tr}_{S} \left(
    e^{-tD^*D} - e^{-tDD^*}\right) \\
   & = \textnormal{tr}_H
   \lim_{ r \rightarrow  \infty} 
   \sinh^{d-1}{r} 
   \int_0^t 
   \int_{S_1(o)}  \\
 & \qquad \qquad i \langle c, \nu_x \rangle \bigg(A(\exp_o(ry)) \tilde{r}_s^{d-1}(
   \exp_o(ry),\exp_o(ry))  \bigg) \Id S_1(y) \Id s.\\
    \end{split}
\end{equation} 
\end{lemma}
\begin{proof}
   Recall from Lemma \ref{eqbeforemain} that
   \begin{equation*}
   \tilde{\kappa}_s^{+,d+1} := \sum_{j=\frac{d+1}{2}}^{\infty} \tilde{r}_s^{2j},
\end{equation*}
where 
\begin{equation*}
\begin{split}
    &\tilde{r}_s^j(x) :=   
P_s(x,x)  \int_{u\in s\Delta_{d-1}} e^{-uA^2(x)} 
 \shuffle (ic\Id B(x))^{\otimes j} \Id s.
    \end{split}
\end{equation*}
From Hypothesis \ref{hypothesis} we have 
that
\begin{equation*}
    \sup_{x \in S_R(o)} e^{R } \norm{ic(\Id B)(x)}   \leq C < \infty,
\end{equation*}
for some constant $C$ not depending on $R$,
hence
\begin{equation*}
\begin{split}
    \norm{\tilde{\kappa}_s^{+,d+1} } 
    &\leq C_d \cdot \text{vol}(S_R(o)) \sum_{j = \frac{d+1}{2}}^{\infty} e^{-2jR} C^{2j}
   \frac{ s^{2j}}{(2j)!}\\
  & \leq C_d e^{(d-1)R} e^{-(d+1)R} e^{Cs}
  \longrightarrow 0, \quad \text{as} \: R \rightarrow \infty,
   \end{split} 
\end{equation*}
which allows us to conclude. \\
\end{proof}

\begin{proposition}\label{prop:xnbounded}
Set 
\begin{align}
 g_1(x)&:=\sup_{|v|_g=1}\|\Id B(x)(v)\langle A_0 \rangle^{-\alpha}\|_{S^d},\nonumber\\
 G(x)&:=|1 - \phi(x)|g_1(x)+|\Id\phi(x)|
                    \|B(x)\langle A_0 \rangle^{-\alpha}\|_{S^d},\label{G}
\end{align}
    $A_n(x) = A_0 + (1- \phi(x)) P_n B(x) P_n$ where $P_n$ are the spectral projections as in Definition \ref{approx} and let, for $0<s\leq t$,
    \[
 X_s^{(n)}=\int_{u \in \Delta_{d-1}}\textnormal{tr}_H
 \bigwedge_{j=1}^d\bigl(\Id A_n e^{-su_jA_n^2}\bigr)\Id u,
 \qquad\mathcal X_s^{(n)}=\int_{\mathbb{H}^d}X_s^{(n)}.
\]
Then, under Hypothesis \ref{hypothesis}, $G\in L^d(\mathbb{H}^d)$, and 
there exists $C_t$, independent of $n,s$, such that
\begin{equation}\label{mainestimate}
 \int_{\mathbb{H}^d}|X_s^{(n)}|\,\Id\textnormal{vol}(x)
 \leq C_t s^{-d\alpha/2}\|G\|_{L^d}^d,
\end{equation}
where the same estimate holds with $A_{\phi,n}$ replaced by $A_\phi$.

\end{proposition}
\begin{proof}
Since $P_n$ commutes with $\langle A_0 \rangle$, we have 
\begin{equation}
\begin{split}
    &\norm{(\Id A_{\phi,n})(x) \langle A_0 \rangle^{-\alpha}}_{S^d} \\
    &= \norm{(1- \phi(x)) P_n (\Id B)(x) \langle A_0 \rangle^{-\alpha} P_n - (\Id \phi)(x) P_n B(x) \langle A_0 \rangle^{-\alpha} P_n}_{S^d} \\
    &\leq \vert 1- \phi(x) \vert \norm{(\Id B)(x) \langle A_0 \rangle^{-\alpha}}_{S^d} + \vert \Id \phi(x) \vert 
    \norm{B(x) \langle A_0 \rangle^{-\alpha}}_{S^d}.
\end{split}
\end{equation}
   For a path $\gamma$ from $x_0$ to $x$, 
 \begin{equation*}
     \begin{split}
         &\norm{B(x) \langle A_0 \rangle^{-\alpha}}_{S^d}  =
         \norm{\int_0^1 (\Id B)(\gamma(t))  \dot \gamma(t)  \langle A_0 \rangle^{-\alpha} \Id t}_{S^d} \\
         & \leq \int_0^1 g_1(\gamma(t)) \vert \dot \gamma(t) \vert \Id t,
     \end{split}
 \end{equation*}  
where the assumption 
$\tau^{\alpha,\alpha,\nabla,d,1}(B) < \infty$ gives local boundedness of $g_1$.
Taking the geodesics to reside in a fixed ball containing both 
$x_0$ and $\text{supp}(\Id \phi)$, we obtain 
$$\sup_{x \in \text{supp}(d \phi)} \norm{B(x) \langle A_0 \rangle^{-\alpha}}_{S^d} < \infty.$$ 
On the other hand, assumption \eqref{hyp3} implies 
$\norm{(\Id B)(x)}_{S^1} \leq C e^{-d(x,o)}$, where 
$d(x,o) \geq R$, so that, together with 
$\| \langle A_0 \rangle^{-\alpha}\|\leq1$ and $S^1\subset S^d$,
\begin{equation*}
    \begin{split}
        G(x)
        \leq C_1 e^{-d(x,o)} + C_2 \one_K(x),
    \end{split}
\end{equation*}
for some compact set $K \subset \mathbb{H}^d$, with the right-hand side being in $L^d(\mathbb{H}^d)$.
For a hyperbolic unit vector $v$, differentiation and commutation with $P_n$
give
\[
 \Id A_{\phi,n}(v) \langle A_0 \rangle^{-\alpha}
 =(1- \phi) P_n\Id B(v) \langle A_0 \rangle^{-\alpha}P_n
 -(\Id\phi)(v)P_nB \langle A_0 \rangle^{-\alpha}P_n,
\]
so that its $S^d$-norm is at most $G(x)$.
Using  $\|\langle A_0 \rangle^\alpha e^{-wA_{\phi,n}^2}\|_{B(H)}\leq C_{a,t}w^{-a/2}$ for $0<w \leq t$, we obtain 
\begin{equation}\label{factorbound}
 \|\Id A_{\phi,n}(v)e^{-su_jA_{\phi,n}^2}\|_{S^d}
 \leq C_{\alpha,t}G(x)s^{-\alpha/2}u_j^{-\alpha/2},
\end{equation}
when $u_j >0$.
Evaluating the exterior product on an oriented orthonormal frame at $x$ and
 applying the Schatten-von Neumann version of the H\"older estimate to each product
gives, for almost every $u$ in the simplex,
together with the fact that there are 
 $d!$ permutations, 
\begin{equation}\label{productbound}
    \left\|\bigwedge_{j=1}^d
       (\Id A_{\phi,n}e^{-su_jA_{\phi,n}^2})\right\|_{S^1}
 \leq d!C_{\alpha,t}^dG(x)^d s^{-d\alpha/2}
                  \prod_{j=1}^d u_j^{-\alpha/2}.
\end{equation}
For $\alpha<1+2/d<2$, the Dirichlet integral is
\[
 \int_{u \in\Delta_{d-1}}\prod_{j=1}^d u_j^{-\alpha/2}\dd u
 =\frac{\Gamma(1-\alpha/2)^d}{\Gamma(d(1-\alpha/2))}<\infty,
\]
so that integration in $x$ and $u$  proves
Equation \eqref{mainestimate}.
The proof for
$A_\phi$ is identical.
\end{proof}

\begin{lemma}\label{lem:limitxn}
    Under Hypothesis \ref{hypothesis} and for every $s>0$,
\[
 \|X_s^{(n)}-X_s\|_{L^1(\mathbb{H}^d;\Lambda^d T^* \mathbb{H}^d \otimes \mathbb{C})}\xrightarrow{n \rightarrow \infty}0,
 \qquad \mathcal{X}_s^{(n)}:=\int_{\mathbb{H}^d}X_s^{(n)}\xrightarrow{n \rightarrow \infty}
 \mathcal{X}_s:=\int_{\mathbb{H}^d}X_s,
\]
where $X_s$ uses $A_\phi$, with no approximation indices. 
\end{lemma}

\begin{proof}
Fix a point $x$ where Hypothesis \ref{hypothesis} holds. 
For  $ \xi \in \text{Dom} A_0$, observe that
$P_n (1- \phi(x)) B(x) P_n\xi\to (1- \phi(x)) B(x)\xi$, since
$P_n(1- \phi(x) B(x)\langle A_0 \rangle^{-1})P_n\to (1- \phi(x)) B(x)\langle A_0 \rangle^{-1}$ strongly and
$P_n\xi\to\xi$ in the graph norm of $A_0$.
Convergence on this common core implies strong resolvent convergence
$A_{\phi,n}(x)\to A_\phi(x)$ and thus convergence of the associated semigroups,
\[
 E_n(v):=e^{-vA_{\phi,n}(x)^2}\longrightarrow
 E(v):=e^{-vA_\phi(x)^2},
\]
strongly for every $v>0$.
For a fixed tangent vector $v_0$, let $T:=\Id ((1- \phi(x)) B(x))(v_0)\langle A_0 \rangle^{-\alpha}\in S^d$, so that $T_n=P_nTP_n\to T$ in $S^d$.
Set $R_n :=\langle A_0 \rangle^\alpha E_n(v)$ and $R :=\langle A_0 \rangle^\alpha E(v)$.
These operators are uniformly bounded by 
$$\sup_{n,x}\|\langle A_0 \rangle^\alpha e^{-vA_{\phi,n}(x)^2}\|_{B(H)}
 \leq C_{\alpha,t}v^{-\alpha/2},\qquad 0<v\leq t,$$
 and 
their adjoints converge strongly since on the dense set $\text{Dom} \langle A_0 \rangle^\alpha$,
\[
 R_n^*\xi=E_n(v)\langle A_0 \rangle^\alpha\xi\longrightarrow E(v)\langle A_0 \rangle^\alpha\xi=R^*\xi,
\]
with uniform boundedness extending this to all $\xi\in H$.
Therefore,
\[
 \Id A_{\phi,n}(x)(v_0)E_n(v)=T_nR_n
 \longrightarrow TR=\Id A_\phi(x)(v_0)E(v)
 \quad\hbox{in }S^d(H).
\]
Now fix $s>0$ and $u$ in the interior of $\Delta_{d-1}$.
Apply the preceding conclusion with $v=su_j$ to each
factor of the exterior product. Expanding on an
orthonormal frame
and applying the Schatten-von Neumann version of  the     H\"older estimate gives
\[
 \bigwedge_{j=1}^d
 \bigl(\mathrm dA_{\phi,n}(x)e^{-su_jA_{\phi,n}(x)^2}\bigr)
 \xrightarrow{n \rightarrow \infty}
 \bigwedge_{j=1}^d
 \bigl(\mathrm dA_\phi(x)e^{-su_jA_\phi(x)^2}\bigr), 
\]
in the trace norm of each coefficient.
By \eqref{productbound}, the trace norms of these 
products admit an $n$-independent bound
\[
 C_sG(x)^d\prod_{j=1}^d u_j^{-\alpha/2},
\]
which is integrable on $\mathbb H^d\times\Delta_{d-1}$
because $G\in L^d(\mathbb H^d)$ and $\alpha<2$.
The same bound holds for the limiting product.
Taking traces, integrating over the simplex, and
applying dominated convergence therefore gives
\[
 \int_{\mathbb H^d}|X_s^{(n)}-X_s|\,
       \mathrm d\operatorname{vol} 
 \longrightarrow0.
\]
Finally,
\[
 |\mathcal X_s^{(n)}-\mathcal X_s|
 \leq
 \int_{\mathbb H^d}|X_s^{(n)}-X_s|\,
       \mathrm d\operatorname{vol},
\]
which proves the remaining assertion.
\end{proof}

\begin{rem}\label{rem:intxn}
Let $Q(s) := s^{d/2}p_s(0)$     and observe 
that $Q$ is a polynomial and $Q'$ is bounded on $[0,t]$, whereby
$s^{d/2}Q'(s)$ is generally not a polynomial when $d$ is odd. 
We then have the estimate
\begin{equation}\label{timedom}
 |s^{d/2}Q'(s)\mathcal{X}_s^{(n)}|
 \leq C_t\|G\|_{L^d}^d s^{-d(\alpha-1)/2},
\end{equation}
where the right-hand side is integrable at zero exactly when
$\alpha<1+2/d$. Lemma \ref{lem:limitxn} and dominated convergence therefore give
\begin{equation}\label{timeconv}
 \int_0^t s^{d/2}Q'(s)\mathcal{X}_s^{(n)}\,\Id s
 \longrightarrow\int_0^t s^{d/2}Q'(s)\mathcal{X}_s\,\Id s.
\end{equation}
\end{rem}

\begin{theorem}\label{maintheorem}
    Let $t >0$ and assume Hypothesis \ref{hypothesis}.
For real-valued $\phi \in C^{\infty}_c(\mathbb{H}^d)$, 
write $A_{\phi} := A_0 + (1- \phi)(A-A_0)$.
Then 
\begin{equation}\label{eq:traceformula}
    \begin{split}
&\textnormal{tr}_{L^2(\mathbb{H}^d,H)}
    \textnormal{tr}_{S} \left(
    e^{-tD^*D} - e^{-tDD^*}\right) \\
    &=  \frac{2}{d} i^d \kappa_c 
     \left[ t^d p_t(0) \mathcal{X}_t -
    \int_0^t s^{d/2} \partial_s \left(s^{d/2} p_s(0)\right) \mathcal{X}_s \Id s  
     \right],
    \end{split}
\end{equation}
 where 
 \begin{equation*}
 \mathcal{X}_s = \int_{\mathbb{H}^d} X_s, \qquad      X_s := \int_{u \in \Delta_{d-1}} 
    \textnormal{tr}_H \bigwedge_{j=1}^d 
  \bigg( (\Id A_{\phi})e^{-su_j A_{\phi}^2}\bigg) \Id u,
 \end{equation*}
and 
\begin{equation*}
        \begin{split}
            p_t(0)
=
&\frac{1}{2^{d} \pi^{d/2}\Gamma(d/2)} \\
& \cdot\sum_{k=0}^{m}
\Gamma\left(k+\frac12\right)
e_{m-k}\left(
\frac14,\frac94,\ldots,\frac{(2m-1)^2}{4}
\right)
t^{-k-\frac12},
        \end{split}
    \end{equation*}
    with $e_j$ denoting the $j$-th elementary symmetric polynomial,
with $e_0=1$.

\end{theorem}
\begin{proof}
Let us write $A_* := A_{n,m,l}$ and 
$D_* := D_{n,m,l}$.
    We start with Equation \eqref{prelformula} from Lemma \ref{highervanish}, and develop the 
    trace using Proposition \ref{cliffordid}
    as to obtain
    \begin{equation}
W_t := \textnormal{tr}_{L^2(\mathbb{H}^d,H)}
    \textnormal{tr}_{S} \left(
    e^{-tD_*^*D_*} - e^{-tD_*D_*^*}\right) = \sum_{k=0}^m I_{m,k},         \end{equation}
where 
\begin{equation}\label{iintegral}
    \begin{split}
        &I_{m,k} :=      \textnormal{tr}_H \bigg(
  \lim_{r \rightarrow \infty}  \sinh^{d-1}r 
  \int_0^t b_{m,k} s^{-k-\frac{1}{2}} i^d \kappa_c 
  \int_{S_1(o)}  \\
  &A_*(\exp_o(ry)) 
  \int_{u \in s \Delta_{d-1}} 
  \sum_{j=1}^d 
  \sum_{\alpha \in \{1, \dots,d\}^{d-1}} \\
& \frac{y^j}{\cosh(r)} \epsilon_{(j,\alpha)} e^{-uA_*^2(\exp_o(ry))}    \shuffle
\bigg((\Id_{\alpha_1}A_*)(\exp_o(ry)),\dots,
(\Id_{\alpha_{d-1}}A_*)(\exp_o(ry)) \bigg) \\
  &\Id u \Id S_1(y) \Id s \bigg).
    \end{split}
\end{equation}
Now, if $x = (\sinh(r) y, \cosh(r)) \in \mathbb{H}^d$ with $y \in S^{d-1}$, then
$\Id r = \frac{ \sum_{j=1}^d y^j \Id x^j}{
 \cosh(r)
}$. 
For $u \in s\Delta_{d-1}$, we therefore have
\begin{equation}
\begin{split}
   & \Id r \wedge \bigg( 
    e^{-suA_*^2} \shuffle_{\wedge} (\Id A_*)^{(d-1)}
    \bigg) = \sum_{j=1}^d \sum_{\alpha \in \{1, \dots,d\}^{d-1}} \\&
  \frac{ y^j}{\cosh(r)} \epsilon_{(j,\alpha)}e^{-uA_*^2}    \shuffle
\bigg((\Id_{\alpha_1}A_*),\dots,
(\Id_{\alpha_{d-1}}A_*) \bigg) 
\Id x.
\end{split}
\end{equation}
Hence, after applying the interior multiplication  and restricting to $S_R(o)$, we obtain
    \begin{equation}\label{stoked}
    \begin{split}
 &e^{-suA_*^2} \shuffle_{\wedge} (\Id A_*)^{(d-1)}
 = \sum_{j=1}^d \sum_{\alpha \in \{1, \dots,d\}^{d-1}} \\
&y^j \epsilon_{(j,\alpha)}e^{-uA_*^2}    \shuffle
\bigg((\Id_{\alpha_1}A_*),\dots,
(\Id_{\alpha_{d-1}}A_*) \bigg) 
\Id S_R(o),
        \end{split}
    \end{equation}
so that, upon inserting this into Equation
\eqref{iintegral}, changing variables 
in the $u$-integral and summing in $k$, yields
\begin{equation}\label{wequation}
\begin{split}
     W_t =  i^d \kappa_c \text{tr}_H
    \bigg( &\lim_{R \rightarrow \infty} \int_0^t 
    s^{d-1} p_s(0) \\
    &\int_{S_R(o)} A_*
    \int_{u \in \Delta_{d-1}} e^{-suA_*^2} \shuffle_{\wedge} (\Id A_*)^{(d-1)}
    \Id u \Id s \bigg).
\end{split}
\end{equation}
Given that $R \rightarrow \infty$, we can 
replace every occurrence of $A_*$ by 
$A_{ \phi,*} :=  A_0 + (1- \phi)(A_*-A_0)$, where
$\phi \in C_c^{\infty}(\mathbb{H}^d)$.
Note moreover that, since $\Id A_*$ is finite rank,
the integrand is Bochner integrable in trace-class.
Commuting the trace with the integral and applying 
Stoke's Theorem gives us then
\begin{equation}\label{weq}
    \begin{split}
        W_t = i^d \kappa_c 
        \int_0^t s^{d-1} p_s(0)
        \big( \mathcal{X}^{(n,m,l)}_s
        + \mathcal{Y}^{(n,m,l)}_s \big) \Id s,
    \end{split}
\end{equation}
where $\mathcal{X}^{(n,m,l)}_s = \int_{\mathbb{H}^d} X^{(n,m,l)}_s$, $\mathcal{Y}^{(n,m,l)}_s = \int_{\mathbb{H}^d} Y^{(n,m,l)}_s$,
\begin{equation}
    X^{(n,m,l)}_s := \int_{u \in \Delta_{d-1}} 
    \text{tr}_H \bigwedge_{j=1}^d 
  \bigg( (\Id A_{\phi,*})e^{-su_j A_{\phi,*}^2}\bigg) \Id u,
\end{equation}
and 
\begin{equation}
    Y^{(n,m,l)}_s := -2s \int_{u \in \Delta_{d}} 
    \text{tr}_H \bigg( A_{\phi,*}^2 e^{-su_0A_{\phi,*}^2} \bigwedge_{j=1}^d 
  \bigg( (\Id A_{\phi,*})e^{-su_j A_{\phi,*}^2}\bigg) \bigg) \Id u.
\end{equation}
Rewriting $Y^{(n,m,l)}_s$ as 
\begin{equation}
     Y^{(n,m,l)}_s = -2s \int_{u \in \Delta_{d-1}} 
  u_0  \text{tr}_H \bigg( (\Id A_{\phi,*}) e^{-su_0A^2} A_{\phi,*}^2 \bigwedge_{j=1}^{d-1} 
  \bigg( (\Id A_{\phi,*})e^{-su_j A_{\phi,*}^2}\bigg) \bigg) \Id u,
\end{equation}
and changing variables, we notice that
\begin{equation}
    Y^{(n,m,l)}_s = \frac{2}{d} s \partial_s X^{(n,m,l)}_s.
\end{equation}
Since we may commute differentiation with integration, we also have 
$\mathcal{Y}^{(n,m,l)}_s = \frac{2}{d} s \partial_s \mathcal{X}^{(n,m,l)}_s$.
Upon inserting this into Equation \eqref{weq}, we are left with
\begin{equation}
    W_t =   i^d \kappa_c
 \int_0^t s^{d-1} p_s(0)
        \big( \mathcal{X}^{(n,m,l)}_s
        + \frac{2}{d} s \partial_s \mathcal{X}^{(n,m,l)}_s \big) \Id s.
\end{equation}
Now, for each fixed $n,m,l \in \mathbb{N}$, we have 
$\sup_{0 < s \leq t} \vert \mathcal{X}_s^{(n,m,l)} \vert < \infty$, by the calculations performed in Proposition 
\ref{prop:xnbounded}.
On the other hand,
$$\mathcal{X}^{(n,m,l)}_s
        + \frac{2}{d} s \partial_s \mathcal{X}^{(n,m,l)}_s = \frac{2}{d} 
        s^{1 - \frac{d}{2}} \partial_s(s^{\frac{d}{2}}\mathcal{X}^{(n,m,l)}_s), $$
and hence 
\begin{equation}
     W_t =  \frac{2}{d} i^d \kappa_c
 \int_0^t s^{d/2} p_s(0)
      \partial_s(s^{d/2}\mathcal{X}^{(n,m,l)}_s) \Id s.
\end{equation}
Integrating by parts on $[\varepsilon,t]$, 
we obtain, together with
\[
\epsilon^dp_\epsilon(0)\mathcal{X}_\epsilon^{(n,m,l)}
 =\epsilon^{d/2}Q(\epsilon)\mathcal{X}_\epsilon^{(n,m,l)}\longrightarrow0,
\]
as $\epsilon \rightarrow 0$, where 
$Q(s):=s^{d/2}p_s(0)$ as in Remark \ref{rem:intxn},
\begin{equation}\label{eq66}
     W_t =  \frac{2}{d} i^d \kappa_c 
     \left[ t^d p_t(0) \mathcal{X}^{(n,m,l)}_t -
    \int_0^t s^{d/2} \partial_s \left(s^{d/2} p_s(0)\right) \mathcal{X}^{(n,m,l)}_s \Id s  
     \right].
\end{equation}

It remains to argue that we can apply the 
$n,m,l$-limit to the trace formula
in Equation \eqref{eq66}. 
For that, set $r(x):=d(x,o)$. 
Furthermore, 
 retain the approximation indices on the left-hand
side by writing $W_t^{(n,m,l)}$. \\
We start by establishing  estimates for fixed $n$.
Since $P_n$ commutes with $\langle A_0\rangle$,
\[
 \|\mathrm dB_n(x)\|_{S^d}
 \leq \langle n\rangle^\alpha
       \|\mathrm dB(x)\langle A_0\rangle^{-\alpha}
        \|_{S^d}
                                      \leq C_n e^{-r(x)}.
\]
Moreover, $\sup_{x\in\mathbb H^d}\|B_n(x)\|_{S^d}<\infty$,
because $\int_0^\infty e^{-r}\,\mathrm dr<\infty$, and hence
\[
 \sup_{m,x}\|B_{n,m}(x)\|_{S^d}\leq C_n,
 \qquad
 \sup_m\|\mathrm dB_{n,m}(x)\|_{S^d,g}
 \leq C_n e^{-r(x)}.
\]
By positivity and the normalization of the mollifiers in Definition \ref{delta}, we have
\[
 \|B_{n,m,l}(x)\|_{S^d}
 \leq\sup_y\|B_{n,m}(y)\|_{S^d}\leq C_n,
\]
where $C_n$ is independent on $m$ and $l$.
Furthermore,
\[
 \int_{\mathbb H^d}\mathrm d_x k_l(x,y)\,
                  \mathrm d\operatorname{vol}(y)=0,
 \qquad
 \int_{\mathbb H^d}|\mathrm d_x k_l(x,y)|\,
                  \mathrm d\operatorname{vol}(y)
 \leq Cl, 
\]
hence
\begin{equation}
     \mathrm dB_{n,m,l}(x)
 =
 \int_{\mathbb H^d}\mathrm d_x k_l(x,y)
       \bigl(B_{n,m}(y)-B_{n,m}(x)\bigr)\,
                  \mathrm d\operatorname{vol}(y).
\end{equation}
On the support of $k_l(x,\cdot)$, we have
$d(x,y)\leq l^{-1}$, and the preceding derivative
estimate gives
\[
 \|B_{n,m}(y)-B_{n,m}(x)\|_{S^d}
 \leq C_n l^{-1}e^{-r(x)}, 
\]
so that $\sup_{m,l}\|\mathrm dB_{n,m,l}(x)\|_{S^d}
 \leq C_n e^{-r(x)}.$
Since $\mathrm d\phi$ is compactly supported,
we can therefore choose $ G_n(x):=C_n e^{-r(x)}\in L^d(\mathbb H^d)$
in order to obtain
\[
 \|\mathrm dA_{\phi,n,m,l}(x)\|_{S^d}
 \leq G_n(x),
\]
for every $m,l$.
After increasing $C_n$ if necessary, the same bound
holds for $A_{\phi,n,m}$ and $A_{\phi,n}$.
For fixed $n,m$, the  properties of
$\delta_l$ give, locally in $x$,
\[
 B_{n,m,l}(x)\xrightarrow{l \rightarrow \infty} B_{n,m}(x),
 \qquad
 \mathrm dB_{n,m,l}(x)\xrightarrow{l \rightarrow \infty}
                         \mathrm dB_{n,m}(x),
\]
in $S^d(H)$.
For fixed $n$, $K_m\to I$ similarly
gives
\[
 B_{n,m}(x)\longrightarrow B_n(x),
 \qquad
 \mathrm dB_{n,m}(x)\longrightarrow\mathrm dB_n(x),
\]
in $S^d(H)$, where we used $B_n(x)\in S^d(H)$, as established above.
Since these convergences also hold in operator norm and  the perturbations of $A_0$
are bounded and symmetric, this implies norm
resolvent convergence of the corresponding
self-adjoint operators. 
Hence, for every 
$v>0$,
\[
 e^{-vA_{\phi,n,m,l}(x)^2}
 \longrightarrow e^{-vA_{\phi,n,m}(x)^2},
 \qquad l\to\infty,
\]
and
\[
 e^{-vA_{\phi,n,m}(x)^2}
 \longrightarrow e^{-vA_{\phi,n}(x)^2},
 \qquad m\to\infty,
\]
in operator norm, by continuous functional calculus.
It follows that, for every $s>0$ and every $u$
in the interior of $\Delta_{d-1}$, each factor $(\mathrm dA_{\phi,n,m,l}(x))(w)\,
             e^{-su_jA_{\phi,n,m,l}(x)^2}$
converges in $S^d(H)$ in the successive $l$- and
$m$-limits, where $w \in T_x \mathbb{H}^d$.
Applying the 
Schatten-von Neumann H\"older estimate also gives the uniform bound
\[
 \left\|
 \bigwedge_{j=1}^d
 \left(
 \mathrm dA_{\phi,n,m,l}(x)
 e^{-su_jA_{\phi,n,m,l}(x)^2}
 \right)
 \right\|_{S^1}
 \leq d!\,G_n(x)^d, 
\]
where the  same estimate holds after taking either limit.
Its right-hand side is integrable on
$\mathbb H^d\times\Delta_{d-1}$ and is independent
of $s,m,l$.
Dominated convergence thus yields,
for every $s>0$,
\[
 \|X_s^{(n,m,l)}-X_s^{(n,m)}\|_{L^1(\mathbb{H}^d;\Lambda^d T^* \mathbb{H}^d \otimes \mathbb{C})}
 \longrightarrow0,
 \qquad l\to\infty,
\]
and
\[
 \|X_s^{(n,m)}-X_s^{(n)}\|_{L^1(\mathbb{H}^d;\Lambda^d T^* \mathbb{H}^d \otimes \mathbb{C})}
 \longrightarrow0,
 \qquad m\to\infty.
\]
The same bound then also justifies passage through the
time integral in \eqref{eq66}.
Indeed, because $Q'$ is bounded on $[0,t]$,
\[
 \left|
 s^{d/2}Q'(s)\mathcal X_s^{(n,m,l)}
 \right|
 \leq C_{n,t}s^{d/2},
 \qquad 0<s\leq t, 
\]
where the bound is integrable and independent of
$m,l$. 
Hence,
\[
 \lim_{m\to\infty}\lim_{l\to\infty}
 \int_0^t s^{d/2}Q'(s)
                  \mathcal X_s^{(n,m,l)}\,\mathrm ds
 =
 \int_0^t s^{d/2}Q'(s)
                  \mathcal X_s^{(n)}\,\mathrm ds,
\]
and
\[
 \lim_{m\to\infty}\lim_{l\to\infty}
 \mathcal X_t^{(n,m,l)}
 =\mathcal X_t^{(n)}.
\]
The final  limit $n\to\infty$ is then 
handled by Proposition \ref{prop:xnbounded}, Lemma \ref{lem:limitxn} and Remark \ref{rem:intxn}.
Therefore, Corollary \ref{cor:semigrouplimit} and continuity of the trace
in the trace norm imply
\[
 \lim_{n\to\infty}\lim_{m\to\infty}\lim_{l\to\infty}
 W_t^{(n,m,l)}
 =
 \operatorname{tr}_{L^2(\mathbb H^d,H)}
 \operatorname{tr}_{S}
 \left(e^{-tD^*D}-e^{-tDD^*}\right),
\]
thus proving \eqref{eq:traceformula}.

\end{proof}

\begin{corollary}
     Let $t >0$ and assume Hypothesis \ref{hypothesis}.
For $\phi \in C^{\infty}_c(\mathbb{H}^d)$, 
write $A_{\phi} := A_0 + (1- \phi)(A-A_0)$.
Then 
\begin{equation}
    \begin{split}
&\textnormal{tr}_{L^2(\mathbb{H}^d,H)}
    \textnormal{tr}_{S} \left(
    e^{-tD^*D} - e^{-tDD^*}\right) \\
    &=
    i^d \kappa_c  \cdot 
     \lim_{R \rightarrow \infty}
  \textnormal{tr}_H \bigg[
  \int_{S_R(o)}
  \int_0^t s^{d-1} p_s(0) \\
  &\int_{u \in \Delta_{d-1}} 
  A_{\phi} e^{-su_0 A_{\phi}^2}
     \bigwedge_{j=1}^{d-1}
  \bigg( (\Id A_{\phi})e^{-su_j A_{\phi}^2}\bigg) \Id u \Id s \bigg],
    \end{split}
\end{equation}
exists as a Pettis integral.
\end{corollary}

\begin{rem}
    In the Euclidean case, the analogous of 
    Equation \eqref{eq:traceformula}
    reads as 
    \begin{equation}\label{eq:euclideancase}
        \begin{split}
   &\textnormal{tr}_{L^2(\mathbb{R}^d,H)}
    \textnormal{tr}_{\mathbb{C}^r} \left(
    e^{-tD^*D} - e^{-tDD^*}\right)  \\
    &= \frac{2}{d} (4 \pi)^{-d/2}i^d \kappa_c
    t^{d/2} \int_{s \in \Delta_{d-1}} 
    \int_{\mathbb{R}^d} 
    \textnormal{tr}_H \bigwedge_{j=0}^{d-1}
  \bigg( (\Id A_{\phi})e^{-ts_j A_{\phi}^2}\bigg) \Id s,
        \end{split}
    \end{equation}
    where $r = 2^{\frac{d-1}{2}}$. Now,
    when $t \downarrow 0$, then, by Remark 
    \ref{rem:smalltimespinorkernel},
    $p_t(0)
\sim
(4 \pi)^{-d/2}
t^{-d/2}$, hence $s^{d/2}p_s(0) \sim a_d$.
Moreover, by Remark \ref{rem:intxn}, 
\[
 \left|\int_0^t s^{d/2}Q'(s)\mathcal{X}_s\,\Id s\right|
 \leq C t^{1-d(\alpha-1)/2}\longrightarrow0.
\]
so that the integral in Equation \eqref{eq:traceformula}
vanishes and what remains is $\mathcal{X}_t$,
multiplied by $\frac{d}{2} (4 \pi)^{-d/2}i^d \kappa_c
    t^{d/2}$, in agreement with  \eqref{eq:euclideancase}.
\end{rem}

\subsection*{Acknowledgements}
I thank Oliver Fürst and Matthias Lesch for supervising the master's thesis out of which this paper originated, and Edward McDonald for telling me the story surrounding  the mean motion Theorem. \\
I gratefully acknowledge the financial support of the
Bonn International Graduate School for Mathematics (BIGS Mathematics) as part of  the Hausdorff Center for Mathematics (HCM),
funded by the Deutsche Forschungsgemeinschaft (DFG, German Research Foundation) under
Germany's Excellence Strategy - GZ 2047/1, Projekt-ID 390685813.

\subsection*{AI disclosure}
ChatGPT 5.6 Sol (developed by OpenAI) was used to screen the literature to obtain an overview of the Callias problem and the use of the Witten index through its historical applications in physics,
 and for numerically checking the calculations involving the diagonal of the hyperbolic spinor heat kernel  (Proposition 5.20) and the bounds 
 in Lemma 6.3.
Furthermore, the  strategy of passing to Fermi coordinates in the proof of Lemma 6.5 was suggested by ChatGPT 5.6 Sol. 
Besides minor stylistic improvements regarding the LaTeX typesetting and notations, as well as proofreading, no AI was used in the writing of the paper.
I have written and verified the manuscript, and take full responsibility
for its content.


\begin{thebibliography}{9}

\bibitem[]{anghel1} {\sc N. Anghel}, \textit{$L^2$-index formulae for perturbed Dirac operators}, Commun. Math. Phys. 128, no.
1, 77-97, (1990).

\bibitem[]{anghel2} {\sc N. Anghel}, \textit{On the index of Callias-type operators}, Geom. Funct. Anal. 3, no. 5, 431-438,
(1993).


\bibitem[]{atiyah1} {\sc M. Atiyah, V. Patodi, I.} and  {\sc Singer}, \textit{ Spectral asymmetry and Riemannian geometry}. I., Math.
Proc. Cambridge Philos. Soc. 77, 43-69, (1975).

\bibitem[]{atiyah2} {\sc M. Atiyah, V. Patodi, I.} and  {\sc Singer}, 
\textit{Spectral asymmetry and Riemannian geometry. II.}, Math.
Proc. Cambridge Philos. Soc. 78, 405-432, (1975).

\bibitem[]{atiyah3} {\sc M. Atiyah, V. Patodi, I.} and  {\sc Singer}, 
\textit{Spectral asymmetry and Riemannian geometry. III.}, Math.
Proc. Cambridge Philos. Soc. 79, 71-99, (1976).

\bibitem[]{lectures} {\sc R. Benedetti } and {\sc Petronio}, \textit{Lectures on hyperbolic geometry}. Springer Science $\&$ Business Media, (1992).

\bibitem[]{bolle} {\sc D. Bollé, F. Gesztesy, H. Grosse, W. Schweiger} and {\sc B. Simon}, \textit{Witten Index, Axial Anomaly, and
Krein’s Spectral Shift Function in Supersymmetric Quantum Mechanics}, J. Math. Phys. 28,
no. 7, p.1512-1525, (1987).




\bibitem[]{braverman} {\sc M. Braverman} and {\sc P. Shi}, \textit{Cobordism invariance of the index of Callias-type operators}, Comm.
Partial Diﬀerential Equations 41, no. 8, 1183-1203, (2016).

\bibitem[]{bruning} {\sc J. Brüning} and { \sc H. Moscovici}, \textit{$L^2$-index for certain Dirac-Schrödinger operators}, Duke Math. J.
66, no. 2, 311-336 (1992).

\bibitem[]{radon} {\sc T. P. Branson; G. Ólafsson} and {\sc Schlichtkrull}, 
\textit{A bundle valued Radon transform, with applications to invariant wave equations}. The Quarterly Journal of Mathematics, 45. Jg., Nr. 4, S. 429-461, (1994).

\bibitem[]{bunke91} {\sc U. Bunke}, 
\textit{The spectrum of the Dirac operator on the hyperbolic space}, Mathematische Nachrichten, 153. Jg., Nr. 1, S. 179-190, (1991).

\bibitem[]{bunke95} {\sc U. Bunke}, \textit{A K-theoretic relative index theorem and Callias-type Dirac operators}, Math. Ann.
303, no. 1, 241-279, (1995).

\bibitem[]{callias} {\sc C. Callias}, \textit{Index Theorems on Open Spaces}. Communications in Mathematical Physics 62, 213–234, (1978).

\bibitem[]{camporesi92} {\sc R. Camporesi}, \textit{The spinor heat kernel in maximally symmetric spaces}. Communications in mathematical physics 148.2, 283-308, (1992).



\bibitem[]{camporesi} {\sc R. Camporesi} and {\sc Pedon}, \textit{Harmonic analysis for spinors on real hyperbolic spaces}. Colloquium Mathematicum. Instytut Matematyczny Polskiej Akademii Nauk, 245-286, (2001).

\bibitem[]{carey77} {\sc R. W. Carey} and {\sc J. D. Pincus}, 
\textit{Mosaics, principal functions, and mean motion in von Neumann algebras}.
Acta Math., 138(3--4), 153--218, (1977). 

\bibitem[]{carey79} {\sc R. W. Carey} and {\sc J. D. Pincus}, 
\textit{Mean motion, principal functions, and the zeros of Dirichlet series}.
Integral Equations and Operator Theory, 2(4), 484--502, (1979).

\bibitem[]{carey14} {\sc A. Carey, F. Gesztesy, D. Potapov, F. Sukochev} and  {\sc Y. Tomilov}, \textit{A Jost–Pais-type reduction of Fredholm determinants and some applications}. Integral Equations and Operator Theory, 79(3), 389-447, (2014).


\bibitem[]{carey16} {\sc A. Carey, F. Gesztesy, H. Grosse, G. Levitina, D. Potapov, F. Sukochev} and {\sc D. Zanin}, \textit{Trace
Formulas for a Class of non-Fredholm Operators: A Review}, Reviews in Mathematical Physics
28, No. 10, (2016).

\bibitem[]{carey162} {\sc A. Carey, F. Gesztesy, G. Levitina, D. Potapov, F. Sukochev} and  {\sc D. Zanin},  \textit{On index theory for non‐Fredholm operators: A (1+ 1)‐dimensional example}. Mathematische Nachrichten, 289(5-6), 575-609, (2016).

\bibitem[]{carey2022} {\sc A. Carey}, and {\sc G. Levitina}. \textit{Index theory beyond the Fredholm case}. Vol. 2323. Cham: Springer, (2022).

\bibitem[]{carvalho} {\sc C. Carvalho} and  {\sc V. Nistor}, \textit{An index formula for perturbed Dirac operators on Lie manifolds}, J.
Geom. Anal. 24, 1808-1843, (2014).

\bibitem[]{cecchini} {\sc S. Cecchini}, \textit{Callias-type operators in $C^*$-algebras and positive scalar curvature on noncompact
manifolds}, J. Topol. Anal. 12, no. 4, 897-939, (2020).

\bibitem[]{coburn} {\sc L. A. Coburn, R. G. Douglas, D. G. Schaeffer}, and {\sc I. M. Singer}, 
\textit{$C^*$-algebras of operators on a half-space. II. Index theory}.
Publications Mathématiques de l'IHÉS, (40), 69--79, (1971).



\bibitem[]{dungen} {\sc K. van den Dungen}, \textit{The index of generalised Dirac-Schrödinger operators}, J. Spectr. Theory
9, 1459-1506, (2019).

\bibitem[]{dungen2} {\sc K. van den Dungen}, \textit{Generalised Dirac-Schrödinger operators and the Callias Theorem},
arXiv:2312.17600, (2023).

\bibitem[]{oliver23} { \sc O. Fürst}, \textit{Trace and Index of Dirac-Schrödinger Operators on Open Space with Operator Potentials}. arXiv preprint, arXiv:2311.02593, (2023).

\bibitem[]{oliver25} { \sc O. Fürst}, 
\textit{The Witten index of massless (d+1)-Dirac-Schrödinger operators}. Journal of Geometry and Physics,  209. Jg., S. 105413, (2025).

\bibitem[]{olivershift} {\sc O. Fürst},  \textit{Higher order spectral shift of Euclidean Callias operators}. arXiv preprint arXiv:2506.01647 (2025).

\bibitem[]{mercermfd} {\sc O. Fürst} 
and {\sc J. Kevo},
\textit{A Mercer-type Theorem for Trace-Class Integral Operators on Direct Integral Hilbert Spaces}, in preparation. 


\bibitem[]{gesztesy11} {\sc F. Gesztesy, Y. Latushkin, K. Makarov, F. Sukochev} and {\sc Y. Tomilov}, \textit{The Index Formula and
the Spectral Shift Function for Relatively Trace Class Perturbations}. Adv.Math. 227, 319-420,
(2011).

\bibitem[]{gesztesy88} {\sc F. Gesztesy} and {\sc B. Simon}, \textit{Topological Invariance of the Witten Index}, Journal of Functional
Analysis 79, p.91-102, (1988)


\bibitem[]{gesztesy16} {\sc F. Gesztesy} and  {\sc M. Waurick}, \textit{The Callias index formula revisited}, Lecture Notes in Mathematics,
Springer International Publishing, (2016).

\bibitem[]{grig} {\sc A. Grigor'yan} and {\sc M. Noguchi}, \textit{ The heat kernel on hyperbolic space}. Bulletin of the London Mathematical Society, 30(6), 643-650, (1998).

\bibitem[]{batu} {\sc B. Güneysu}, 
\textit{Covariant Schrödinger semigroups on Riemannian manifolds}. Springer International Publishing, (2017).

\bibitem[]{helgason2} {\sc S. Helgason}, \textit{Groups and geometric analysis: integral geometry, invariant differential operators, and spherical functions}. American Mathematical Society, (2022).

\bibitem[]{helgason} {\sc S. Helgason}, \textit{Geometric Analysis on Symmetric Spaces}. American Mathematical Society, (2024).

\bibitem[]{hooft} {\sc G. 't Hooft}, \textit{Magnetic monopoles in unified theories}. Nucl. Phys. B, 79 (CERN-TH-1876), 276-284, (1974).

\bibitem[]{jackiw} {\sc R. Jackiw} and {\sc C. Rebbi}, \textit{Solitons with Fermion Number $1/2$}, Phys. Rev. D 13, 3398, (1976).

\bibitem[]{knapp86} {\sc A. W. Knapp}, \textit{Representation theory of semisimple groups: An Overview Based on Examples}. 
Princeton Math. Ser. 36,  Princeton Univ. Press, (1986).

\bibitem[]{kottke} {\sc C. Kottke}, \textit{A Callias-type index theorem with degenerate potentials}, Comm. Partial Diﬀerential
Equations 40, no. 2, 219-264, (2015).

\bibitem[]{kucero} {\sc D. Kucerovsky}, \textit{A short proof of an index theorem}, Proc. Amer. Math. Soc. 129, no. 12,
3729-3736, (2001).

\bibitem[]{lichno} {\sc A. Lichnerowicz},  \textit{Géometrie différentielle-spineurs harmoniques}. Comptes Rendus Hebdomadaires Des Séances De L'Académie Des Sciences 257.1 : 7,
(1963).

\bibitem[]{mckean} {\sc H. P. McKean, Jr.} and {\sc I. M. Singer}, 
\textit{Curvature and the eigenvalues of the Laplacian}.
Journal of Differential Geometry, 1(1), 43--69, (1967).


\bibitem[]{partha} {\sc R. Parthasarathy}, \textit{Dirac operator and the discrete series}. Annals of mathematics,  96. Jg., Nr. 1, S. 1-30, (1972).


\bibitem[]{peller} {\sc V. V. Peller}, \textit{Multiple operator integrals and higher operator derivatives}. Journal of Functional Analysis, 233(2), 515-544, (2006).


\bibitem[]{polyakov} {\sc A. M. Polyakov},  \textit{Particle spectrum in quantum field theory}. JETP lett, 20(194-195), 300, (1974).

\bibitem[]{pushni} {\sc A. Pushnitski}, \textit{The Spectral Flow, the Fredholm Index, and the Spectral Shift Function}, T.
Suslina, D. Yafaev (Eds.), Spectral Theory of Diﬀerential Operators: M.Sh. Birman 80th
Anniversary Collection, Amer. Math. Soc. Transl. Ser. 2, Adv. Math. Sci., vol. 225, Amer.
Math. Soc., Providence, RI,  p. 141-155, (2008).

\bibitem[]{rade} {\sc J. Rade}, \textit{Callias’ index theorem, elliptic boundary conditions, and cutting and gluing}, Com-
mun. Math. Phys. 161, no. 1, 51-61, (1994).

\bibitem[]{schulz}
{\sc H. Schulz-Baldes} and {\sc T. Stoiber}, \textit{Callias-type operators associated to spectral triples}, J. Noncom-
mut. Geom. 17, 527-572, (2023).

\bibitem[]{shi} {\sc P. Shi}, \textit{The index of Callias-type operators with Atiyah-Patodi-Singer boundary conditions},
Ann. Glob. Anal. Geom. 52, 465-482, (2017).

\bibitem[]{simon} {\sc B. Simon}, \textit{Trace ideals and their applications}. No. 120. American Mathematical Soc., (2005).

\bibitem[]{wallach} {\sc N. R. Wallach}, 
\textit{Harmonic analysis on homogeneous spaces}. Courier Dover Publications, (2018).

\bibitem[]{weinberg} {\sc E. J.  Weinberg}, 
\textit{Fundamental monopoles and multimonopole solutions for arbitrary simple gauge groups}. Nuclear Physics B, 167(3), 500-524, 
(1980).


\bibitem[]{witten81} {\sc E. Witten},  \textit{Dynamical breaking of supersymmetry}. Nuclear Physics B, 188(3), 513-554, (1981).



     



\end{thebibliography}
\end{document}